\documentclass[12pt,leqno]{article}
\usepackage{amssymb}
\usepackage[mathscr]{eucal}
\usepackage{amsmath,amssymb,latexsym,theorem,bbm}
\usepackage{color}
\usepackage{appendix}

\newcommand{\SC}{\scriptstyle}

\newcommand{\CC}{\mathsf{C}}
\newcommand{\DD}{\mathsf{D}}

\newcommand{\NN}{\mathbb{N}}

\newcommand{\RR}{\mathbb{R}}

\newcommand{\ZZ}{\mathbb{Z}}

\newcommand{\bA}{{\boldsymbol{A}}}

\newcommand{\tbA}{\widetilde{\bA}}
\newcommand{\bb}{{\boldsymbol{b}}}

\newcommand{\tc}{\widetilde{c}}
\newcommand{\bd}{{\boldsymbol{d}}}
\newcommand{\tbd}{\widetilde{\bd}}
\newcommand{\be}{{\boldsymbol{e}}}
\newcommand{\bF}{{\boldsymbol{F}}}
\newcommand{\bg}{{\boldsymbol{g}}}
\newcommand{\bI}{{\boldsymbol{I}}}
\newcommand{\tK}{\widetilde{K}}

\newcommand{\tQ}{\widetilde{Q}}
\newcommand{\bu}{{\boldsymbol{u}}}
\newcommand{\tbu}{\widetilde{\bu}}
\newcommand{\bv}{{\boldsymbol{v}}}
\newcommand{\tbv}{\widetilde{\bv}}
\newcommand{\bx}{{\boldsymbol{x}}}
\newcommand{\bX}{{\boldsymbol{X}}}
\newcommand{\by}{{\boldsymbol{y}}}
\newcommand{\bY}{{\boldsymbol{Y}}}

\newcommand{\bz}{{\boldsymbol{z}}}
\newcommand{\bZ}{{\boldsymbol{Z}}}
\newcommand{\bU}{{\boldsymbol{U}}}
\newcommand{\bgamma}{{\boldsymbol{\gamma}}}
\newcommand{\bmu}{{\boldsymbol{\mu}}}

\newcommand{\bzero}{{\boldsymbol{0}}}
\newcommand{\bone}{{\boldsymbol{1}}}

\newcommand{\cA}{{\mathcal A}}
\newcommand{\cB}{{\mathcal B}}
\newcommand{\tcB}{\widetilde{\cB}}

\newcommand{\cF}{{\mathcal F}}

\newcommand{\cL}{{\mathcal L}}
\newcommand{\cM}{{\mathcal M}}

\newcommand{\cN}{{\mathcal N}}
\newcommand{\cP}{{\mathcal P}}
\newcommand{\cS}{{\mathcal S}}
\newcommand{\cU}{{\mathcal U}}
\newcommand{\cT}{{\mathcal T}}
\newcommand{\bcU}{\boldsymbol{\cU}}

\newcommand{\cX}{{\mathcal X}}
\newcommand{\cY}{{\mathcal Y}}
\newcommand{\cZ}{{\mathcal Z}}
\newcommand{\cW}{{\mathcal W}}
\newcommand{\bcW}{\boldsymbol{\cW}}
\newcommand{\bcY}{\boldsymbol{\cY}}

\newcommand{\bcZ}{\boldsymbol{\cZ}}
\newcommand{\tcW}{\widetilde{\cW}}

\newcommand{\dd}{\mathrm{d}}

\newcommand{\slu}{{\SC\mathrm{lu}}}

\newcommand{\ARtwo}{\textup{AR(2)}}
\newcommand{\INARtwo}{\textup{INAR(2)}}
\newcommand{\INARp}{\textup{INAR($p$)}}

\newcommand{\EE}{\operatorname{\mathbb{E}}}
\newcommand{\PP}{\operatorname{\mathbb{P}}}
\newcommand{\OO}{\operatorname{O}}

\newcommand{\halpha}{\widehat{\alpha}}
\newcommand{\hbeta}{\widehat{\beta}}
\newcommand{\hvarrho}{\widehat{\varrho}}

\newcommand{\hmu}{\widehat{\mu}}

\newcommand{\tV}{\widetilde{V}}

\newcommand{\tPhi}{\widetilde{\Phi}}

\newcommand{\vare}{\varepsilon}

\renewcommand{\mid}{\,|\,}
\newcommand{\bmid}{\,\big|\,}

\renewcommand{\leq}{\leqslant}
\renewcommand{\geq}{\geqslant}

\newcommand{\stoch}{\stackrel{\PP}{\longrightarrow}}
\newcommand{\distr}{\stackrel{\cL}{\longrightarrow}}
\newcommand{\distre}{\stackrel{\cL}{=}}

\newcommand{\lu}{\stackrel{\slu}{\longrightarrow}}
\newcommand{\as}{\stackrel{{\mathrm{a.s.}}}{\longrightarrow}}

\newcommand{\bbone}{\mathbbm{1}}
\newcommand{\ns}{{\lfloor ns\rfloor}}
\newcommand{\nt}{{\lfloor nt\rfloor}}
\newcommand{\nT}{{\lfloor nT\rfloor}}
\newcommand{\nhalf}{{\lfloor n/2\rfloor}}
\newcommand{\proofend}{\hfill\mbox{$\Box$}}

\numberwithin{equation}{section}

\theoremstyle{change} \theorembodyfont{\em}
\newtheorem{Lem}{Lemma.}[section]
\newtheorem{Thm}{Theorem.}[section]
\newtheorem{Pro}{Proposition.}[section]
\newtheorem{Cor}{Corollary.}[section]
\newtheorem{Def}{Definition.}[section]

\theorembodyfont{\rm}
\newtheorem{Rem}{Remark.}[section]

\begin{document}

\begin{center}
 {\bfseries\Large Asymptotic behavior of CLS estimators \\[2mm]
                   for unstable INAR(2) models} \\[5mm]

 {\sc\large M\'aty\'as $\text{Barczy}^{*,\diamond}$,
            \ M\'arton $\text{Isp\'any}^*$, \ Gyula $\text{Pap}^{\star}$}
\end{center}

\vskip0.2cm

\noindent * Faculty of Informatics, University of Debrecen,
            Pf.~12, H--4010 Debrecen, Hungary.

\noindent $\star$ Bolyai Institute, University of Szeged,
            Aradi v\'ertan\'uk tere 1, H--6720 Szeged, Hungary.

\noindent e--mails: barczy.matyas@inf.unideb.hu (M. Barczy),
                    ispany.marton@inf.unideb.hu (M. Isp\'any),
                    papgy@math.u-szeged.hu (G. Pap).

\noindent $\diamond$ Corresponding author.

%\vskip0.2cm

%\centerline{\sl February 2, 2012.}

%\renewcommand{\thefootnote}{}
\vspace*{0.2cm}
\footnote{The authors have been supported by the Hungarian
 Chinese Intergovernmental S \& T Cooperation Programme for 2011-2013
 under Grant No.\ 10-1-2011-0079.
M. Isp\'any has been partially supported by the
 T\'AMOP-4.2.2.C-11/1/KONV-2012-0001 project. 
The project has been supported by the European Union, co-financed by the
 European Social Fund.
The research of M. Barczy was realized in the frames of
 T\'AMOP 4.2.4.\ A/2-11-1-2012-0001 ,,National Excellence Program --
 Elaborating and operating an inland student and researcher personal support
 system''.
The project was subsidized by the European Union and co-financed by the
 European Social Fund.}

\vspace*{-10mm}

\begin{abstract}
In this paper the asymptotic behavior of the conditional least squares
 estimators of the autoregressive parameters \ $(\alpha, \beta)$, \
 of the stability parameter \ $\varrho := \alpha + \beta$, \ and of the
 mean \ $\mu$ \ of the innovation \ $\vare_k$, $k \in \NN$,
 \ for an unstable integer-valued autoregressive process
 \ $X_k = \alpha \circ X_{k-1} + \beta \circ X_{k-2} + \vare_k$, \ $k \in \NN$,
 \ is described.
The limit distributions and the scaling factors are different according to the
 following three cases: (i) decomposable, (ii) indecomposable but not
 positively regular, and (iii) positively regular models.
\end{abstract}

\section{Introduction}
\label{section_intro}

The theory and practice of statistical inference for integer-valued time series
 models are rapidly developing and important topics of the modern theory of
 statistics.
A number of results are now available in specialized monographs and review
 papers, to name a few, see, e.g., Steutel and van Harn \cite{SteHar} and
 Wei{\ss} \cite{Wei1}.
Among the most successful integer-valued time series models proposed in the
 literature we mention the INteger-valued AutoRegressive model of order \ $p$
 \ (\INARp).
This model was first introduced by McKenzie \cite{McK} and Al-Osh and Alzaid
 \cite{AloAlz1} for the case \ $p=1$.
The INAR(1) model has been investigated by several authors.
The more general \INARp\ processes were first introduced by Alzaid and Al-Osh
 \cite{AloAlz2}.
In their setup the autocorrelation structure of the process corresponds to
 that of an ARMA($p,p-1$) process.
Another definition of an \INARp\ process was proposed independently by Du and
 Li \cite{DuLi} and by Gauthier and Latour \cite{GauLat} and
 Latour \cite{Lat2}, and is different from that of
 Alzaid and Al-Osh \cite{AloAlz2}.
In Du and Li's setup the autocorrelation structure of an \INARp\ process is
 the same as that of an AR($p$) process.
The setup of Du and Li \cite{DuLi} has been followed by most of the authors,
 and our approach will also be the same.
In Barczy et al.\ \cite{BarIspPap0} we investigated the asymptotic behavior of
 unstable INAR($p$) processes, i.e., when the characteristic polynomial has a
 unit root.
Under some natural assumptions we proved that the sequence of appropriately
 scaled random step functions formed from an unstable \INARp\ process
 converges weakly towards a squared Bessel process.
This limit process is a continuous time branching process with immigration
 also known as the square-root process or the Cox--Ingersoll--Ross process.

Parameter estimation for \INARp\ models has a long history.
Franke and Seligmann \cite{FraSel} analyzed conditional maximum likelihood
 estimator of some parameters (including the autoregressive parameter) for
 stable INAR(1) models with Poisson innovations.
Du and Li \cite[Theorem 4.2]{DuLi} proved asymptotic normality of the
 conditional least squares (CLS) estimator of the autoregressive parameters
 for stable \INARp\ models (see also Latour \cite[Proposition 6.1]{Lat2}),
 Br\"{a}nn\"{a}s and Hellstr\"{o}m \cite{BraHel} considered generalized method
 of moment estimation.
Silva and Oliveira \cite{SilOli} proposed a frequency domain based estimator
 of the autoregressive parameters for stable \INARp\ models with Poisson
 innovations.
Isp\' any et al.\ \cite{IspPapZui0}, \cite{IspPapZui1} derived asymptotic
 inference for nearly unstable INAR(1) models which has been refined by
 Drost et al.\ \cite{DroAkkWer2} later.
In \cite{IspPapZui0} the mean of the innovation was supposed to be known,
 while in \cite{IspPapZui1} both the autoregressive parameter and the mean of
 the innovation have been estimated jointly.
Drost et al.\ \cite{DroAkkWer1} studied asymptotically efficient estimation of
 the parameters for stable \INARp\ models.
The stability parameter \ $\varrho := \alpha_1 + \cdots + \alpha_p$ \ of an
 \INARp\ model with autoregressive parameters \ $(\alpha_1, \ldots, \alpha_p)$
 \ has not been treated yet, but this stability parameter is well investigated
 in case of unstable AR($p$) processes, see the unit root tests, e.g., in
 Hamilton \cite[Section 17, Table 17.3, Case 1]{Ham}.
Namely, for the simplicity in case of \ $p = 1$, \ if \ $(Y_k)_{k\geq 0}$ \ is
 an AR(1) process, i.e., \ $Y_k = \varrho Y_{k-1} +\zeta_k$, \ $k\geq 1$, \ with
 \ $Y_0 := 0$ \ and an i.i.d. sequence \ $(\zeta_k)_{k\geq 1}$ \ having mean
 \ $0$ \ and positive variance, then the ordinary least squares estimator of
 the stability parameter \ $\varrho$ \ based on the sample
 \ ${\bY}_n := (Y_1, \ldots, Y_n)$ \ takes the form
 \[
   \hvarrho_n(\bY_n) = \frac{\sum_{k=1}^n Y_{k-1}Y_k}{\sum_{k=1}^n Y_k^2}, \qquad
   n \geq 1 ,
 \]
 see, e.g., Hamilton \cite[17.4.2]{Ham}, and, by Hamilton \cite[17.4.7]{Ham},
 in the unstable case, i.e., when \ $\varrho = 1$,
 \[
   n (\hvarrho_n(\bY_n) - 1)
   \distr \frac{\int_0^1 \cW_t \, \dd \cW_t}{\int_0^1 \cW_t^2 \, \dd t}
   \qquad \text{as \ $n \to \infty$,}
 \]
 where \ $(\cW_t)_{t\geq0}$ \ is a standard Wiener process and \ $\distr$
 \ denotes convergence in distribution.
Here \ $n (\hvarrho_n(\bY_n) - 1)$ \ is known as the Dickey--Fuller statistics.
In this paper the asymptotic behavior of the CLS estimators of the
 autoregressive and the stability parameters and of the mean of the innovation for unstable \INARtwo\ models is
 described (see our main results in Section \ref{section_main_results}) which can be considered as a first step of
 examining this question for general unstable INAR(p) processes and more generally for
 critical multitype branching processes.
We call the attention that in case of unstable INAR(2) processes new types of limit distribution occur
 (see Theorem \ref{main}) compared to those of unstable AR($p$) processes.

First we recall \INARtwo\ models.
Let \ $\ZZ_+$, \ $\NN$, \ $\RR$ \ and \ $\RR_+$ \ denote the set of
 non-negative integers, positive integers, real numbers and non-negative real
 numbers, respectively.
Every random variable will be defined on a fixed probability space
 \ $(\Omega, \cA, \PP)$.

\begin{Def}
Let \ $(\vare_k)_{k\in\NN}$ \ be an independent and identically distributed
 (i.i.d.) sequence of non-negative integer-valued random variables, and let
 \ $(\alpha, \beta) \in [0,1]^2$.
\ An \INARtwo\ time series model with autoregressive parameters
 \ $(\alpha, \beta)$ \ and innovations \ $(\vare_k)_{k\in\NN}$ \ is a stochastic
 process \ $(X_k)_{k \geq -1}$ \ given by
 \begin{align}\label{INAR2}
   X_k = \sum_{j=1}^{X_{k-1}} \xi_{k,j}
         + \sum_{j=1}^{X_{k-2}} \eta_{k,j} + \vare_k , \qquad k \in \NN ,
 \end{align}
 where for all \ $k\in\NN$, \ $(\xi_{k,j})_{j\in\NN}$ \ and
 \ $(\eta_{k,j})_{j\in\NN}$ \ are sequences of i.i.d.\ Bernoulli random variables
 with mean \ $\alpha$ \ and \ $\beta$, \ respectively, \ such that these
 sequences are mutually independent and independent of the sequence
 \ $(\vare_k)_{k\in\NN}$, \ and \ $X_0$ \ and \ $X_{-1}$ \ are non-negative
 integer-valued random variables independent of the sequences
 \ $(\xi_{k,j})_{j\in\NN}$, \ $(\eta_{k,j})_{j\in\NN}$, \ $k\in\NN$, \ and
 \ $(\vare_k)_{k\in\NN}$.
\end{Def}

The \INARtwo\ model \eqref{INAR2} can be written in another way using the
 binomial thinning operator \ $\circ$
 \ (due to Steutel and van Harn \cite{SteHar}) which we recall now.
Let \ $X$ \ be a non-negative integer-valued random variable.
Let \ $(\xi_j)_{j\in\NN}$ \ be a sequence of i.i.d.\ Bernoulli random variables
 with mean \ $\alpha\in[0,1]$.
\ We assume that the sequence \ $(\xi_j)_{j\in\NN}$ \ is independent of \ $X$.
\ The non-negative integer-valued random variable \ $\alpha\,\circ X$
 \ is defined by
 \[
   \alpha\circ X
     :=\begin{cases}
        \sum\limits_{j=1}^X\xi_j, & \quad \text{if \ $X>0$},\\[2mm]
         0, & \quad \text{if \ $X=0$}.
       \end{cases}
 \]
The sequence \ $(\xi_j)_{j\in\NN}$ \ is called a counting sequence.
Then the \INARtwo\ model \eqref{INAR2} takes the form
 \[
    X_k = \alpha \circ X_{k-1} + \beta \circ X_{k-2} + \vare_k ,
    \qquad k \in \NN .
 \]
Note that the above form of the \INARtwo\ model is quite analogous with a
 usual \ARtwo\ process (another slight link between them is the similarity of
 some conditional expectations, see \eqref{seged1}).

For the sake of simplicity we consider a zero start \INARtwo\ process, that is
 we suppose \ $X_0 = X_{-1} = 0$.
\ The general case of nonzero initial values may be handled in a similar way,
 but we renounce to consider it.

In the sequel we always assume \ $\EE(\vare_1^2) < \infty$.
\ Let us denote the mean and variance of \ $\vare_1$ \ by \ $\mu$ \ and
 \ $\sigma^2$.
\ Further, we assume \ $\mu > 0$, \ otherwise \ $X_k = 0$ \ for all
 \ $k \in \NN$.

Based on the asymptotic behavior of \ $\EE(X_k)$ \ as \ $k \to \infty$
 \ described in Barczy et al.\ \cite[Proposition 2.6]{BarIspPap0}, we
 distinguish three types of \INARtwo\ models.
The asymptotic behavior of \ $\EE(X_k)$ \ as \ $k \to \infty$ \ is determined
 by the spectral radius \ $r$ \ of the matrix
 \begin{equation}\label{bA}
  A :=\begin{bmatrix}
       \alpha & \beta \\
       1 & 0 \\
      \end{bmatrix},
 \end{equation}
 i.e., by the maximum of the modulus of the eigenvalues of \ $A$.
\ The case \ $r < 1$, \ when \ $\EE(X_k)$ \ converges to a finite
 limit as \ $k \to \infty$, \ is called \emph{stable} or
 \emph{asymptotically stationary}, whereas the cases \ $r = 1$, \ when
 \ $\EE(X_k)$ \ tends linearly to \ $\infty$, \ and \ $r > 1$, \ when
 \ $\EE(X_k)$ \ converges to \ $\infty$ \ with an exponential rate, are called
 \emph{unstable} and \emph{explosive}, respectively.
It is easy to check that \ $r < 1$, \ $r = 1$, \ and \ $r > 1$ \ are
 equivalent with \ $\varrho < 1$, \ $\varrho = 1$, \ and \ $\varrho > 1$,
 \ respectively, where \ $\varrho := \alpha + \beta$ \ is called the
 \emph{stability parameter}, \ see
 Barczy et al.~\cite[Proposition 2.2]{BarIspPap0}.

We also note that an \INARtwo\ process can be considered as a special 2-type
 branching process with immigration.
Namely, by \eqref{INAR2},
  \[
     \begin{bmatrix}
       X_k \\
       X_{k-1} \\
     \end{bmatrix}
     =\sum_{j=1}^{X_{k-1}}
       \begin{bmatrix}
          \xi_{k,j} \\
           1 \\
     \end{bmatrix}
      + \sum_{j=1}^{X_{k-2}}
       \begin{bmatrix}
          \eta_{k,j} \\
           0 \\
       \end{bmatrix}
      + \begin{bmatrix}
          \vare_k \\
           0 \\
       \end{bmatrix} ,
       \qquad k\in\NN,
  \]
 and hence the so-called mean matrix of an \INARtwo\ process with
 autoregressive parameters \ $(\alpha, \beta)$ \ (considered as a 2-type
 branching process) is nothing else but \ $A$.
\ This process is called \emph{positively regular} if there is a positive
 integer \ $k \in \NN$ \ such that the entries of \ $A^k$ \ are positive
 (see Kesten and Stigum \cite{KesSti1}), which is equivalent with
 \ $\alpha > 0$ \ and \ $\beta > 0$.
\ The model is called \emph{decomposable} if the matrix \ $A$ \ is
 decomposable (see Kesten and Stigum \cite{KesSti3}), which is equivalent with
 \ $\beta = 0$.
\ If \ $\alpha = 0$ \ and \ $\beta > 0$, \ then the process is
 \emph{indecomposable but not positively regular}
 (see Kesten and Stigum \cite{KesSti2}).
If \ $\alpha > 0$ \ and \ $\beta = 0$, \ then the decomposable process
 \ $(X_k)_{k\geq -1}$ \ is an INAR(1) process with autoregressive parameter
 \ $\alpha$.
\ If \ $\alpha = 0$ \ and \ $\beta > 0$, \ then the indecomposable process
 \ $(X_k)_{k\geq -1}$ \ takes the form
 \[
   X_k = \beta \circ X_{k-2} + \vare_k, \qquad k \in \NN ,
 \]
 and hence the subsequences \ $(X_{2k-j})_{k\geq0}$, \ $j \in \{0, 1\}$, \ form
 independent positively regular INAR(1) processes with autoregressive
 parameter \ $\beta$ \ such that \ $X_{-j}=0$, \ $j \in \{0, 1\}$.
\ For more details of this classification of \INARtwo\ processes, see Appendix
 \ref{app_A}.

Next we give an overview of the structure of the paper.
Section \ref{section_main_results} contains our main results,
 see Theorem \ref{main} for unstable and positively regular \INARtwo\ processes,
 Theorem \ref{10main} for unstable and decomposable \INARtwo\ processes,
 and Theorem \ref{01main} for unstable, indecomposable but not positively regular ones.
In order to highlight our main results, the preliminaires and (technical) details on CLS estimators
 are presented only after our main results, see Section \ref{section_estimators}.
In Theorems \ref{main_Ad}, \ref{10main_Ad} and \ref{01main_Ad} of Section \ref{section_proof}
 we present joint asymptotic behaviours of the building blocks of the CLS estimators
 (according to the above mentioned three cases), and by applying a version of the continuous mapping theorem
 (which is formulated for completeness in Appendix \ref{app_B})
 we show how one can derive Theorems \ref{main}, \ref{10main} and \ref{01main} using these theorems.
Section \ref{section_proof_main} is devoted to the proof of Theorem \ref{main_Ad} which is based on
 Lemma \ref{XV_main_VV} and Theorem \ref{main_conv}.
Due to its length, the proof of Theorem \ref{main_conv} is given separately in Section \ref{section_proof_main_conv}.
Sections \ref{section_proof_10main} and \ref{section_proof_01main} are devoted to the proofs of Theorem \ref{10main_Ad}
 and Theorem \ref{01main_Ad}, respectively.
In Section \ref{section_moments} we present estimates for the moments of the
 processes involved, these estimates are used throughout the paper.
In Appendix \ref{section_conv_step_processes} we recall a result about convergence of
 random step processes noting that the proof of Theorem \ref{main_conv} is based on this result.

\section{Main results}
\label{section_main_results}

In what follows we always assume \ $\varrho = \alpha + \beta = 1$, \ that is,
 the process \ $(X_k)_{k\geq-1}$ \ is unstable.

For each \ $n \in \NN$, \ any CLS estimator
 \ $(\halpha_n(\bX_n), \hbeta_n(\bX_n), \hmu_n(\bX_n))$
 \ of the autoregressive parameters
 \ $(\alpha, \beta)$
 \ and of the mean \ $\mu$ \ of the innovation based on a sample
 \ $\bX_n := (X_1,\ldots,X_n)$ \ has
 the form
 \[
   \begin{bmatrix}
    \halpha_n(\bX_n) \\
    \hbeta_n(\bX_n) \\
    \hmu_n(\bX_n)
   \end{bmatrix}
   = \left( \sum_{k=1}^n
             \begin{bmatrix}
              X_{k-1}^2 & X_{k-1} X_{k-2} & X_{k-1} \\
              X_{k-1} X_{k-2} & X_{k-2}^2 & X_{k-2} \\
              X_{k-1} & X_{k-2} & 1
             \end{bmatrix} \right) ^{-1}
     \sum_{k=1}^n
             \begin{bmatrix}
              X_k X_{k-1} \\
              X_k X_{k-2} \\
              X_k
             \end{bmatrix}
 \]
 on the set \ $\{ \omega \in \Omega : \sum_{k=1}^n X_{k-2}(\omega)^2 > 0\}$
 \ with \ $\lim_{n \to \infty} \PP\left( \sum_{k=1}^n X_{k-2}^2 > 0 \right) = 1$,
 \ see Proposition \ref{ExUn}.
Moreover, for each \ $n \in \NN$, \ any CLS estimator of the stability
 parameter \ $\varrho$ \ takes the form
 \[
   \hvarrho_n(\bX_n) = \halpha_n(\bX_n) + \hbeta_n(\bX_n)
 \]
 on the set \ $\{ \omega \in \Omega : \sum_{k=1}^n X_{k-2}(\omega)^2 > 0\}$,
 \ see Section \ref{section_estimators}.

\begin{Thm}\label{main}
Let \ $(X_k)_{k \geq -1}$ \ be an \INARtwo\ process with autoregressive
 parameters \ $(\alpha, \beta) \in (0,1)^2$ \ such that \ $\alpha + \beta = 1$
 \ (hence it is unstable and positively regular).
Suppose that \ $X_0 = X_{-1} = 0$, \ $\EE(\vare_1^8) < \infty$ \ and
 \ $\mu > 0$.
\ Then
 \begin{equation}\label{rho}
   n\bigl(\hvarrho_n(\bX_n) - 1\bigr)
   \distr
   \frac{\sqrt{2\alpha\beta} \int_0^1 \cX_t^{3/2} \, \dd\cW_t
         - [(1 + \beta) \cX_1 - \mu] \int_0^1 \cX_t \, \dd t}
        {\int_0^1 \cX_t^2 \, \dd t - \bigl(\int_0^1 \cX_t \, \dd t\bigr)^2}
 \end{equation}
 \begin{equation}\label{alpha,beta}
   \begin{bmatrix}
    n^{1/2} (\halpha_n(\bX_n) - \alpha) \\
    n^{1/2} (\hbeta_n(\bX_n) - \beta)
   \end{bmatrix}
   \distr
   \sqrt{\alpha (1 + \beta)}
    \frac{\int_0^1 \cX_t \, \dd\tcW_t}
         {\int_0^1 \cX_t \, \dd t}
    \begin{bmatrix} -1 \\ 1 \end{bmatrix}
 \end{equation}
 and
 \begin{equation}\label{mu}
   \hmu_n(\bX_n) - \mu
   \distr
   \frac{-\sqrt{2\alpha\beta} \int_0^1 \cX_t \, \dd t
         \int_0^1 \cX_t^{3/2} \, \dd\cW_t
         + [(1 + \beta) \cX_1 - \mu] \int_0^1 \cX_t^2 \, \dd t}
        {\int_0^1 \cX_t^2 \, \dd t - \bigl(\int_0^1 \cX_t \, \dd t\bigr)^2}
 \end{equation}
 as \ $n \to \infty$,
 \ where \ $(\cX_t)_{t\in\RR_+}$ \ is the unique strong
 solution of the stochastic
 differential equation (SDE)
 \begin{equation}\label{X}
   \dd\cX_t
   = \frac{1}{1 + \beta}
     \Big( \mu \, \dd t
     + \sqrt{2 \alpha \beta \cX_t^+} \, \dd\cW_t \Big),
   \qquad t \in \RR_+ ,
 \end{equation}
 with initial value \ $\cX_0 = 0$, \ where \ $(\cW_t)_{t\in\RR_+}$,
 \ $(\tcW_t)_{t\in\RR_+}$ \ are independent standard Wiener processes, and
 \ $x^+$ \ denotes the positive part of \ $x \in \RR$.
\end{Thm}

\begin{Rem}\label{REMARK5}
The moment condition \ $\EE(\vare_1^8) < \infty$ \ in Theorem \ref{main} seems
 to be too strong, but we call the attention that the process
 \ $(X_k)_{k\geq-1}$ \ can be considered as a heteroscedastic time series.
Indeed, \ $X_k = \alpha X_{k-1} + \beta X_{k-2} + M_k + \mu$, \ see \eqref{regr},
 and by \eqref{Mcond},
 \ $\EE(M_k^2\mid\cF_{k-1})
    = \alpha(1-\alpha) X_{k-1} + \beta(1-\beta) X_{k-2} + \sigma^2$,
 \ $k \in \NN$.
\ That is why we think that the behavior of the process \ $(X_k)_{k\geq-1}$ \ is
 similar to GARCH models, where, even in the stable case, high moment
 conditions are needed for convergence of estimators such as the quasi-maximum
 likelihood estimator in Hall and Yao \cite{HalYao} or the Whittle estimator
 in Mikosch and Straumann \cite{MikStr}. 
\proofend
\end{Rem}

\begin{Rem}\label{REMARK1}
The SDE \eqref{X} has a unique strong solution \ $(\cX_t^{(x)})_{t\geq 0}$ \ for
 all initial values  \ $\cX_0^{(x)} = x \in \RR$.
\ Indeed, since \ $| \sqrt{x} - \sqrt{y}|\leq \sqrt{|x - y|}$
 \ for \ $x, y \geq 0$, \ the coefficient functions
 \ $\RR \ni x \mapsto \mu /(1+\beta)$ \ and
 \ $\RR \ni x \mapsto \sqrt{2\alpha\beta x^+}/(1+\beta)$ \ satisfy conditions
 of part (ii) of Theorem 3.5 in Chapter IX in Revuz and Yor \cite{RevYor} or
 the conditions of Proposition 5.2.13 in Karatzas and Shreve \cite{KarShr}.
Further, by the comparison theorem
 (see, e.g., Revuz and Yor \cite[Theorem 3.7, Chapter IX]{RevYor}), if the
 initial value \ $\cX_0^{(x)} = x$ \ is nonnegative, then \ $\cX_t^{(x)}$ \ is
 nonnegative for all \ $t \in \RR_+$ \ with probability one.
Hence \ $\cX_t^+$ \ may be replaced by \ $\cX_t$ \ under the square root in
 \eqref{X}.
The unique strong solution of the SDE \eqref{X} is known as a squared Bessel
 process, a squared-root process or a Cox--Ingersoll--Ross (CIR) process.
\proofend
\end{Rem}

\begin{Rem}\label{REMARK2}
By It\^o's formula and Remark \ref{REMARK1},
 \ $\cM_t := (1 + \beta) \cX_t - \mu t$, \ $t \in \RR_+$, \ is the
 unique
 strong solution of the SDE
 \begin{equation}\label{M}
  \dd\cM_t
  = \sqrt{\frac{2\alpha\beta}{1+\beta} (\cM_t + \mu t)^+} \, \dd\cW_t ,
  \qquad t \in \RR_+ ,
 \end{equation}
 with initial value \ $\cM_0 = 0$, \ and \ $(\cM_t + \mu t)^+$ \ may be
 replaced by \ $\cM_t + \mu t$ \ under the square root in \eqref{M}.
Hence \ $\dd\cM_t = \sqrt{2 \alpha \beta \cX_t} \, \dd\cW_t$, \ and the
 convergences \eqref{rho} and \eqref{mu}
 can also be formulated as
 \begin{equation}\label{primitive_rho_M}
  n (\hvarrho_n(\bX_n) - 1)
  \distr
  \frac{\int_0^1 \cX_t \, \dd\cM_t - \cM_1 \int_0^1 \cX_t \, \dd t}
       {\int_0^1 \cX_t^2 \, \dd t - \bigl(\int_0^1 \cX_t \, \dd t\bigr)^2}
  \qquad \text{as \ $n \to \infty$,}
 \end{equation}
 \begin{equation}\label{primitive_mu_M}
  \hmu_n(\bX_n) - \mu
  \distr
  \frac{-\int_0^1 \cX_t \, \dd t \int_0^1 \cX_t \, \dd\cM_t
        + \cM_1 \int_0^1 \cX_t^2 \, \dd t}
       {\int_0^1 \cX_t^2 \, \dd t - \bigl(\int_0^1 \cX_t \, \dd t\bigr)^2}
  \qquad \text{as \ $n \to \infty$.}
\vspace*{-6.5mm}
 \end{equation}
\proofend
\end{Rem}

The next theorem contains our result for decomposable unstable
 \INARtwo\ processes.

\begin{Thm}\label{10main}
Let \ $(X_k)_{k\geq-1}$ \ be an \INARtwo\ process with autoregressive parameters
 \ $(1, 0)$ \ (hence it is unstable and decomposable).
Suppose that \ $X_0 = X_{-1} = 0$, \ $\EE(\vare_1^4) < \infty$ \ and
 \ $\mu > 0$.
\ Then
 \begin{align}\label{10rho}
  n^{3/2} \bigl( \hvarrho_n(\bX_n) - 1 \bigr)
  \distr \cN_1\left(0, \, \frac{12 \sigma^2}{\mu^2} \right)
  \qquad \text{as \ $n \to \infty$,}
 \end{align}
 \begin{align}\label{10,alpha,beta}
  \begin{bmatrix}
   n^{1/2} (\halpha_n(\bX_n) - 1) \\
   n^{1/2} \hbeta_n(\bX_n)
  \end{bmatrix}
  \distr Z \begin{bmatrix} -1 \\ 1 \end{bmatrix}
  \qquad \text{as \ $n \to \infty$,}
 \end{align}
 and
 \begin{align}\label{10mu}
  n^{1/2} \bigl( \hmu_n(\bX_n) - \mu \bigr) \distr \cN_1(0, \mu^2 + 4 \sigma^2)
  \qquad \text{as \ $n \to \infty$,}
 \end{align}
 where \ $Z$ \ is a standard normally distributed random variable.
\end{Thm}

\begin{Rem}
Note that an unstable and decomposable \INARtwo\ process has autoregressive
 parameters \ $(1,0)$, \ i.e.,
 it is actually an unstable \textup{INAR(1)} process.
However, we call the attention that the asymptotic behaviour of the estimators
 \ $\hvarrho_n(\bX_n)$, \ $(\halpha_n(\bX_n),\hbeta_n(\bX_n))$ \ and \ $\hmu_n(\bX_n)$ \ as \ $n\to\infty$ \ 
 in Theorem \ref{10main} can not be derived from the
 corresponding results for an unstable \textup{INAR(1)} process, since the
 CLS estimator of the coefficient (which can also be considered as the stability parameter) of an \textup{INAR(1)} process is different from \ $\hvarrho_n(\bX_n)$, \ see, e.g., Isp\'any et al.~\cite{IspPapZui1}.
Remark also that the CLS estimator of the coefficient of an unstable \textup{INAR(1)} process is also asymptotically normal 
 with the same scaling \ $n^{3/2}$, \ but the asymptotic variance \ $3\sigma^2/\mu^2$ \ is different from the corresponding one
 \ $12\sigma^2/\mu^2$ \ for an unstable and decomposable \INARtwo\ process, see Isp\'any et al.~\cite[Theorem 2.1]{IspPapZui1}.\proofend
\end{Rem}

The last theorem contains our result for unstable, indecomposable but not
 positively regular \INARtwo\ processes.

\begin{Thm}\label{01main}
Let \ $(X_k)_{k \geq -1}$ \ be an \INARtwo\ process with autoregressive
 parameters \ $(0,1)$
 \ (hence it is unstable, indecomposable but not positively regular).
Suppose that \ $X_0 = X_{-1} = 0$, \ $\EE(\vare_1^2) < \infty$ \ and
 \ $\mu > 0$.
\ Then
 \begin{align}\label{01rho}
  n^{3/2} \bigl( \hvarrho_n(\bX_n) - 1 \bigr)
  \distr
  \cN_1\left(0, \, \frac{48 \sigma^2}{\mu^2} \right)
  \qquad \text{as \ $n \to \infty$,}
 \end{align}
 \begin{align}\label{01,alpha,beta}
  \begin{bmatrix}
   n \halpha_n(\bX_n) \\
   n (\hbeta_n(\bX_n) - 1)
  \end{bmatrix}
  \distr
  \frac{\int_0^1 \cW_t \, \dd \cW_t}{\int_0^1 \cW_t^2 \, \dd t}
  \begin{bmatrix} -1 \\ 1 \end{bmatrix}
  \qquad \text{as \ $n \to \infty$,}
 \end{align}
 and
 \begin{align}\label{01mu}
  n^{1/2} \bigl( \hmu_n(\bX_n) - \mu \bigr)
  \distr
  \cN_1\left(0, \, 4 \sigma^2 \right)
  \qquad \text{as \ $n \to \infty$,}
 \end{align}
 where \ $(\cW_t)_{t \in \RR_+}$ \ is a standard Wiener process.
\end{Thm}

\begin{Rem}\label{REMARK3}
We note that in all unstable cases the limit distributions for the estimators
 of the autoregressive parameters are concentrated on the same line
 \ $\{(x, y) \in \RR^2 : x + y = 0\}$.
\ However, these limit distributions are pairwise different.
Surprisingly, both in the unstable positively regular case and in the
 unstable decomposable case the scaling factor is \ $\sqrt{n}$, \ while in the
 unstable, indecomposable but not positively regular case it is \ $n$.
\ In the stable case this factor is again \ $\sqrt{n}$
 \ (see Du and Li \cite[Theorem 4.2]{DuLi} or Latour
 \cite[Proposition 6.1]{Lat2}).
The reason of this strange phenomena can be understood from the asymptotic
 behavior of the sequence \ $(\bA_n, \bd_n)_{n\in\NN}$ \ of random vectors
 defined and analyzed in Sections \ref{section_estimators},
 \ref{section_proof}, \ref{section_proof_main}, \ref{section_proof_10main} and
 \ref{section_proof_01main}.
Namely, the scaling factor for the entries of the matrices \ $(\bA_n)_{n\in\NN}$
 \ as well as for the entries of the vectors \ $(\bd_n)_{n\in\NN}$ \ are different.
In order to get over these difficulties, we use the canonical form of the
 process \ $(X_k)_{k\in\NN}$ \ due to Sims, Stock and Watson \cite{SimStoWat}.
Further, one of the decisive tools in deriving the needed asymptotic behavior is a good
 bound for the moments of the involved processes, see Corollary
 \ref{EEX_EEU_EEV}.
\proofend
\end{Rem}

\begin{Rem}\label{REMARK4}
We recall that the distribution of
 \ $\int_0^1 \cW_t \, \dd \cW_t \big/ \int_0^1 \cW_t^2 \, \dd t$
 \ in Theorem \ref{01main} agrees with the limit distribution of the
 Dickey--Fuller statistics for unit root test of AR(1) time series, see, e.g.,
 Hamilton \cite[17.4.2 and 17.4.7]{Ham} or Tanaka
 \cite[(7.14) and Theorem 9.5.1]{Tan}.
The limit distribution in \eqref{alpha,beta} is also a fraction of two
 stochastic integrals, but it contains two independent standard Wiener
 processes.
This phenomena is very similar to the appearing of two independent standard
 Wiener processes in limit theorems for CLS estimators of the variance of the
 offspring and immigration distributions for critical branching processes with
 immigration in Winnicki \cite[Theorems 3.5 and 3.8]{Win}.
Finally, note that the limit distribution of the CLS estimator of the
 autoregressive parameters \ $(\alpha, \beta)$ \ is symmetric in Theorems
 \ref{main} and \ref{10main}, and non-symmetric in Theorem \ref{01main}.
Indeed, since \ $(\cW_t)_{t\in\RR_+}$ \ and  \ $(\tcW_t)_{t\in\RR_+}$ \ are
 independent, by the SDE \eqref{X}, the processes \ $(\cX_t)_{t\in\RR_+}$ \ and
 \ $(\tcW_t)_{t\in\RR_+}$ \ are also independent, which yields that the limit
 distribution of the CLS estimator of \ $(\alpha, \beta)$ \ is symmetric in
 Theorem \ref{main}.
\proofend
\end{Rem}

\begin{Rem}
We note that the CLS estimator \ $\hvarrho_n(\bX_n)$ \ of \ $\varrho$, 
 \ and the CLS estimator \ $(\halpha_n(\bX_n),\hbeta_n(\bX_n)$ \ of \ $(\alpha,\beta)$ \  
 are asymptotically weakly consistent as \ $n\to\infty$ \ in Theorems \ref{main}, \ref{10main} and \ref{01main}.   
The CLS estimator \ $\hmu_n(\bX_n)$ \ of \ $\mu$ \ in Theorems 
 \ref{10main} and \ref{01main} is also asymptotically weakly consistent as \ $n\to\infty$, \ 
 however in Theorem \ref{main} it is not asymptotically weakly consistent.
Note that in the case of an unstable INAR(1) model the CLS estimator of the
 mean of the innovation is asymptotically weakly consistent, see
 Isp\'any et al.\ \cite{IspPapZui1}.
Further, we remark that in Theorem \ref{main} the variance \ $\sigma^2$ \ of
 the innovation does not show up in the limit distributions, while in Theorems
 \ref{10main} and \ref{01main} it appears.
Finally, in Theorems \ref{main}, \ref{10main} and \ref{01main} one could prove
 joint convergence as well. 
\proofend
\end{Rem}

\section{CLS estimators}
\label{section_estimators}

For all \ $k \in \ZZ_+$, \ let us denote by \ $\cF_k$ \ the \ $\sigma$-algebra
 generated by the random variables \ $X_{-1},X_0,X_1,\ldots,X_k$.
\ (Note that \ $\cF_0 = \{ \Omega, \emptyset \}$, \ since
 \ $X_0 = X_{-1} = 0$.)
\ By \eqref{INAR2},
 \begin{align}\label{seged1}
   \EE(X_k\mid\cF_{k-1})
   = \alpha X_{k-1} + \beta X_{k-2} + \mu , \qquad k \in \NN .
 \end{align}
Let us introduce the sequence
 \begin{equation}\label{Mk}
  M_k := X_k - \EE(X_k \mid \cF_{k-1})
       = X_k - \alpha X_{k-1} - \beta X_{k-2} - \mu ,
  \qquad k \in \NN ,
 \end{equation}
 of martingale differences with respect to the filtration
 \ $(\cF_k)_{k \in \ZZ_+}$.
\ The process \ $(X_k)_{k \geq -1}$ \ satisfies the recursion
 \begin{equation}\label{regr}
  X_k = \alpha X_{k-1} + \beta X_{k-2} + M_k + \mu ,
  \qquad k \in \NN .
 \end{equation}

For each \ $n \in \NN$, \ a CLS estimator
 \ $(\halpha_n(\bX_n), \hbeta_n(\bX_n), \hmu_n(\bX_n))$
 \ of the parameters
 \ $(\alpha, \beta, \mu)$
 \ based on a sample \ $\bX_n = (X_1, \ldots, X_n)$ \ can be obtained by
 minimizing the sum of squares
 \begin{equation}\label{SumSquares}
   \sum_{k=1}^n \big( X_k - \EE(X_k \mid \cF_{k-1}) \big)^2
   = \sum_{k=1}^n (X_k - \alpha X_{k-1} - \beta X_{k-2} - \mu)^2
 \end{equation}
 with respect to
 \ $(\alpha, \beta, \mu)$ \ over \ $\RR^3$.
\ For all \ $n \in \NN$ \ and \ $x_1, \ldots, x_n \in \RR$, \ let us put
 \[
  \bx_n := (x_1, \ldots, x_n),
 \]
 and in what follows we use the convention
 \[
   x_{-1} := x_0 := 0.
 \]
For all \ $n \in \NN$, \ we define the
 function
 \ $Q_n : \RR^n \times \RR^3 \to \RR$ \ by
 \[
   Q_n(\bx_n ; \alpha', \beta', \mu')
   := \sum_{k=1}^n (x_k - \alpha' x_{k-1} - \beta' x_{k-2} - \mu')^2
 \]
 for all \ $\alpha', \beta', \mu' \in \RR$ \ and \ $\bx_n \in \RR^n$.
\ By definition, for all \ $n \in \NN$, \ a CLS estimator of the parameters
 \ $(\alpha, \beta, \mu)$ \ is a measurable function
 \ $(\halpha_n, \hbeta_n, \hmu_n) : \RR^n \to \RR^3$ \ such that
 \[
   Q_n(\bx_n ; \halpha_n(\bx_n), \hbeta_n(\bx_n), \hmu_n(\bx_n))
   = \inf_{(\alpha', \beta', \mu') \in \RR^3} Q_n(\bx_n ; \alpha', \beta', \mu')
   \qquad \forall\; \bx_n \in \RR^n .
 \]
Since the variance \ $\sigma^2$ \ of the innovation does not appear in the
 conditional expectation \ $\EE(X_k \mid \cF_{k-1})$ \ given in \eqref{seged1},
 and hence, in the definition of \ $Q_n$, \ we do not need to know the value
 of \ $\sigma^2$ \ for the calculation of the CLS estimator of the parameters
 \ $(\alpha, \beta, \mu)$.
 
Next we give the solutions of this extremum problem.

\begin{Lem}\label{CLSE1}
For each \ $n\geq 2$, \ $n \in \NN$, \ any CLS estimator of the parameters
 \ $(\alpha, \beta, \mu)$ \ is a measurable function
 \ $(\halpha_n, \hbeta_n, \hmu_n) : \RR^n \to \RR^3$ \ for which
 \begin{equation}\label{CLS1}
  \begin{bmatrix}
   \halpha_n(\bx_n) \\
   \hbeta_n(\bx_n) \\
   \hmu_n(\bx_n)
  \end{bmatrix}
  = F_n(\bx_n)^{-1} g_n(\bx_n)
 \end{equation}
 if \ $\sum_{k=1}^n x_{k-2}^2 > 0$, \ where
 \[
   F_n(\bx_n)
   := \sum_{k=1}^n
       \begin{bmatrix}
        x_{k-1}  \\
        x_{k-2} \\
        1
       \end{bmatrix}
       \begin{bmatrix}
        x_{k-1}  \\
        x_{k-2} \\
        1
       \end{bmatrix}^\top ,
   \qquad
   g_n(\bx_n)
   := \sum_{k=1}^n
       x_k
       \begin{bmatrix}
        x_{k-1}  \\
        x_{k-2} \\
        1
       \end{bmatrix} ,
 \]
 \begin{equation}\label{CLS2}
  \halpha_n(\bx_n) = \frac{x_n}{x_{n-1}} - \frac{1}{n-1} , \qquad
  \hmu_n(\bx_n) = \frac{x_{n-1}}{n-1}
 \end{equation}
 if \ $x_1 = \dots = x_{n-2} = 0$ \ and \ $x_{n-1} \ne 0$, \ and
 \begin{equation}\label{CLS3}
  \hmu_n(\bx_n) = \frac{x_n}{n}
 \end{equation}
 if \ $x_1 = \dots = x_{n-1} = 0$.
\end{Lem}

Note that
 \ $(\halpha_n, \hbeta_n, \hmu_n)$
 \ is not defined uniquely on the set
 \ $\{ \bx_n \in \RR^n : x_1 = \dots = x_{n-2} = 0 \}$.
\ Namely, if \ $x_1 = \dots = x_{n-2} = 0$ \ and \ $x_{n-1} \ne 0$, \ then
 \ $\hbeta_n$ \ can be chosen as an arbitrary measurable function, while if
 \ $x_1 = \dots =  x_{n-1} = 0$, \ then the same holds for
 \ $(\halpha_n, \hbeta_n)$.
We call the attention that Lemma \ref{CLSE1} holds for all types of \INARtwo\
 processes, i.e., it covers the stable, unstable and explosive cases as well.

\noindent
\textbf{Proof of Lemma \ref{CLSE1}.}
For any fixed \ $\bx_n \in \RR^n$ \ with \ $\sum_{k=1}^n x_{k-2}^2 > 0$, \ the
 quadratic function
 \ $\RR^3 \ni (\alpha', \beta', \mu')
          \mapsto Q_n(\bx_n ; \alpha', \beta', \mu')$
 \ can be written in the form
 \[
   Q_n(\bx_n ; \alpha', \beta', \mu')
   = \left( \begin{bmatrix} \alpha' \\ \beta' \\ \mu' \end{bmatrix}
            - F_n(\bx_n)^{-1} g_n(\bx_n) \right)^\top
     \!\!\! F_n(\bx_n)
     \left( \begin{bmatrix} \alpha' \\ \beta' \\ \mu' \end{bmatrix}
            - F_n(\bx_n)^{-1} g_n(\bx_n) \right)
     + \tQ_n(\bx_n) ,
 \]
 where
 \[
   \tQ_n(\bx_n)
   := \sum_{k=1}^n x_k^2 - g_n(\bx_n)^\top F_n(\bx_n)^{-1} g_n(\bx_n) .
 \]
We check that the matrix \ $F_n(\bx_n)$ \ is strictly positive definite.
For this, it is enough to show that \ $\sum_{k=1}^n x_{k-2}^2 > 0$ \ implies that the rank of the system of vectors
 \begin{align}\label{dyadic_polynomial}
   \begin{bmatrix}
    x_{k-1}  \\
    x_{k-2} \\
    1
   \end{bmatrix} , \qquad k \in \{1, \ldots, n\} ,
 \end{align}
 equals 3.
Indeed, if \ $a_i\in\RR^3$, $i\in\{1,\ldots,n\}$, and the rank of the system \ $\{a_1,\ldots,a_n\}$ \ is \ $3$, 
 \ then \ $A:=\sum_{i=1}^n a_ia_i^\top$ \ is strictly positive definite
 which can be checked as follows.
For any \ $z\in\RR^3$, 
 \begin{align*}
  \langle Az,z\rangle 
    = \sum_{i=1}^n \langle a_ia_i^\top z,z\rangle
    = \sum_{i=1}^n \langle a_i^\top z,a_i^\top z\rangle
    \geq 0, 
 \end{align*}  
 and \ $\langle Az,z\rangle = 0$ \ holds if and only if \ $a_i^\top z = 0$, $i=1,\ldots,n$.
\ Since the rank of the system \ $\{a_1,\ldots,a_n\}$ \ is \ $3$, \ we have \ $z=0$.   
     
The rank of the system of vectors in \eqref{dyadic_polynomial} is \ $3$, \ since the rank of the matrix
 \[
   \begin{bmatrix}
    x_{n-1} & x_{n-2} & \cdots & x_2 & x_1 & 0 \\
    x_{n-2} & x_{n-3} & \cdots & x_1 & 0   & 0 \\
    1      & 1      & 1      & 1   & 1   & 1
   \end{bmatrix}
 \]
 equals 3.
Indeed, there exists some \ $i \in \{1, \ldots, n-2\}$ \ such that
 \ $x_i \ne 0$ \ and \ $x_{i-1} = 0$, \ and hence there exists a submatrix with
 negative determinant
 \[
   \begin{vmatrix}
    x_{i+1} & x_i    & 0 \\
    x_i    & x_{i-1} & 0 \\
    1      & 1      & 1
   \end{vmatrix}
   = x_{i-1} x_{i+1} - x_i^2 = - x_i^2 < 0 .
 \]
Hence we obtain \eqref{CLS1}.

For any fixed \ $\bx_n \in \RR^n$ \ with \ $x_1 = \dots = x_{n-2} = 0$ \ and
 \ $x_{n-1} \ne 0$, \ the quadratic function
 \ $\RR^3 \ni (\alpha', \beta', \mu')
    \mapsto Q_n(\bx_n ; \alpha', \beta', \mu')$
 \ can be written in the form
 \[
   Q_n(\bx_n ; \alpha', \beta', \mu')
   = (x_n - \alpha' x_{n-1} - \mu')^2 + (x_{n-1} - \mu')^2 + (n-2) (\mu')^2 ,
   \qquad (\alpha', \beta', \mu') \in \RR^3 .
 \]
The system of equation consisting of the first order partial derivates of \ $Q_n$ \ with respect to 
 \ $\alpha'$ \ and \ $\mu'$ \ takes the form
 \begin{align*}
   &x_n - \alpha'x_{n-1} - \mu' =0,\\
   &x_n - \alpha'x_{n-1}-\mu' + x_{n-1} - \mu' - (n-2)\mu' =0.
 \end{align*}   
Using that \ $n\geq 2$, \ by an easy computation, we conclude \eqref{CLS2}.

If \ $\bx_n \in \RR^n$ \ with \ $x_1 = \dots = x_{n-1} = 0$, \ then
 \ $Q_n(\bx_n ; \alpha', \beta', \mu') = ( x_n - \mu')^2 + (n-1) (\mu')^2$,
 \ which implies \eqref{CLS3}.
\proofend

We note that one could give a different proof of Lemma \ref{CLSE1} as in
 Barczy et al.\ \cite[Lemma 2.1]{BarIspPap1}.

Next we present a result about the existence and uniqueness of
 \ $(\halpha_n(\bX_n), \hbeta_n(\bX_n), \hmu_n(\bX_n))$.

\begin{Pro}\label{ExUn}
Let \ $(X_k)_{k \geq -1}$ \ be an \INARtwo\ process with autoregressive
 parameters \ $(\alpha, \beta) \in [0,1]^2$ \ such that \ $\alpha + \beta = 1$
 \ (hence it is unstable).
Suppose that \ $X_0 = X_{-1} = 0$, \ $\EE(\vare_1^2) < \infty$ \ and
 \ $\mu > 0$.
\ Then
 \[
   \lim_{n \to \infty} \PP\left( \sum_{k=1}^n X_{k-2}^2 > 0 \right) = 1 ,
 \]
 and hence the probability of the existence of a unique CLS estimator
 \ $(\halpha_n(\bX_n), \hbeta_n(\bX_n), \hmu_n(\bX_n))$ \ converges to 1 as
 \ $n \to \infty$, \ and this CLS estimator has the form
 \begin{equation}\label{CLSE}
  \begin{bmatrix}
   \halpha_n(\bX_n) \\
   \hbeta_n(\bX_n) \\
   \hmu_n(\bX_n)
  \end{bmatrix}
  = \bF_n^{-1} \bg_n
 \end{equation}
 on the set \ $\{ \omega \in \Omega : \sum_{k=1}^n X_{k-2}(\omega)^2 > 0\}$,
 \ where
 \[
   \bF_n := F_n(\bX_n)
          = \sum_{k=1}^n
             \begin{bmatrix}
              X_{k-1}^2 & X_{k-1} X_{k-2} & X_{k-1} \\
              X_{k-1} X_{k-2} & X_{k-2}^2 & X_{k-2} \\
              X_{k-1} & X_{k-2} & 1
             \end{bmatrix} , \qquad
   \bg_n := g_n(\bX_n)
          = \sum_{k=1}^n
             \begin{bmatrix}
              X_k X_{k-1} \\
              X_k X_{k-2} \\
              X_k
             \end{bmatrix} .
 \]
\end{Pro}

\noindent
\textbf{Proof.}
First we prove the statements for \ $(\alpha, \beta) \in (0, 1)^2$.
\ For each \ $n \in \NN$, \ consider the random step process
 \[
   \cX^{(n)}_t := n^{-1} X_\nt , \qquad t \in \RR_+ ,
 \]
 where \ $\lfloor x \rfloor$ \ denotes the integer part of a real number
 \ $x \in \RR$.
\ By Barczy et al.\ \cite[Theorem 3.1]{BarIspPap0} we have
 \begin{equation}\label{convX}
   \cX^{(n)} \distr \cX \qquad \text{as \ $n \to \infty$,}
 \end{equation}
 where the process \ $(\cX_t)_{t\in\RR_+}$ \ is the unique strong solution of
 the SDE \eqref{X} with initial value \ $\cX_0=0$.
\ Next we show that
 \begin{align}\label{seged2}
   \frac{1}{n^3} \sum_{k=1}^n X_{k-2}^2 \distr \int_0^1 \cX_t^2 \, \dd t \qquad
   \text{as \ $n \to \infty$.}
 \end{align}
Let us apply Lemmas \ref{Conv2Funct} and \ref{Marci} with the special choices
 \ $d := p := q := 1,$ \ $h : \RR \to \RR$, \ $h(x) := x$, \ $x \in \RR$,
 \ $K : [0,1] \times \RR^2 \to \RR$,
  \[
    K(s, x_1, x_2) := x_1^2, \qquad (s, x_1, x_2) \in [0, 1] \times \RR^2 ,
  \]
 and \ $\cU := \cX$, \ $\cU^{(n)} := \cX^{(n)}$, \ $n \in \NN$.
\ Then
 \begin{align*}
  | K(s, x_1, x_2) - K(t, y_1, y_2) |
    & = | x_1^2 - y_1^2|
      \leq (| x_1 | + | y_1 |) | x_1 - y_1|
      \leq 2 R(| t-s| + | x_1 - y_1 |)\\
    & \leq 2R\big( | t - s| + \| (x_1, x_2) - (y_1, y_2) \| \big)
 \end{align*}
 for all \ $s, t \in [0, 1]$ \ and \ $(x_1, x_2), (y_1, y_2) \in \RR^2$ \ with
 \ $\| (x_1, x_2) \| \leq R$ \ and  \ $\| (y_1, y_2) \| \leq R$, \ where
 \ $R > 0$.
\ Further, using the definition of \ $\Phi$ \ and \ $\Phi_n$, \ $n \in \NN$,
 \ given in Lemma \ref{Marci},
 \begin{align*}
  \Phi_n(\cX^{(n)})
  &= \left(\cX_1^{(n)},\frac{1}{n}\sum_{k=1}^n\bigl(\cX_{k/n}^{(n)}\bigr)^2 \right)
   = \left(\frac{1}{n} X_n,\frac{1}{n^3} \sum_{k=1}^n X_k^2\right) , \\
  \Phi(\cX)
  &= \left(\cX_1,\int_0^1\cX_u^2\,\dd u\right) .
 \end{align*}
Since the process \ $(\cX_t)_{t\in\RR_+}$ \ admits continuous paths with
 probability one, Lemma \ref{Conv2Funct}
 (with the choice \ $C := \CC(\RR_+,\RR)$) \ and Lemma \ref{Marci} yield
 \eqref{seged2}.
Since \ $\mu > 0$, \ by the SDE \eqref{X}, we have
 \ $\PP\bigl( \cX_t =0, \, t \in [0, 1] \bigr) = 0$, \ which implies that
 \ $\PP\bigl( \int_0^1 \cX_t^2 \, \dd t > 0 \bigr) = 1$.
\ Consequently, the distribution function of \ $\int_0^1 \cX_t^2 \, \dd t$ \ is
 continuous at 0, and hence, by \eqref{seged2},
 \[
   \PP\left( \sum_{k=1}^n X_{k-2}^2 > 0 \right)
   = \PP\left( \frac{1}{n^3} \sum_{k=1}^n X_{k-2}^2 > 0 \right)
     \to \PP\left( \int_0^1 \cX_t^2 \, \dd t > 0 \right) =1
     \qquad \text{as \ $n \to \infty$.}
 \]
Clearly, \eqref{CLSE} also holds, hence we obtain the statement in the case of
 \ $(\alpha, \beta) \in (0, 1)^2$.

Next we consider the case of \ $(\alpha, \beta) = (1, 0)$.
\ In this case equation \eqref{INAR2} has the form \ $X_k = X_{k-1} + \vare_k$,
 \ $k \in \NN$, \ and hence \ $X_n = \sum_{k=1}^n \vare_k$, \ $n \in \NN$.
\ By the strong law of large numbers we have
 \begin{equation}\label{SLLN1}
  n^{-1} X_n \as \mu ,
 \end{equation}
 and hence
 \[
   n^{-2} X_n^2 \as \mu^2 ,
 \]
 where \ $\as$ \ denotes almost sure convergence.
Then, by Toeplitz theorem, we conclude
 \begin{equation}\label{SLLN2}
  n^{-3} \sum_{k=1}^n X_k^2 = \sum_{k=1}^n \frac{k^2}{n^3} k^{-2}X_k^2 \as \frac{1}{3} \mu^2 ,
 \end{equation}
 where we used that 
 \begin{align*}
  \lim_{n\to\infty} \sum_{k=1}^n  \frac{k^2}{n^3} 
   = \lim_{n\to\infty} \frac{n(n+1)(2n+1)}{6n^3} 
   =\frac{1}{3}.
 \end{align*}

Since \ $\mu > 0$, \ this implies the existence of an event
 \ $\Omega_0 \in \cA$ \ such that \ $\PP(\Omega_0) = 1$, \ and for all
 \ $\omega\in\Omega_0$ \ there exists an \ $n_0(\omega) \in \NN$ \ such that
 \ $\sum_{k=1}^n X_{k-2}(\omega)^2 > 0$ \ for \ $n \geq n_0(\omega)$.
\ This is equivalent with
 \ $\PP(\bigcup_{n=1}^\infty \{ \sum_{k=1}^n X_{k-2}^2 > 0 \}) = 1$,
 \ and, by continuity of probability, is also equivalent with
 \ $\lim_{n\to\infty} \PP(\{\sum_{k=1}^n X_{k-2}^2 > 0\}) = 1$.
\ Clearly \eqref{CLSE} also holds, hence we obtain the statement in case
 \ $(\alpha, \beta) = (1, 0)$.

Finally, we consider the case \ $(\alpha, \beta) = (0, 1)$.
\ In this case equation \eqref{INAR2} has the form \ $X_k = X_{k-2} + \vare_k$,
 \ $k \in \NN$, \ and hence \ $X_{2n} = \sum_{k=1}^n \vare_{2k}$,
 \ $X_{2n-1} = \sum_{k=1}^n \vare_{2k-1}$, \ $n \in \NN$.
\ By the strong law of large numbers we have
 \[
   n^{-1} X_{2n} \as \mu , \qquad n^{-1} X_{2n-1} \as \mu ,
 \]
 which yield that
 \begin{equation}\label{SLLN1_01}
  n^{-1} X_n \as \frac{1}{2} \mu .
 \end{equation}
Using Toeplitz theorem, as in case \ $(\alpha, \beta) = (1,0)$, \ we get
 \begin{equation}\label{SLLN2_01}
  n^{-3} \sum_{k=1}^n X_k^2 \as \frac{1}{12} \mu^2 .
 \end{equation}
One can finish the proof as in case \ $(\alpha, \beta) = (1, 0)$.
\proofend

The recursion \eqref{regr} can also be written in the form
 \begin{equation}\label{can}
  X_k = \varrho X_{k-1} - \beta (X_{k-1}-X_{k-2}) + M_k + \mu ,
  \qquad k \in \NN .
 \end{equation}
The representation \eqref{can} is called the canonical form of Sims, Stock and
 Watson \cite{SimStoWat}, see also Hamilton \cite[17.7.6]{Ham}.
A natural CLS estimator of the stability parameter \ $\varrho$ \ takes the
 form \ $\hvarrho_n(\bX_n) = \halpha_n(\bX_n) + \hbeta_n(\bX_n)$, \ since, for
 each \ $n \in \NN$, \ a CLS estimator
 \ $(\hvarrho_n(\bX_n), \hbeta_n(\bX_n),\hmu_n(\bX_n))$ \ of \ $(\varrho, \beta, \mu)$ \ based on
 a sample \ $\bX_n = (X_1, \ldots, X_n)$ \ can be obtained by minimizing the
 sum of squares
 \begin{equation}\label{SumSquares_varrho}
   \sum_{k=1}^n \big( X_k - \EE(X_k \mid \cF_{k-1}) \big)^2
   = \sum_{k=1}^n
      \big(X_k - \varrho X_{k-1} + \beta (X_{k-1} - X_{k-2}) - \mu\big)^2
 \end{equation}
 with respect to
 \ $(\varrho, \beta, \mu)$ \ over \ $\RR^3$.
\ One can easily argue that any CLS estimator
 \ $(\hvarrho_n, \hbeta_n, \hmu_n) : \RR^n \to \RR^3$ \ of
 \ $(\varrho, \beta, \mu)$ \ is of the form
 \begin{align}\label{CLSE_rho_beta}
   \begin{bmatrix}
    \hvarrho_n(\bx_n) \\
    \hbeta_n(\bx_n) \\
    \hmu_n(\bx_n)
   \end{bmatrix}
   = \begin{bmatrix} 1 & 1 & 0 \\ 0 & 1 & 0 \\ 0 & 0 & 1 \end{bmatrix}
     \begin{bmatrix}
      \halpha_n(\bx_n) \\
      \hbeta_n(\bx_n) \\
      \hmu_n(\bx_n)
     \end{bmatrix} ,
   \qquad \bx_n \in \RR^n ,
 \end{align}
 where \ $(\halpha_n, \hbeta_n, \hmu_n)$ \ is a CLS estimator of
 \ $(\alpha, \beta, \mu)$.
\ Namely, if \ $\psi : \RR^3 \to \RR^3$ \ is a bijective measurable function
 such that
 \[
   \RR^3 \ni (\alpha', \beta', \mu') \mapsto \psi(\alpha', \beta', \mu')
   := \begin{bmatrix}
       \alpha' + \beta'\\
       h(\alpha', \beta', \mu') \\
      \end{bmatrix}
   =: \begin{bmatrix}
       \varrho' \\
       \gamma' \\
       \delta'
      \end{bmatrix}
 \]
 with some function \ $h : \RR^3 \to \RR^2$, \ then there is a bijection
 between the set of CLS estimators of the parameters \ $(\alpha, \beta, \mu)$
 \ and the set of CLS estimators of the parameters
 \ $\psi(\alpha, \beta, \mu)$.
\ Indeed, for all \ $n \in \NN$, \ $(x_1, \ldots, x_n) \in \RR^n$ \ and
 \ $(\alpha', \beta', \mu') \in \RR^3$,
 \begin{align*}
  \sum_{k=1}^n (x_k - \alpha' x_{k-1} - \beta' x_{k-2} - \mu')^2
  & = \sum_{k=1}^n
       \left(x_k - \begin{bmatrix}
                    \alpha' \\
                    \beta' \\
                    \mu'
                   \end{bmatrix}^\top
                   \begin{bmatrix}
                    x_{k-1} \\
                    x_{k-2} \\
                    1
                   \end{bmatrix} \right)^2  \\
  & = \sum_{k=1}^n
       \left(x_k - \left(\psi^{-1}(\varrho', \gamma', \delta')\right)^\top
                   \begin{bmatrix}
                    x_{k-1} \\
                    x_{k-2} \\
                    1
                   \end{bmatrix}\right)^2 ,
 \end{align*}
 hence \ $(\halpha_n, \hbeta_n, \hmu_n) : \RR^n \to \RR^3$ \ is a CLS estimator
 of \ $(\alpha, \beta, \mu)$ \ if and only if
 \ $\psi(\halpha_n, \hbeta_n, \hmu_n)$ \ is a CLS estimator of
 \ $\psi(\alpha, \beta, \mu)$.
\ With the special choice \ $h : \RR^3 \to \RR^2$,
 \ $h(\alpha, \beta, \mu) := (\beta, \mu)$,
 \ $(\alpha, \beta, \mu) \in \RR^3$, \ we get \eqref{CLSE_rho_beta}.
In what follows, by speaking about the CLS estimator \ $\hvarrho_n$ \ of
 \ $\varrho$ \ we mean the first coordinate of \ $\psi(\halpha_n,\hbeta_n)$.
\ Hence, by Proposition \ref{ExUn}, the probability of the
 existence of a unique CLS estimator
 \ $(\hvarrho_n(\bX_n), \hbeta_n(\bX_n), \hmu_n(\bX_n))$ \ converges to 1 as
 \ $n \to \infty$, \ and this CLS estimator has the form
 \begin{equation}\label{CLSE_rho}
  \begin{bmatrix}
   \hvarrho_n(\bX_n) \\
   \hbeta_n(\bX_n) \\
   \hmu_n(\bX_n)
  \end{bmatrix}
  = \bA_n^{-1} \bb_n
 \end{equation}
 on the set \ $\{ \omega \in \Omega : \sum_{k=1}^n X_{k-2}(\omega)^2 > 0\}$,
 \ where
 \[
   \bA_n := \sum_{k=1}^n
             \begin{bmatrix}
              X_{k-1}^2 & - X_{k-1} V_{k-1} & X_{k-1} \\
              - X_{k-1} V_{k-1} & V_{k-1}^2 & - V_{k-1} \\
              X_{k-1} & - V_{k-1} & 1
             \end{bmatrix} , \qquad
   \bb_n := \sum_{k=1}^n
             \begin{bmatrix}
              X_k X_{k-1} \\
              - X_k V_{k-1} \\
              X_k
             \end{bmatrix}
 \]
 with
 \[
   V_{k-1} := X_{k-1} - X_{k-2} , \qquad k \in \NN .
 \]
(In Appendix \ref{app_A}, in Remark \ref{REM_Putzer} one can find a detailed
 motivation of the definition of \ $V_k$, \ $k \in \NN$.)
\ Indeed, by \eqref{CLSE},
 \begin{align*}
  \begin{bmatrix} \hvarrho_n(\bX_n) \\ \hbeta_n(\bX_n) \\  \hmu_n(\bX_n) \end{bmatrix}
  = \begin{bmatrix} 1 & 1 & 0 \\ 0 & 1 & 0 \\ 0 & 0 & 1 \end{bmatrix}
     \bF_n^{-1} \bg_n
   = \left( \begin{bmatrix} 1 & 0 & 0 \\ -1 & 1 & 0 \\ 0 & 0 & 1 \end{bmatrix}
            \bF_n
            \begin{bmatrix}
             1 & 1 & 0 \\
             0 & 1 & 0 \\
             0 & 0 & 1
            \end{bmatrix}^{-1}\right)^{-1}
     \begin{bmatrix} 1 & 0 & 0 \\ -1 & 1 & 0 \\ 0 & 0 & 1 \end{bmatrix} \bg_n
   =  \bA_n^{-1} \bb_n,
 \end{align*}
 which also shows that \ $\bA_n^{-1}$ \ exists on the set
 \ $\{ \omega \in \Omega : \sum_{k=1}^n X_{k-2}(\omega)^2 > 0\}$.

Alternatively, the CLS estimator  $\hvarrho_n(\bX_n)$ \ of the stability
 parameter \ $\varrho$ \ could also be obtained via a CLS estimator
 \ $(\halpha_n(\bX_n), \hvarrho_n(\bX_n), \hmu_n(\bX_n))$ \ of
 \ $(\alpha, \varrho, \mu)$.

Note also that in case of an unstable INAR(2) process, i.e., when
 \ $\varrho = 1$, \  we have
 \begin{equation}\label{rec_V}
  V_k = - \beta V_{k-1} + M_k + \mu , \qquad k \in \NN ,
 \end{equation}
 hence \ $(V_k)_{k\in\NN}$ \ is a stable AR(1) process with heteroscedastic
 innovations \ $(M_k)_{k\in\NN}$ \ and with positive drift \ $\mu$
 \ whenever \ $0\leq \beta < 1$.

\section{Proof of the main results}
\label{section_proof}

In case of an unstable INAR(2) process, i.e., when \ $\varrho=\alpha+\beta=1$,
by \eqref{CLSE_rho}, we have
 \begin{equation}\label{CLSE_rho-1}
   \begin{bmatrix}
    \hvarrho_n(\bX_n) - 1 \\
    \hbeta_n(\bX_n) - \beta \\
    \hmu_n(\bX_n) - \mu
   \end{bmatrix} = \bA_n^{-1} \bd_n ,
   \qquad n \in \NN ,
 \end{equation}
 on the set \ $\{ \omega \in \Omega : \sum_{k=1}^n X_{k-2}(\omega)^2 > 0\}$,
 \ where
 \[
   \bd_n
   := \sum_{k=1}^n
       \begin{bmatrix} M_k X_{k-1} \\ - M_k V_{k-1} \\ M_k \end{bmatrix} , \qquad
   n \in \NN .
 \]
Theorems \ref{main}, \ref{10main}, and \ref{01main} will follow from Theorems
 \ref{main_Ad}, \ref{10main_Ad}, and \ref{01main_Ad}, respectively
 (see the details below).

\begin{Thm}\label{main_Ad}
Under the assumptions of Theorem \ref{main} we have
 \ $(\tbA_n , \tbd_n) \distr (\tbA, \tbd)$ \ as \ $n \to \infty$, \ where
 \begin{equation*}%\label{inv}
  \tbA_n := \begin{bmatrix}
             n^{-3/2} & 0 & 0 \\
             0 & n^{-1} & 0 \\
             0 & 0 & n^{-1/2}
            \end{bmatrix}
            \bA_n
            \begin{bmatrix}
             n^{-3/2} & 0 & 0 \\
             0 & n^{-1} & 0 \\
             0 & 0 & n^{-1/2}
            \end{bmatrix} , \qquad
  \tbd_n := \begin{bmatrix}
             n^{-2} & 0 & 0 \\
             0 & n^{-3/2} & 0 \\
             0 & 0 & n^{-1}
            \end{bmatrix} \bd_n,
 \end{equation*}
 \[
   \tbA :=
   \begin{bmatrix}
    \int_0^1 \cX_t^2 \, \dd t & 0 & \int_0^1 \cX_t \, \dd t \\
    0 & \frac{2\beta}{1 + \beta} \int_0^1 \cX_t \, \dd t & 0 \\
    \int_0^1 \cX_t \, \dd t & 0 & 1
   \end{bmatrix} , \qquad
   \tbd :=
   \begin{bmatrix}
    \sqrt{2\alpha\beta} \int_0^1 \cX_t^{3/2} \, \dd \cW_t \\[2mm]
    - \frac{2\beta\sqrt{\alpha}}{\sqrt{1+\beta}}
      \int_0^1 \cX_t \, \dd \tcW_t \\[2mm]
    (1 + \beta) \cX_1 - \mu
   \end{bmatrix} ,
 \]
 where \ $(\cW_t)_{t\in\RR_+}$ \ and \ $(\tcW_t)_{t\in\RR_+}$ \ are independent
 standard Wiener processes.
\end{Thm}

\begin{Thm}\label{10main_Ad}
Under the assumptions of Theorem \ref{10main} we have
 \ $(\tbA_n , \tbd_n) \distr (\tbA, \tbd)$ \ as \ $n \to \infty$, \ where
 \begin{equation*}%\label{inv}
  \tbA_n := \begin{bmatrix}
             n^{-3/2} & 0 & 0 \\
             0 & n^{-1/2} & 0 \\
             0 & 0 & n^{-1/2}
            \end{bmatrix}
            \bA_n
            \begin{bmatrix}
             n^{-3/2} & 0 & 0 \\
             0 & n^{-1/2} & 0 \\
             0 & 0 & n^{-1/2}
            \end{bmatrix} , \quad
  \tbd_n := \begin{bmatrix}
             n^{-3/2} & 0 & 0 \\
             0 & n^{-1/2} & 0 \\
             0 & 0 & n^{-1/2}
            \end{bmatrix}
            \bd_n ,
 \end{equation*}
 \[
   \tbA :=
   \begin{bmatrix}
    \frac{1}{3} \mu^2 & - \frac{1}{2} \mu^2 &  \frac{1}{2} \mu \\
    - \frac{1}{2} \mu^2 & \mu^2 + \sigma^2 & - \mu \\
    \frac{1}{2} \mu & - \mu & 1
   \end{bmatrix} , \qquad
   \tbd \distre
   \cN_3(\bzero, \sigma^2 \tbA) .
 \]
\end{Thm}

\begin{Thm}\label{01main_Ad}
Under the assumptions of Theorem \ref{01main} we have
 \ $(\tbA_n , \tbd_n) \distr (\tbA, \tbd)$ \ as \ $n \to \infty$, \ where
 \begin{equation*}%\label{inv}
  \tbA_n := \begin{bmatrix}
             n^{-3/2} & 0 & 0 \\
             0 & n^{-1} & 0 \\
             0 & 0 & n^{-1/2}
            \end{bmatrix}
            \bA_n
            \begin{bmatrix}
             n^{-3/2} & 0 & 0 \\
             0 & n^{-1} & 0 \\
             0 & 0 & n^{-1/2}
            \end{bmatrix} , \qquad
  \tbd_n := \begin{bmatrix}
             n^{-3/2} & 0 & 0 \\
             0 & n^{-1} & 0 \\
             0 & 0 & n^{-1/2}
            \end{bmatrix}
            \bd_n ,
 \end{equation*}
 \[
   \tbA :=
   \begin{bmatrix}
    \frac{1}{12} \mu^2 & 0 & \frac{1}{4} \mu \\
    0 & \sigma^2 \int_0^1 \cW_t^2 \, \dd t & 0 \\
    \frac{1}{4} \mu & 0 & 1
   \end{bmatrix} , \qquad
   \tbd :=
   \begin{bmatrix}
    \frac{1}{2} \mu \sigma \int_0^1 t \, \dd \tcW_t \\[1mm]
    \sigma^2 \int_0^1 \cW_t \, \dd \cW_t \\[1mm]
    \sigma \tcW_1
   \end{bmatrix} ,
 \]
 where \ $(\cW_t)_{t\in\RR_+}$ \ and \ $(\tcW_t)_{t\in\RR_+}$ \ are independent
 standard Wiener processes.
\end{Thm}

Now we briefly summarize how Theorem \ref{main_Ad} yields Theorem \ref{main}.
The function \ $g : \RR^{3\times3} \times \RR^{3\times1} \to \RR^{3\times1}$,
 \ defined by
 \begin{align}\label{01seged12}
  g(\bX, \by)
  :=\begin{cases}
     \bX^{-1} \by, & \text{if \ $\exists$ $\bX^{-1}$,}\\
     \bzero, & \text{otherwise,}
    \end{cases}
 \end{align}
 is continuous on the set
 \ $\{\bX \in \RR^{3\times3} : \exists \, \bX^{-1}\} \times \RR^{3\times1}$, \ and
 \ the limit distribution in Theorem \ref{main_Ad} is concentrated on this
 set, since, 
 \[
   \det(\tbA)
    = \frac{2\beta}{1+\beta} \int_0^1 \cX_t\,\dd t 
        \left( \int_0^1 \cX_t^2\,\dd t  - \left(\int_0^1 \cX_t\,\dd t\right)^2 \right),
 \]
 and, by Remark \ref{REMARK1} and the proof of Proposition \ref{ExUn},
 \[
   \PP\left( \int_0^1 \cX_t\,\dd t >0\right)
   = \PP\left( \int_0^1 \cX_t^2\,\dd t >0\right)
   = 1 ,
 \]
 and, by Lemma 4.3 in Barczy et al. \cite{BarDorLiPap},
 \[
   \PP\left( \int_0^1 \cX_t^2 \, \dd t
             - \left(\int_0^1 \cX_t \, \dd t \right)^2 >0 \right)
   = 1 .
 \]
Hence the continuous mapping theorem
 (see, e.g., Theorem 2.3 in van der Vaart \cite{Vaa}) yields that
 \begin{align*}
  g\bigl( \tbA_n , \, \tbd_n \bigr)
  \distr
  g\bigl( \tbA , \, \tbd \bigr)
 \end{align*}
 as \ $n\to\infty$.
Under the conditions of Proposition \ref{ExUn}, by \eqref{seged2} and \eqref{CLSE_rho-1}, we have
 \begin{multline*}
  \PP\left( \begin{bmatrix}
             n & 0 & 0 \\
             0 & n^{1/2} & 0 \\
             0 & 0 & 1
            \end{bmatrix}
            \begin{bmatrix}
             \hvarrho_n(\bX_n) - 1 \\
             \hbeta_n(\bX_n) -\beta \\
             \hmu_n(\bX_n) - \mu
            \end{bmatrix}
            = g\bigl( \tbA_n , \, \tbd_n \bigr) \right)
  \geq \PP\left( \exists \, \tbA_n^{-1} \right)
  = \PP\left( \exists \, \bA_n^{-1} \right) \\
  \geq \PP\left(\sum_{k=1}^n X_{k-2}^2 > 0\right)
  = \PP\left(\frac{1}{n^3} \sum_{k=1}^n X_{k-2}^2 > 0\right)
  \to \PP\left( \int_0^1 \cX_t^2\,\dd t > 0 \right)
  = 1
 \end{multline*}
 as \ $n \to \infty$.
\ Clearly, if \ $\xi_n$, \ $\eta_n$, \ $n \in \NN$, \ and \ $\xi$ \ are random
 variables such that \ $\xi_n\distr \xi$ \ as \ $n \to \infty$ \ and
 \ $\lim_{n\to\infty} \PP(\xi_n = \eta_n) = 1$, \ then \ $\eta_n \distr \xi$ \ as
 \ $n \to \infty$, \ see, e.g., Barczy et al.\ \cite[Lemma 3.1]{BarIspPap1}.
Consequently, under the conditions of Theorem \ref{main}, Theorem
 \ref{main_Ad} yields that
 \[
   \begin{bmatrix}
    n(\hvarrho_n(\bX_n) - 1) \\
    n^{1/2}(\hbeta_n(\bX_n) -\beta) \\
    \hmu_n(\bX_n) - \mu
   \end{bmatrix}
   \distr g\bigl( \tbA , \, \tbd \bigr) \qquad
   \text{as \ $n \to \infty$,}
 \]
 where
 \begin{multline*}
  g\bigl( \tbA , \, \tbd \bigr)
  = \tbA^{-1}\tbd \\
  \begin{aligned}
   & =\frac{1}
           {\int_0^1 \cX_t^2 \, \dd t - \bigl(\int_0^1 \cX_t \, \dd t\bigr)^2}
      \begin{bmatrix}
       1 & 0 & - \int_0^1 \cX_t \, \dd t \\
       0 & \frac{1+\beta}{2\beta}
           \frac{\int_0^1 \cX_t^2 \, \dd t
                 - \left( \int_0^1 \cX_t \, \dd t \right)^2}
                {\int_0^1 \cX_t\,\dd t}
         & 0 \\
       - \int_0^1 \cX_t \, \dd t & 0 & \int_0^1 \cX_t^2 \, \dd t
      \end{bmatrix}\!\!\!
      \begin{bmatrix}
       \sqrt{2\alpha\beta} \int_0^1 \cX_t^{3/2} \, \dd\cW_t \\
       -\frac{2\beta\sqrt{\alpha}}{\sqrt{1+\beta}}
        \int_0^1 \cX_t \, \dd\tcW_t \\
       (1 + \beta) \cX_1 - \mu
      \end{bmatrix} \\
   &= \begin{bmatrix}
       \frac{\sqrt{2\alpha\beta} \int_0^1 \cX_t^{3/2} \, \dd\cW_t
             - [(1 + \beta) \cX_1 - \mu] \int_0^1 \cX_t \, \dd t}
            {\int_0^1 \cX_t^2 \, \dd t
             - \bigl(\int_0^1 \cX_t \, \dd t\bigr)^2} \\[5mm]
       -\frac{\sqrt{\alpha(1+\beta)} 
              \int_0^1\cX_t \, \dd\tcW_t}
             {\int_0^1 \cX_t \, \dd t} \\[3mm]
       \frac{-\sqrt{2\alpha\beta} \int_0^1 \cX_t \, \dd t
             \int_0^1 \cX_t^{3/2} \, \dd\cW_t
             + [(1 + \beta) \cX_1 - \mu] \int_0^1 \cX_t^2 \, \dd t}
            {\int_0^1 \cX_t^2 \, \dd t
             - \bigl(\int_0^1 \cX_t \, \dd t\bigr)^2}
      \end{bmatrix} .
  \end{aligned}
 \end{multline*}
Hence we obtain \eqref{rho} and \eqref{mu}, and, using that
 \ $\int_0^1 \cX_t \, \dd \tcW_t$ \ is symmetric, also the convergence of the
 second coordinate in \eqref{alpha,beta}.
By Slutsky's lemma, convergence \eqref{rho} implies
 \ $n^{1/2}(\hvarrho_n(\bX_n) - 1) \stoch 0$ \ as well, where \ $\stoch$
 \ denotes convergence in probability, hence
 \[
   \begin{bmatrix}
    n^{1/2} (\halpha_n(\bX_n) -\alpha) \\
    n^{1/2} (\hbeta_n(\bX_n) -\beta) \\
   \end{bmatrix}
   = n^{1/2} (\hbeta_n(\bX_n) -\beta)
     \begin{bmatrix} - 1 \\ 1 \end{bmatrix}
     + n^{1/2}(\hvarrho_n(\bX_n) - 1)
       \begin{bmatrix} 1 \\ 0 \end{bmatrix}
 \]
 yields \eqref{alpha,beta}.

Next we briefly summarize how Theorem \ref{10main_Ad} yields Theorem
 \ref{10main}.
Similarly to the previous case, under the conditions of Proposition
 \ref{ExUn}, by \eqref{seged2} and \eqref{CLSE_rho-1}, we have
 \begin{align*}
  \PP\left( \begin{bmatrix}
             n^{3/2} & 0 & 0 \\
             0 & n^{1/2} & 0 \\
             0 & 0 & n^{1/2}
            \end{bmatrix}
            \begin{bmatrix}
             \hvarrho_n(\bX_n) - 1 \\
             \hbeta_n(\bX_n) \\
             \hmu_n(\bX_n) - \mu
            \end{bmatrix}
            = g\bigl( \tbA_n , \, \tbd_n \bigr) \right)
  \geq \PP\left(\sum_{k=1}^n X_{k-2}^2 > 0\right)
  \to 1
 \end{align*}
 as \ $n\to \infty$.
Consequently, under the conditions of Theorem \ref{10main}, Theorem
 \ref{10main_Ad} yields that
 \[
   \begin{bmatrix}
    n^{3/2} (\hvarrho_n(\bX_n) - 1) \\
    n^{1/2} \hbeta_n(\bX_n) \\
    n^{1/2} (\hmu_n(\bX_n) - \mu)
   \end{bmatrix}
   \distr g\bigl( \tbA , \, \tbd \bigr) \qquad
   \text{as \ $n \to \infty$,}
 \]
 where
 \begin{align*}
  g\bigl( \tbA , \, \tbd \bigr)
  &= \tbA^{-1} \tbd
   \distre \tbA^{-1} \cN_3(\bzero, \sigma^2\tbA)
   \distre \cN_3\bigl(\bzero, \sigma^2 \tbA^{-1}\bigr) \\
  &\distre
   \cN_3\left( \begin{bmatrix} 0 \\ 0 \\ 0 \end{bmatrix} , \,
             \frac{1}{\mu^2}
             \begin{bmatrix}
              12 \sigma^2 & 0 & -6 \mu \sigma^2 \\
              0 & \mu^2 & \mu^3 \\
              -6 \mu \sigma^2 & \mu^3 & \mu^2 (\mu^2 + 4 \sigma^2)
             \end{bmatrix} \right) .
 \end{align*}
Hence we obtain \eqref{10rho}, \eqref{10mu}, and convergence of the second coordinate in
 \eqref{10,alpha,beta}.
By Slutsky's lemma, convergence \eqref{10rho} implies
 \ $n^{1/2} (\hvarrho_n(\bX_n) - 1) \stoch 0$ \ as well, hence
 \[
   \begin{bmatrix}
    n^{1/2} (\halpha_n(\bX_n) - 1) \\
    n^{1/2} \hbeta_n(\bX_n)  \\
   \end{bmatrix}
   = n^{1/2} \hbeta_n(\bX_n)
     \begin{bmatrix} - 1 \\ 1 \end{bmatrix}
     + n^{1/2}(\hvarrho_n(\bX_n) - 1)
       \begin{bmatrix} 1 \\ 0 \end{bmatrix}
 \]
 yields \eqref{10,alpha,beta}.

Finally, we briefly summarize how Theorem \ref{01main_Ad} yields Theorem
 \ref{01main}.
Similarly as above, under the conditions of Proposition \ref{ExUn}, by
 \eqref{seged2} and \eqref{CLSE_rho-1}, we have
 \begin{align*}
  \PP\left( \begin{bmatrix}
             n^{3/2} & 0 & 0 \\
             0 & n & 0 \\
             0 & 0 & n^{1/2}
            \end{bmatrix}
            \begin{bmatrix}
             \hvarrho_n(\bX_n) - 1 \\
             \hbeta_n(\bX_n) - 1 \\
             \hmu_n(\bX_n) - \mu
            \end{bmatrix}
            = g\bigl( \tbA_n , \, \tbd_n \bigr) \right)
  \geq \PP\left(\sum_{k=1}^n X_{k-2}^2 > 0\right)
  \to 1
 \end{align*}
 as \ $n \to \infty$.
\ Consequently, under the conditions of Theorem \ref{01main}, Theorem
 \ref{01main_Ad} yields that
 \[
   \begin{bmatrix}
    n^{3/2} (\hvarrho_n(\bX_n) - 1) \\
    n (\hbeta_n(\bX_n) - 1) \\
    n^{1/2} (\hmu_n(\bX_n) - \mu)
   \end{bmatrix}
   \distr g\bigl( \tbA , \, \tbd \bigr) \qquad
   \text{as \ $n \to \infty$,}
 \]
 where
 \begin{align*}
  g\bigl( \tbA , \, \tbd \bigr)
  &= \tbA^{-1} \tbd
   = \begin{bmatrix}
      \frac{48}{\mu^2} & 0 & - \frac{12}{\mu} \\
      0 & \frac{1}{\sigma^2 \int_0^1\cW_t^2\,\dd t} & 0 \\
      - \frac{12}{\mu} & 0 & 4
     \end{bmatrix}
     \begin{bmatrix}
      \frac{1}{2} \mu \sigma \int_0^1 t \, \dd \tcW_t \\
      \sigma^2 \int_0^1 \cW_t \, \dd \cW_t \\
      \sigma \tcW_1
     \end{bmatrix} \\
  &= \begin{bmatrix}
      \frac{24 \sigma}{\mu} \int_0^1 t \, \dd \tcW_t
       - \frac{12 \sigma}{\mu} \tcW_1 \\[1mm]
      \frac{\int_0^1 \cW_t \, \dd \cW_t}{\int_0^1 \cW_t^2 \, \dd t} \\
      - 6 \sigma \int_0^1 t \, \dd \tcW_t + 4 \sigma \tcW_1
     \end{bmatrix}
   = \begin{bmatrix}
      \frac{12 \sigma}{\mu} \int_0^1 (2t - 1) \, \dd \tcW_t \\[1mm]
      \frac{\int_0^1 \cW_t \, \dd \cW_t}{\int_0^1 \cW_t^2 \, \dd t} \\
      2 \sigma \int_0^1 (2 -3t) \, \dd \tcW_t
     \end{bmatrix} .
 \end{align*}
Since \ $\int_0^1 (2t -1) \, \dd \tcW_t$ \ and
 \ $\int_0^1 (2 -3t) \, \dd \tcW_t$ \ are normally distributed random variables
 with mean \ $0$ \ and with variance
 \[
   \int_0^1 (2t - 1)^2 \, \dd t = \frac{1}{3} , \qquad
   \int_0^1 (2 -3t)^2 \, \dd t = 1 ,
 \]
 respectively, we obtain \eqref{01rho}, \eqref{01mu}, and convergence of the second
 coordinate in \eqref{01,alpha,beta}.
By Slutsky's lemma, convergence \eqref{01rho} implies
 \ $n (\hvarrho_n(\bX_n) - 1) \stoch 0$ \ as well, hence
 \[
   \begin{bmatrix}
    n \halpha_n(\bX_n) \\
    n (\hbeta_n(\bX_n) -1)  \\
   \end{bmatrix}
   = n (\hbeta_n(\bX_n) -1)
     \begin{bmatrix} - 1 \\ 1 \end{bmatrix}
     + n (\hvarrho_n(\bX_n) - 1)
       \begin{bmatrix} 1 \\ 0 \end{bmatrix}
 \]
 yields \eqref{01,alpha,beta}.

\section{Proof of Theorem \ref{main_Ad}}
\label{section_proof_main}

We have
 \begin{equation}\label{tbA_tbd}
  \begin{split}
   \tbA_n&= \sum_{k=1}^n
             \begin{bmatrix}
              n^{-3} X_{k-1}^2 & - n^{-5/2} X_{k-1} V_{k-1} & n^{-2} X_{k-1} \\
              - n^{-5/2} X_{k-1} V_{k-1} & n^{-2} V_{k-1}^2 & - n^{-3/2} V_{k-1} \\
              n^{-2} X_{k-1} & - n^{-3/2} V_{k-1} & n^{-1}
             \end{bmatrix} , \\
   \tbd_n&= \sum_{k=1}^n
      \begin{bmatrix}
       n^{-2} M_k X_{k-1} \\
       - n^{-3/2} M_k V_{k-1} \\
       n^{-1} M_k
      \end{bmatrix} .
  \end{split}
 \end{equation}

\begin{Lem}\label{XV_main_VV}
Under the assumptions of Theorem \ref{main} we have
 \begin{gather}\label{V}
  n^{-3/2} \sum_{k=1}^n V_k \stoch 0 \qquad \text{as \ $n \to \infty$,} \\
  \label{XV}
  n^{-5/2} \sum_{k=1}^n X_k V_k \stoch 0 \qquad \text{as \ $n \to \infty$,} \\
  n^{-2} \left( \sum_{k=1}^n V_k^2
               - \frac{2\beta}{1+\beta} \sum_{k=1}^n X_{k-1} \right)
  \stoch 0 \qquad \text{as \ $n \to \infty$.} \label{main_VV}
 \end{gather}
\end{Lem}

\noindent
\textbf{Proof.}
We have \ $\sum_{k=1}^n V_k = X_n \geq 0$, $n\in\NN$, \ and, by Corollary
 \ref{EEX_EEU_EEV}, \ $\EE(X_n) = \OO(n)$, \ hence we conclude \eqref{V}.
We have
 \begin{align*}
  \sum_{k=1}^n (X_k-X_{k-1})^2
  &= \sum_{k=1}^n X_k^2 - 2 \sum_{k=1}^n X_k X_{k-1} + \sum_{k=1}^n X_{k-1}^2
   = 2\sum_{k=1}^n X_k^2 - 2 \sum_{k=1}^n X_k X_{k-1} - X_n^2 \\
  &= 2\sum_{k=1}^n X_k (X_k - X_{k-1}) - X_n^2 ,
 \end{align*}
 thus
 \begin{align}\label{XV_VV}
   \sum_{k=1}^n X_k V_k
   = \frac{1}{2} X_n^2 + \frac{1}{2} \sum_{k=1}^n V_k^2
   \geq 0 .
 \end{align}
Corollary \ref{EEX_EEU_EEV} implies
 \[
   \EE\left(\sum_{k=1}^n X_k V_k \right)
   = \frac{1}{2} \EE(X_n^2) + \frac{1}{2} \sum_{k=1}^n \EE(V_k^2)
   = \OO(n^2) ,
 \]
 hence we obtain \eqref{XV}.

In order to prove \eqref{main_VV} we derive a decomposition of
 \ $\sum_{k=1}^n V_k^2$ \ as a sum of a martingale and some negligible terms.
Using recursion \eqref{rec_V} and Lemma \ref{Moments}, we obtain
 \begin{align*}
  \EE(V_k^2 \mid \cF_{k-1})
  &= \EE\big((- \beta V_{k-1} + M_k + \mu)^2 \mid \cF_{k-1}\big) \\
  &= \beta^2 V_{k-1}^2 - 2 \beta \mu V_{k-1} + \mu^2
     + \EE(M_k^2 \mid \cF_{k-1}) \\
  &= \beta^2 V_{k-1}^2 - 2 \beta \mu V_{k-1} + \mu^2
     + \alpha \beta (X_{k-1} + X_{k-2}) + \sigma^2 \\
  &= \beta^2 V_{k-1}^2 + 2 \alpha \beta X_{k-1} + \mu^2 + \sigma^2
     - (2 \beta \mu + \alpha \beta ) V_{k-1} \\
  &= \beta^2 V_{k-1}^2 + 2 \alpha \beta X_{k-1}
     + \text{constant}
     + \text{constant $\times$ $V_{k-1}$,}
 \end{align*}
 where we used that \ $X_{k-1} + X_{k-2} = 2 X_{k-1} - V_{k-1}$, \ $k \in \NN$.
\ Thus
 \begin{align*}
  &\sum_{k=1}^n V_k^2
   = \sum_{k=1}^n \big[ V_k^2 - \EE(V_k^2 \mid \cF_{k-1}) \big]
     + \sum_{k=1}^n \EE(V_k^2 \mid \cF_{k-1}) \\
  &= \sum_{k=1}^n \big[ V_k^2 - \EE(V_k^2 \mid \cF_{k-1}) \big]
     + \beta^2 \sum_{k=1}^n V_{k-1}^2
     + 2 \alpha \beta \sum_{k=1}^n X_{k-1}
     + \OO(n)
     + \text{constant $\times$ $\sum_{k=1}^n V_{k-1}$.}
 \end{align*}
Consequently,
 \begin{align}\label{sum_Vk2}
  \begin{split}
  \sum_{k=1}^n V_k^2
   &= \frac{1}{1 - \beta^2}
      \sum_{k=1}^n \big[ V_k^2 - \EE(V_k^2 \mid \cF_{k-1}) \big]
      + \frac{2 \beta}{1 + \beta} \sum_{k=1}^n X_{k-1} \\
   &\quad
      - \frac{\beta^2}{1 - \beta^2} V_n^2 + \OO(n)
      + \text{constant $\times$ $\sum_{k=1}^n V_{k-1}$.}
  \end{split}
 \end{align}
By the tower rule of conditional expectation,
 \ $V_k^2 - \EE(V_k^2 \mid \cF_{k-1})$ \ and
 \ $V_\ell^2 - \EE(V_\ell^2 \mid \cF_{\ell-1})$ \ are uncorrelated if
 \ $k \ne \ell$, \ so
 \[
   \EE\left(\left(\sum_{k=1}^n \big[ V_k^2 - \EE(V_k^2 \mid \cF_{k-1}) \big]
            \right)^2\right)
   = \sum_{k=1}^n
      \EE\left( \big[ V_k^2 - \EE(V_k^2 \mid \cF_{k-1}) \big]^2 \right)
   \leq \sum_{k=1}^n \EE( V_k^4 )
   = \OO(n^3) ,
 \]
 where we also used Corollary \ref{EEX_EEU_EEV} and
 \begin{align}\label{cond_var}
  \EE\left(\big[\xi - \EE(\xi \mid \cF)\big]^2 \right)
    = \EE(\xi^2) - \EE\left(\EE(\xi \mid \cF)^2\right) \leq \EE(\xi^2)
 \end{align}
 for an arbitrary random variable \ $\xi$ \ with \ $\EE(\xi^2) < \infty$
 \ and $\sigma$-algebra \ $\cF \subset \cA$.
\ Hence
 \begin{align*}
  \frac{1}{n^2}
  \sum_{k=1}^n \big[ V_k^2 - \EE(V_k^2 \mid \cF_{k-1}) \big]
  \stoch 0 \qquad \text{as \ $n \to \infty$.}
 \end{align*}
We note that this convergence follows also by \eqref{seged_UV_UNIFORM4} with
 the choice \ $(\ell, i, j) = (8, 0, 2)$.
\ Again, by Corollary \ref{EEX_EEU_EEV}, we obtain \ $\EE(V_n^2) = \OO(n)$
 \ and \ $\EE(X_{n-1}^2)=\OO(n^2)$, \ and since \ $\sum_{k=1}^n V_{k-1} = X_{n-1}$,
 \ $n \in \NN$, \ we get \ $n^{-2} V_n^2 \stoch 0$ \ and
 \ $n^{-2} \sum_{k=1}^n V_{k-1} \stoch 0$ \ as \ $n \to \infty$ \ (we note that
 the second convergence follows also by \eqref{seged_UV_UNIFORM1} with the
 choice \ $(\ell, i, j) = (8, 0, 1)$).
\ Consequently, by \eqref{sum_Vk2}, we obtain \eqref{main_VV}.
\proofend

Now let
 \[
   U_k := X_k + \beta X_{k-1} , \qquad k \in \ZZ_+ ,
 \]
 with the convention \ $U_{-1} := U_0 := 0$.
\ In Appendix \ref{app_A}, in Remark \ref{REM_Putzer} one can find a detailed
 motivation of the definition of \ $U_k$, \ $k \in \NN$.
\ One can observe that \ $U_k \geq 0$ \ for all \ $k \in \ZZ_+$, \ and, by
 \ $\alpha + \beta = 1$,
 \begin{equation}\label{rec_U}
  U_k = U_{k-1} + M_k + \mu , \qquad k \in \ZZ_+ ,
 \end{equation}
 hence \ $(U_k)_{k\in\ZZ_+}$ \ is a nonnegative unstable AR(1) process with
 positive drift \ $\mu$ \ sharing the innovation \ $(M_k)_{k\in\NN}$
 \ with the stable AR(1) process \ $(V_k)_{k\in\ZZ_+}$.

Consider the sequence of stochastic processes
 \[
   \bcZ^{(n)}_t
   := \begin{bmatrix}
       \cM_t^{(n)} \\
       \cN_t^{(n)} \\
       \cP_t^{(n)}
      \end{bmatrix}
   := \sum_{k=1}^\nt
       \bZ^{(n)}_k
   \qquad \text{with} \qquad
   \bZ^{(n)}_k
   := \begin{bmatrix}
       n^{-1} M_k \\
       n^{-2} M_k U_{k-1} \\
       n^{-3/2} M_k V_{k-1}
      \end{bmatrix} , \qquad t \in \RR_+ , \quad k, n \in \NN .
 \]
Theorem \ref{main_Ad} will follow from Lemma \ref{XV_main_VV} and the
 following theorem (which will be detailed after Theorem \ref{main_conv}).

\begin{Thm}\label{main_conv}
Under the assumptions of Theorem \ref{main} we have
 \begin{equation}\label{conv_Z}
   \bcZ^{(n)} \distr \bcZ \qquad \text{as \ $n\to\infty$,}
 \end{equation}
 where the process \ $(\bcZ_t)_{t\in\RR_+}$ \ with values in \ $\RR^3$ \ is the
 unique strong solution of the SDE
 \begin{equation}\label{ZSDE}
   \dd \bcZ_t = \gamma(t, \bcZ_t) \, \dd \bcW_t ,
   \qquad t \in \RR_+ ,
 \end{equation}
 with initial value \ $\bcZ_0 = \bzero$, \ where
 \ $\bcW_t := \begin{bmatrix} \cW_t & \tcW_t \end{bmatrix}^\top$,
 \ $t \in \RR_+$, \ being a 2-dimensional standard Wiener process, and
 \ $\gamma : \RR_+ \times \RR^3 \to \RR^{3\times 2}$ \ is defined by
 \[
   \gamma(t, \bx)
   := \begin{bmatrix}
       \sqrt{\frac{2 \alpha \beta}{1 + \beta}}[(x_1 + \mu t)^+]^{1/2}
        & 0 \\
       \sqrt{\frac{2 \alpha \beta}{1 + \beta}}[(x_1 + \mu t)^+]^{3/2}
        & 0 \\
       0 & \frac{2 \beta \sqrt{\alpha}}{(1 + \beta)^{3/2}}
           (x_1 + \mu t)
      \end{bmatrix}
 \]
 for \ $t \in \RR_+$ \ and \ $\bx = (x_1 , x_2 , x_3) \in \RR^3$.
\end{Thm}

Indeed, the unique strong solution of \eqref{ZSDE} with initial value
 \ $\bcZ_0 = \bzero$ \ can be written in form
 \[
   \bcZ_t
   := \begin{bmatrix}
       \cM_t \\
       \cN_t \\
       \cP_t
      \end{bmatrix}
   := \begin{bmatrix}
       (1 + \beta) \cX_t - \mu t \\[1mm]
       (1 + \beta)
       \sqrt{2 \alpha \beta} \int_0^t \cX_s^{3/2} \, \dd \cW_s \\[1mm]
       \frac{2 \beta \sqrt{\alpha}}{\sqrt{1 + \beta}}
       \int_0^t \cX_s \, \dd \tcW_s
      \end{bmatrix} , \qquad t\in\RR_+ ,
 \]
 since, by Remark \ref{REMARK2},
 \begin{align*}
  \dd \bcZ_t
  = \begin{bmatrix}
           \dd\cM_t \\
           \dd\cN_t \\
           \dd\cP_t
    \end{bmatrix}
  = \begin{bmatrix}
     \sqrt{\frac{2\alpha\beta}{1+\beta}}[(\cM_t + \mu t)^+]^{1/2}
      \, \dd\cW_t \\[1mm]
     \sqrt{\frac{2\alpha\beta}{1+\beta}}[(\cM_t + \mu t)^+]^{3/2}
      \, \dd\cW_t \\[1mm]
     \frac{2\beta\sqrt{\alpha}}{(1+\beta)^{3/2}}
      (\cM_t + \mu t) \, \dd\tcW_t
    \end{bmatrix}
  = \begin{bmatrix}
           \sqrt{2 \alpha \beta}\cX_t^{1/2}\,\dd \cW_t \\[1mm]
           (1+\beta)\sqrt{2\alpha\beta}\cX_t^{3/2}\,\dd\cW_t  \\[1mm]
           \frac{2 \beta \sqrt{\alpha}}{\sqrt{1 + \beta}}\cX_t\,\dd\tcW_t
          \end{bmatrix},
        \qquad t\in\RR_+ .
 \end{align*}
By the method of the proof of \ $\cX^{(n)} \distr \cX$ \ in Theorem 3.1 in
 Barczy et al.\ \cite{BarIspPap0} one can easily derive
 \begin{align}\label{convXZ}
   \begin{bmatrix} \cX^{(n)} \\ \bcZ^{(n)} \end{bmatrix}
   \distr \begin{bmatrix} \cX \\ \bcZ \end{bmatrix} \qquad
   \text{as \ $n \to \infty$.}
 \end{align}
More precisely, using that
 \[
   X_k = \sum_{j=1}^k
         (M_j + \mu) \be_1^\top A^{k-j} \be_1, \qquad k \in \NN ,
  \qquad \text{where \ $\be_1:=\begin{bmatrix} 1 \\ 0 \end{bmatrix}$,}
 \]
 see, e.g., Barczy et al. \cite[(3.11)]{BarIspPap0}, we have
 \[
   \begin{bmatrix}
    \cX^{(n)} \\
    \bcZ^{(n)}
   \end{bmatrix}
   = \psi_n(\bcZ^{(n)}), \qquad n \in \NN ,
 \]
 where the mapping \ $\psi_n : \DD(\RR_+, \RR^3) \to \DD(\RR_+, \RR^4)$ \ is
 given by
 \[
   \psi_n(f_1,f_2,f_3)(t)
   := \begin{bmatrix}
       \sum_{j=1}^\nt
        \left(f_1\left(\frac{j}{n}\right) - f_1\left(\frac{j-1}{n}\right)
              + \frac{\mu}{n}\right) \be_1^\top A^{\nt-j} \be_1 \\
       f_1(t) \\
       f_2(t) \\
       f_3(t)
      \end{bmatrix}
 \]
 for \ $f_1, f_2, f_3 \in \DD(\RR_+, \RR)$, \ $t \in \RR_+$, \ $n \in \NN$.
\ Further, using that, by Remark \ref{REMARK2},
 \[
   \cX_t = \frac{1}{1+\beta}(\cM_t + \mu t) , \qquad t \in \RR_+ ,
 \]
 we have
 \[
   \begin{bmatrix}
    \cX \\
    \bcZ
   \end{bmatrix}
   = \psi(\bcZ) ,
 \]
 where the mapping \ $\psi : \DD(\RR_+, \RR^3) \to \DD(\RR_+, \RR^4)$ \ is
 given by
 \[
   \psi(f_1,f_2,f_3)(t) := \begin{bmatrix}
                           \frac{1}{1+\beta}(f_1(t) + \mu t) \\
                           f_1(t) \\
                           f_2(t) \\
                           f_3(t)
                          \end{bmatrix}
 \]
 for \ $f_1, f_2, f_3 \in \DD(\RR_+, \RR)$ \ and \ $t \in \RR_+$.
\ By page 603 in Barczy et al. \cite{BarIspPap0}, the mappings \ $\psi_n$,
 \ $n \in \NN$, \ and \ $\psi$ \ are measurable (the latter one is continuous
 too), since the coordinate functions are measurable.
Using page 604 in Barczy et al. \cite{BarIspPap0}, we get the set
 \[
   C := \big\{f \in \CC(\RR_+,\RR^3) : f(0) = \bzero \in \RR^3 \big\}
 \]
 has the properties \ $C\subseteq C_{\psi,(\psi_n)_{n\in\NN}}$ \ with
 \ $C \in \cB(\DD(\RR_+, \RR^3))$ \ and \ $\PP(\bcZ \in C) = 1$, \ where
 \ $C_{\psi,(\psi_n)_{n\in\NN}}$ \ is defined in Appendix \ref{app_B}.
Hence, by \eqref{conv_Z} and Lemma \ref{Conv2Funct}, we have
 \[
   \begin{bmatrix}
     \cX^{(n)} \\
     \bcZ^{(n)}
   \end{bmatrix}
    =\psi_n(\bcZ^{(n)})
     \distr
     \psi(\bcZ)
     =\begin{bmatrix}
        \cX \\
        \bcZ
       \end{bmatrix}
       \qquad \text{as \ $n\to\infty$,}
 \]
 as desired.
Next, similarly to the proof of \eqref{seged2}, by Lemmas \ref{Conv2Funct} and
 \ref{Marci}, convergence \eqref{convXZ} implies
 \begin{align}\label{help1}
   \sum_{k=1}^n
    \begin{bmatrix}
     n^{-1} M_k \\
     n^{-3} X_{k-1}^2 \\
     n^{-2} X_{k-1} \\
     n^{-2} M_k U_{k-1} \\
     n^{-3/2} M_k V_{k-1}
    \end{bmatrix}
   \distr \begin{bmatrix}
           (1 + \beta) \cX_1 - \mu \\
           \int_0^1 \cX_t^2 \, \dd t \\
           \int_0^1 \cX_t \, \dd t \\
           (1+\beta) \sqrt{2\alpha \beta} \int_0^1 \cX_t^{3/2} \, \dd \cW_t \\
           \frac{2\beta\sqrt{\alpha}}{\sqrt{1+\beta}}
           \int_0^1 \cX_t \, \dd \tcW_t
          \end{bmatrix} \qquad
   \text{as \ $n \to \infty$.}
 \end{align}
Namely,
 \begin{align*}
  \sum_{k=1}^\nt
    \begin{bmatrix}
     n^{-1} M_k \\
     n^{-3} X_{k-1}^2 \\
     n^{-2} X_{k-1} \\
     n^{-2} M_k U_{k-1} \\
     n^{-3/2} M_k V_{k-1}
    \end{bmatrix}
    = \widetilde\psi_n
      \left( \begin{bmatrix}
      \cX^{(n)} \\
      \bcZ^{(n)}
     \end{bmatrix} \right)
     (t),
    \qquad t\in\RR_+,\;\; n\in\NN, 
 \end{align*}
 where \ $\widetilde\psi_n : \DD(\RR_+,\RR^4) \to \DD(\RR_+,\RR^5)$ \ is given by 
 \[
   \widetilde\psi_n(f_1,f_2,f_3,f_4)(t)
      := \begin{bmatrix}
           f_2(t) \\[1mm]
       n^{-1} \sum_{k=1}^\nt \left( f_1\left(\frac{k-1}{n}\right)\right)^2 \\[1mm]
       n^{-1} \sum_{k=1}^\nt f_1\left(\frac{k-1}{n}\right) \\[1mm]
      f_3(t) \\
      f_4(t)
    \end{bmatrix}
 \] 
 for \ $f_1,f_2,f_3,f_4\in \DD(\RR_+,\RR)$, $t\in\RR_+$, $n\in\NN$. 
\ Further, 
 \[
  \begin{bmatrix}
           (1 + \beta) \cX_t - \mu \\
           \int_0^t \cX_s^2 \, \dd s \\
           \int_0^t \cX_s \, \dd s \\
           (1+\beta) \sqrt{2\alpha \beta} \int_0^t \cX_s^{3/2} \, \dd \cW_s \\
           \frac{2\beta\sqrt{\alpha}}{\sqrt{1+\beta}}
           \int_0^t \cX_s \, \dd \tcW_s
          \end{bmatrix}
   = \widetilde\psi\left( \begin{bmatrix}
        \cX \\
        \bcZ
       \end{bmatrix} \right)(t),
     \qquad t\in\RR_+,         
 \]
 where \ $\widetilde\psi : \DD(\RR_+,\RR^4) \to \DD(\RR_+,\RR^5)$ \ is given by 
 \[
   \widetilde\psi(f_1,f_2,f_3,f_4)(t)
      := \begin{bmatrix}
           f_2(t) \\
           \int_0^1 (f_1(s))^2 \, \dd s \\
           \int_0^1 f_1(s) \, \dd s \\
           f_3(t) \\
           f_4(t)
          \end{bmatrix}
 \]  
  for \ $f_1,f_2,f_3,f_4\in \DD(\RR_+,\RR)$, $t\in\RR_+$. 
\ As in the proof of Lemma \ref{Marci}, one can check that the set  
 \[
   \widetilde C := \Big\{ f\in\CC(\RR_+,\RR^4) : f(0) = \bzero\in\RR^4 \Big\}
 \]
 has the properties \ $\widetilde C \subseteq C_{\widetilde\psi, (\widetilde \psi_n)_{n\in\NN}}$ \ with 
 \ $\widetilde C\in\cB(\DD(\RR_+,\RR^4))$ \ and 
 \[
  \PP\left( \begin{bmatrix}
        \cX \\
        \bcZ
       \end{bmatrix} 
       \in \widetilde C \right) =1.
 \]  
Hence, by \eqref{convXZ} and Lemma \ref{Conv2Funct}, we have 
 \[
    \widetilde\psi_n
      \left( \begin{bmatrix}
      \cX^{(n)} \\
      \bcZ^{(n)}
     \end{bmatrix} \right)
      \distr
      \widetilde\psi
      \left( \begin{bmatrix}
      \cX \\
      \bcZ
     \end{bmatrix} \right)
    \quad \text{as \ $n\to\infty$,} 
 \] 
 which yields \eqref{help1}.

Using \ $U_{k-1} = (1+\beta) X_{k-1} - \beta V_{k-1}$ \ and convergence of the
 third coordinates in \ $\bcZ^{(n)} \distr \bcZ$ \ as \ $n\to\infty$ \ we obtain
 \[
   n^{-2} \left( \sum_{k=1}^n M_k X_{k-1}
                - \frac{1}{1+\beta} \sum_{k=1}^n M_k U_{k-1} \right)
   = \frac{\beta}{(1+\beta) n^2} \sum_{k=1}^n M_k V_{k-1}
   \stoch 0  \qquad
   \text{as \ $n \to \infty$.}
 \]
Using \eqref{tbA_tbd}, the above two convergences and Lemma \ref{XV_main_VV}
 we obtain Theorem \ref{main_Ad} by Slutsky's lemma.

\section{Proof of Theorem \ref{main_conv}}
\label{section_proof_main_conv}

In order to show convergence \ $\bcZ^{(n)} \distr \bcZ$, \ we apply Theorem
 \ref{Conv2DiffThm} with the special choices \ $\bcU := \bcZ$,
 \ $\bU^{(n)}_k := \bZ^{(n)}_k$, \ $n,k\in\NN$,
 \ $(\cF_k^{(n)})_{k\in\ZZ_+}:=(\cF_k)_{k\in\ZZ_+}$ \ and the function \ $\gamma$
 \ which is defined in Theorem \ref{main_conv}.
Note that the arguments in Section \ref{section_proof_main} and Remark
 \ref{REMARK1} show that the SDE \eqref{ZSDE} admits a unique strong solution
 \ $(\bcZ_t^\bz)_{t\in\RR_+}$ \ for all initial values
 \ $\bcZ_0^\bz = \bz \in \RR^3$.

Now we show that conditions (i) and (ii) of Theorem \ref{Conv2DiffThm} hold.
The conditional variances have the form
 \[
   \EE\bigl(\bZ^{(n)}_k (\bZ^{(n)}_k)^\top \mid \cF_{k-1}\bigr) =
   \EE( M_k^2 \mid \cF_{k-1} )
   \begin{bmatrix}
    n^{-2}
    & n^{-3} U_{k-1}
    & n^{-5/2} V_{k-1} \\
    n^{-3} U_{k-1}
    & n^{-4} U_{k-1}^2
    & n^{-7/2} U_{k-1} V_{k-1} \\
    n^{-5/2} V_{k-1}
    & n^{-7/2} U_{k-1} V_{k-1}
    & n^{-3} V_{k-1}^2
   \end{bmatrix}
 \]
 for \ $n \in \NN, \ k\in\{1,\ldots,n\}$, \ and
 \[
   \gamma(s,\bcZ_s^{(n)}) \gamma(s,\bcZ_s^{(n)})^\top
   = \begin{bmatrix}
      \frac{2 \alpha \beta}{1 + \beta} (\cM_s^{(n)} + \mu s)
       & \frac{2 \alpha \beta}{1 + \beta} (\cM_s^{(n)} + \mu s)^2 & 0 \\
      \frac{2 \alpha \beta}{1 + \beta} (\cM_s^{(n)} + \mu s)^2
       & \frac{2 \alpha \beta}{1 + \beta} (\cM_s^{(n)} + \mu s)^3 & 0 \\
      0 & 0
       & \frac{4 \alpha \beta^2}{(1 + \beta)^3} (\cM_s^{(n)} + \mu s)^2
     \end{bmatrix}
 \]
 for \ $s \in \RR_+$, \ where we used that
 \ $(\cM^{(n)}_s + \mu s)^+ = \cM^{(n)}_s + \mu s$, \ $s \in \RR_+$,
 \ $n \in \NN$, \ see Barczy et al. \cite[page 598]{BarIspPap0} or \eqref{M+}
 later on.
In order to check condition (i) of Theorem \ref{Conv2DiffThm}, we need  to
 prove that for each \ $T>0$,
 \begin{gather}
  \sup_{t\in[0,T]}
   \bigg| \frac{1}{n^2} \sum_{k=1}^{\nt} \EE( M_k^2 \mid \cF_{k-1} )
          - \frac{2 \alpha \beta}{1 + \beta}
            \int_0^t (\cM_s^{(n)} + \mu s) \, \dd s \bigg|
  \stoch 0 , \label{Zcond1} \\
  \sup_{t\in[0,T]}
   \bigg| \frac{1}{n^3} \sum_{k=1}^{\nt} U_{k-1} \EE( M_k^2 \mid \cF_{k-1} )
          - \frac{2 \alpha \beta}{1 + \beta}
            \int_0^t (\cM_s^{(n)} + \mu s)^2 \, \dd s \bigg|
  \stoch 0 , \label{Zcond2} \\
  \sup_{t\in[0,T]}
   \bigg| \frac{1}{n^4} \sum_{k=1}^{\nt} U_{k-1}^2 \EE( M_k^2 \mid \cF_{k-1} )
          - \frac{2 \alpha \beta}{1 + \beta}
            \int_0^t (\cM_s^{(n)} + \mu s)^3 \, \dd s \bigg|
  \stoch 0 , \label{Zcond3} \\
  \sup_{t\in[0,T]}
   \bigg| \frac{1}{n^3} \sum_{k=1}^{\nt} V_{k-1}^2 \EE( M_k^2 \mid \cF_{k-1} )
          - \frac{4 \alpha \beta^2}{(1 + \beta)^3}
            \int_0^t (\cM_s^{(n)} + \mu s)^2 \, \dd s \bigg|
  \stoch 0 , \label{Zcond4} \\
  \sup_{t\in[0,T]}
   \bigg| \frac{1}{n^{5/2}}
          \sum_{k=1}^{\nt} V_{k-1} \EE( M_k^2 \mid \cF_{k-1} ) \bigg|
  \stoch 0 , \label{Zcond5} \\
  \sup_{t\in[0,T]}
   \bigg| \frac{1}{n^{7/2}}
          \sum_{k=1}^{\nt} U_{k-1} V_{k-1} \EE( M_k^2 \mid \cF_{k-1} ) \bigg|
  \stoch 0  \label{Zcond6}
 \end{gather}
 as \ $n \to \infty$.
\ Covergence \eqref{Zcond1} follows from (5.1) in Barczy et al.\
 \cite{BarIspPap0} with the special choices \ $p=2$, \ $\alpha_1 = \alpha$
 \ and \ $\alpha_2 = \beta$.

Next we turn to prove \eqref{Zcond2}.
Since \ $\alpha + \beta = 1$, \ by \eqref{Mk}, we get
 \begin{align}\label{M+}
  \begin{split}
   \cM_s^{(n)} + \mu s
   &= \frac{1}{n}
      \sum_{k=1}^\ns
       \big( X_k - \alpha X_{k-1} - \beta X_{k-2} - \mu \big)
      + \mu s \\
   &= \frac{1}{n}
      \big( X_\ns + \beta X_{\ns -1} \big)
      + \frac{ns - \ns}{n} \mu
    = \frac{1}{n} U_\ns + \frac{ns - \ns}{n} \mu
  \end{split}
 \end{align}
 for \ $s\in\RR_+$, \ $n\in\NN$.
 \ Thus
 \begin{align*}
  \int_0^t (\cM_s^{(n)} + \mu s)^2 \, \dd s
  &= \frac{1}{n^3} \sum_{k=1}^{\nt-1} U_k^2
     + \frac{\mu}{n^3} \sum_{k=1}^{\nt-1} U_k
     + \frac{nt - \nt}{n^3} U_\nt^2 \\
  &\phantom{\quad}
     + \frac{\mu (nt - \nt)^2}{n^3} U_\nt
     + \frac{\nt + (nt - \nt)^3}{3n^3} \mu^2 .
 \end{align*}
Since
 \begin{align}\label{XUV}
  X_{k-1} = \frac{1}{1+\beta} (U_k - V_k) , \qquad
  X_k = \frac{1}{1+\beta} (U_k + \beta V_k) , \qquad k \in \NN ,
 \end{align}
 using Lemma \ref{Moments}, we obtain
 \begin{align}\nonumber
  \sum_{k=1}^{\nt} U_{k-1} \EE( M_k^2 \mid \cF_{k-1} )
  &= \sum_{k=1}^{\nt}
      U_{k-1} \big[ \alpha \beta (X_{k-1} + X_{k-2}) + \sigma^2 \big] \\
  &= \sum_{k=1}^{\nt}
      U_{k-1} \left[ \frac{\alpha \beta}{1 + \beta} ( 2 U_{k-1} - \alpha V_{k-1})
                    + \sigma^2 \right] \label{UM2F} \\
  &= \frac{2 \alpha \beta}{1 + \beta} \sum_{k=1}^\nt U_{k-1}^2
     - \frac{\alpha^2 \beta}{1 + \beta} \sum_{k=1}^\nt U_{k-1} V_{k-1}
     + \sigma^2 \sum_{k=1}^\nt U_{k-1} .\nonumber
 \end{align}
Thus, in order to show \eqref{Zcond2}, it suffices to prove
 \begin{gather}
  n^{-3} \sum_{k=1}^{\nT} |U_k V_k| \stoch 0 , \label{2supsumUV} \\
  n^{-3} \sum_{k=1}^{\nT} U_k \stoch 0 , \label{2supsumU} \\
  n^{-3/2} \sup_{t \in [0,T]} U_\nt \stoch 0, \label{2supU} \\
  n^{-3} \sup_{t \in [0,T]} \left[ \nt + (nt - \nt)^3 \right] \to 0
    \label{2supnt}
 \end{gather}
 as \ $n \to \infty$.
\ Using \eqref{seged_UV_UNIFORM1} with \ $(\ell, i, j) = (8, 1, 1)$ \ and
 \ $(\ell, i, j) = (8, 1, 0)$, \ we have \eqref{2supsumUV} and
 \eqref{2supsumU}, respectively.
Using \eqref{seged_UV_UNIFORM2} with \ $(\ell, i, j) = (8, 1, 0)$, \ we have
 \eqref{2supU}.
Clearly, \eqref{2supnt} follows from \ $|nt - \nt| \leq 1$, \ $n \in \NN$,
 \ $t \in \RR_+$, \ thus we conclude \eqref{Zcond2}.

Now we turn to check \eqref{Zcond3}.
Again by \eqref{M+}, we have
 \begin{align*}
  \int_0^t (\cM_s^{(n)} + \mu s)^3 \, \dd s
  &= \frac{1}{n^4} \sum_{k=1}^{\nt-1} U_k^3
     + \frac{3 \mu}{2n^4} \sum_{k=1}^{\nt-1} U_k^2
     + \frac{\mu^2}{n^4}
       \sum_{k=1}^{\nt-1} U_k
     + \frac{nt - \nt}{n^4} U_\nt^3 \\
  &\phantom{\quad}
     + \frac{3 \mu (nt - \nt)^2}{2n^4} U_\nt^2
     + \frac{\mu^2 (nt - \nt)^3}{n^4} U_\nt
     + \frac{\nt + (nt - \nt)^4}{4n^4} \mu^3 .
 \end{align*}
Using Lemma \ref{Moments}, we obtain
 \begin{align*}
  \sum_{k=1}^{\nt} U_{k-1}^2 \EE( M_k^2 \mid \cF_{k-1} )
  &= \sum_{k=1}^{\nt}
      U_{k-1}^2 \big[ \alpha \beta (X_{k-1} + X_{k-2}) + \sigma^2 \big] \\
  &= \sum_{k=1}^{\nt}
      U_{k-1}^2
      \left[ \frac{\alpha \beta}{1 + \beta} ( 2 U_{k-1} - \alpha V_{k-1})
             + \sigma^2 \right] \\
  &= \frac{2 \alpha \beta}{1 + \beta} \sum_{k=1}^\nt U_{k-1}^3
     - \frac{\alpha^2 \beta}{1 + \beta} \sum_{k=1}^\nt U_{k-1}^2 V_{k-1}
     + \sigma^2 \sum_{k=1}^\nt U_{k-1}^2 .
 \end{align*}
Thus, in order to show \eqref{Zcond3}, it suffices to prove
 \begin{gather}
  n^{-4} \sum_{k=1}^{\nT} | U_k^2 V_k | \stoch 0 , \label{3supsumUUV} \\
  n^{-4} \sum_{k=1}^{\nT} U_k^2 \stoch 0 , \label{3supsumUU} \\
  n^{-4} \sum_{k=1}^{\nT} U_k \stoch 0 , \label{3supsumU} \\
  n^{-4/3} \sup_{t \in [0,T]} U_\nt \stoch 0 , \label{3supU} \\
  n^{-4} \sup_{t \in [0,T]} \left[ \nt + (nt - \nt)^4 \right] \to 0
    \label{3supnt}
 \end{gather}
 as \ $n \to \infty$.
\ Using \eqref{seged_UV_UNIFORM1} with \ $(\ell, i, j) = (8, 2, 1)$,
 \ $(\ell, i, j) = (8, 2, 0)$ \ and \ $(\ell, i, j) = (8, 1, 0)$, \ we have
 \eqref{3supsumUUV}, \eqref{3supsumUU} and \eqref{3supsumU}, respectively.
Using \eqref{seged_UV_UNIFORM2} with \ $(\ell, i, j) = (8, 1, 0)$, \ we have
 \eqref{3supU}.
Clearly, \eqref{3supnt} follows again from \ $|nt - \nt| \leq 1$,
 \ $n \in \NN$, \ $t \in \RR_+$, \ thus we conclude \eqref{Zcond3}.

Next we turn to prove \eqref{Zcond4}.
By \eqref{UM2F}, \eqref{2supsumUV} and \eqref{2supsumU} we get
 \begin{align}\label{Zcond2a}
  n^{-3}
  \sup_{t \in [0,T]}
   \left| \sum_{k=1}^{\nt} U_{k-1} \EE( M_k^2 \mid \cF_{k-1} )
          - \frac{2 \alpha \beta}{1 + \beta} \sum_{k=1}^{\nt} U_{k-1}^2 \right|
  \stoch 0 \qquad \text{as \ $n \to \infty$}
 \end{align}
 for all \ $T>0$.
\ Using \eqref{Zcond2}, in order to prove \eqref{Zcond4}, it is sufficient to
 show that
 \begin{align}\label{Zcond4a}
  n^{-3}
  \sup_{t \in [0,T]}
   \left| \sum_{k=1}^{\nt} V_{k-1}^2 \EE( M_k^2 \mid \cF_{k-1} )
          - \frac{4 \alpha \beta^2}{(1 + \beta)^3}
            \sum_{k=1}^{\nt} U_{k-1}^2 \right|
  \stoch 0 \qquad \text{as \ $n \to \infty$}
 \end{align}
 for all \ $T>0$.
\ As in the previous case, using Lemma \ref{Moments} and \eqref{XUV}, we
 obtain
 \begin{align*}
  \sum_{k=1}^{\nt} V_{k-1}^2 \EE( M_k^2 \mid \cF_{k-1} )
  &= \sum_{k=1}^{\nt}
      V_{k-1}^2 [\alpha \beta (X_{k-1} + X_{k-2}) + \sigma^2] \\
  &= \sum_{k=1}^{\nt}
      V_{k-1}^2
      \left[\frac{\alpha\beta}{1+\beta}(2U_{k-1} - \alpha V_{k-1})
            + \sigma^2 \right] \\
  &= \frac{2 \alpha \beta}{1 + \beta} \sum_{k=1}^\nt U_{k-1} V_{k-1}^2
     - \frac{\alpha^2 \beta}{1 + \beta} \sum_{k=1}^\nt V_{k-1}^3
     + \sigma^2 \sum_{k=1}^\nt V_{k-1}^2 .
 \end{align*}
Using \eqref{seged_UV_UNIFORM1} with \ $(\ell, i, j) = (8, 0, 3)$ \ and
 \ $(\ell, i, j) = (8, 0, 2)$, \ we have
 \begin{align*}
  n^{-3} \sum_{k=1}^{\nT} | V_k |^3 \stoch 0 , \qquad
  n^{-3} \sum_{k=1}^{\nT} V_k^2 \stoch 0 \qquad \text{as \ $n \to \infty$,}
 \end{align*}
 hence \eqref{Zcond4a} will follow from
 \begin{align}\label{Zcond4b}
  n^{-3}
  \sup_{t \in [0,T]}
   \left| \sum_{k=1}^{\nt} U_{k-1} V_{k-1}^2
          - \frac{2 \beta}{(1 + \beta)^2} \sum_{k=1}^{\nt} U_{k-1}^2 \right|
  \stoch 0 \qquad \text{as \ $n \to \infty$}
 \end{align}
 for all \ $T>0$.

The aim of the following discussion is to decompose
 \ $\sum_{k=1}^{\nt} U_{k-1} V_{k-1}^2
    - 2 \beta (1 + \beta)^{-2} \sum_{k=1}^{\nt} U_{k-1}^2$
 \ as a sum of a martingale and some negligible terms.
Using recursions \eqref{rec_V}, \eqref{rec_U} and Lemma \ref{Moments} (formula
 \eqref{Mcond}), we obtain
 \begin{align*}
  \EE(U_{k-1} V_{k-1}^2 \mid \cF_{k-2})
  &= \EE\big((U_{k-2} + M_{k-1} + \mu)
             (- \beta V_{k-2} + M_{k-1} + \mu)^2 \mid \cF_{k-2}\big) \\
  &= \beta^2 U_{k-2} V_{k-2}^2
     + \alpha \beta (X_{k-2} + X_{k-3}) (U_{k-2} - 2 \beta V_{k-2} + 3 \mu)
     + \EE(M_{k-1}^3 \mid \cF_{k-2}) \\
  &\quad
     + \text{constant}
     + \text{linear combination of \ $U_{k-2} V_{k-2}$, \ $V_{k-2}^2$, \ $U_{k-2}$
             \ and \ $V_{k-2}$.}
 \end{align*}
Using again Lemma \ref{Moments} (formula \eqref{M3cond}) and \eqref{XUV}, we
 get
 \begin{align}\label{cond_UVV}
  \begin{split}
   &\EE(U_{k-1} V_{k-1}^2 \mid \cF_{k-2})\\
   &\quad
    = \beta^2 U_{k-2} V_{k-2}^2
      + \frac{\alpha \beta}{1 + \beta} ( 2 U_{k-2} - \alpha V_{k-2} )
        (U_{k-2} - 2\beta V_{k-2} + 3 \mu)
      + \EE(M_{k-1}^3\mid \cF_{k-2}) \\
   &\quad\quad
     + \text{constant}
     + \text{linear combination of \ $U_{k-2} V_{k-2}$, \ $V_{k-2}^2$, \ $U_{k-2}$
              \ and \ $V_{k-2}$}\\
   &\quad= \beta^2 U_{k-2} V_{k-2}^2 + \frac{2 \alpha \beta}{1 + \beta} U_{k-2}^2
      + \text{constant}\\
   &\quad\quad
      + \text{linear combination of \ $U_{k-2} V_{k-2}$, \ $V_{k-2}^2$, \ $U_{k-2}$
              \ and \ $V_{k-2}$.}
  \end{split}
 \end{align}
Thus
 \begin{align*}
  &\sum_{k=1}^{\nt} U_{k-1} V_{k-1}^2
   = \sum_{k=2}^{\nt}
      \big[U_{k-1} V_{k-1}^2 - \EE(U_{k-1} V_{k-1}^2 \mid \cF_{k-2}) \big]
     + \sum_{k=2}^{\nt} \EE(U_{k-1} V_{k-1}^2 \mid \cF_{k-2}) \\
  &= \sum_{k=2}^{\nt}
      \big[U_{k-1} V_{k-1}^2 - \EE(U_{k-1} V_{k-1}^2 \mid \cF_{k-2}) \big]
     + \beta^2 \sum_{k=2}^{\nt} U_{k-2} V_{k-2}^2
     + \frac{2 \alpha \beta}{1 + \beta} \sum_{k=2}^{\nt} U_{k-2}^2\\
  &\quad
     + \OO(n)
     + \text{linear combination of \ $\sum_{k=2}^{\nt} U_{k-2} V_{k-2}$,
             \ $\sum_{k=2}^{\nt} V_{k-2}^2$, \ $\sum_{k=2}^{\nt} U_{k-2}$
             \ and \ $\sum_{k=2}^{\nt} V_{k-2}$.}
 \end{align*}
Consequently,
 \begin{align*}
  \sum_{k=1}^{\nt} U_{k-1} V_{k-1}^2
  &= \frac{1}{1 - \beta^2}
     \sum_{k=2}^{\nt}
      \big[U_{k-1} V_{k-1}^2 - \EE(U_{k-1} V_{k-1}^2 \mid \cF_{k-2}) \big]
     + \frac{2 \alpha \beta}{(1 + \beta)(1 - \beta^2)}
       \sum_{k=2}^{\nt} U_{k-2}^2 \\
  &\quad
     - \frac{\beta^2}{1 - \beta^2} U_{\nt - 1} V_{\nt - 1}^2 + \OO(n) \\
  &\quad
     + \text{linear combination of \ $\sum_{k=2}^{\nt} U_{k-2} V_{k-2}$,
             \ $\sum_{k=2}^{\nt} V_{k-2}^2$, \ $\sum_{k=2}^{\nt} U_{k-2}$
             \ and \ $\sum_{k=2}^{\nt} V_{k-2}$.}
 \end{align*}
Using \eqref{seged_UV_UNIFORM4} with \ $(\ell, i, j) = (8, 1, 2)$ \ we have
 \begin{align*}
    n^{-3}\sup_{t \in [0,T]}\,
           \Biggl\vert \sum_{k=2}^\nt
                    \big[U_{k-1} V_{k-1}^2
                         - \EE(U_{k-1} V_{k-1}^2 \mid \cF_{k-2}) \big]
           \Biggr\vert
\stoch 0 \qquad \text{as \ $n\to\infty$.}
 \end{align*}
Thus,  in order to show \eqref{Zcond4b}, it suffices to prove
 \begin{gather}
  n^{-3} \sum_{k=1}^{\nT} | U_k V_k | \stoch 0 , \label{4supsumUV} \\
  n^{-3} \sum_{k=1}^{\nT} V_k^2 \stoch 0 , \label{4supsumVV} \\
  n^{-3} \sum_{k=1}^{\nT} U_k \stoch 0 , \label{4supsumU} \\
  n^{-3} \sum_{k=1}^{\nT} |V_k| \stoch 0 , \label{4supsumV} \\
  n^{-3} \sup_{t \in [0,T]} U_\nt V_\nt^2 \stoch 0 , \label{4supUVV}\\
  n^{-3} \sup_{t \in [0,T]} U^2_{\nt} \stoch 0 \label{4supUU}
   \end{gather}
 as \ $n \to \infty$.
\ Using \eqref{seged_UV_UNIFORM1} with \ $(\ell, i, j) = (8, 1, 1)$,
 \ $(\ell, i, j) = (8, 0, 2)$, \ $(\ell, i, j) = (8, 1, 0)$, \ and
 \ $(\ell, i, j) = (8, 0, 1)$, \ we have \eqref{4supsumUV}, \eqref{4supsumVV},
 \eqref{4supsumU} and \eqref{4supsumV}.
Using \eqref{seged_UV_UNIFORM2} with \ $(\ell, i, j) = (8, 1, 2)$ \ and
 \ $(\ell, i, j) = (8, 2, 0)$, \ we have \eqref{4supUVV} and \eqref{4supUU}.
Thus we conclude \eqref{Zcond4}.

For \eqref{Zcond5}, consider
 \begin{align*}
  \sum_{k=1}^{\nt} V_{k-1} \EE( M_k^2 \mid \cF_{k-1} )
  &= \sum_{k=1}^{\nt}
      V_{k-1} \big( \alpha \beta (X_{k-1} + X_{k-2}) + \sigma^2 \big) \\
  &= \sum_{k=1}^{\nt}
      V_{k-1} \left( \frac{\alpha\beta}{1+\beta} (2U_{k-1} - \alpha V_{k-1})
                    + \sigma^2 \right) \\
  &= \frac{2 \alpha \beta}{1 + \beta} \sum_{k=1}^\nt U_{k-1} V_{k-1}
     - \frac{\alpha^2 \beta}{1 + \beta} \sum_{k=1}^\nt V_{k-1}^2
     + \sigma^2 \sum_{k=1}^\nt V_{k-1} ,
 \end{align*}
 where we used Lemma \ref{Moments} and \eqref{XUV}.
Using \eqref{seged_UV_UNIFORM1} with \ $(\ell, i, j) = (8, 0, 2)$, \ and
 \ $(\ell, i, j) = (8, 0, 1)$, \ we have
 \begin{align*}
  n^{-5/2} \sum_{k=1}^{\nT} V_k^2 \stoch 0 , \qquad
  n^{-5/2} \sum_{k=1}^{\nT} |V_k| \stoch 0 \qquad \text{as \ $n \to \infty$,}
 \end{align*}
 hence \eqref{Zcond5} will follow from
 \begin{align}\label{Zcond5a}
  n^{-5/2} \sup_{t \in [0,T]} \left| \sum_{k=1}^{\nt} U_{k-1} V_{k-1} \right|
  \stoch 0 .
 \end{align}
The aim of the following discussion is to decompose
 \ $\sum_{k=1}^{\nt} U_{k-1} V_{k-1}$ \ as a sum of a martingale and some
 negligible terms.
Using the recursions \eqref{rec_V}, \eqref{rec_U} and Lemma \ref{Moments}, we
 obtain
 \begin{align*}
  \EE(U_{k-1} V_{k-1} \mid \cF_{k-2})
  &= \EE\big((U_{k-2} + M_{k-1} + \mu)
             (- \beta V_{k-2} + M_{k-1} + \mu)
             \mid \cF_{k-2}\big) \\
  &= - \beta U_{k-2} V_{k-2} + \mu U_{k-2} - \beta \mu V_{k-2}
     + \mu^2 + \EE(M_{k-1}^2 \mid \cF_{k-2}) \\
  &= - \beta U_{k-2} V_{k-2}
     + \text{constant}
     + \text{linear combination of \ $U_{k-2}$ \ and \ $V_{k-2}$.}
 \end{align*}
Thus
 \begin{align*}
  \sum_{k=1}^{\nt} U_{k-1} V_{k-1}
  &= \sum_{k=2}^{\nt} \big[ U_{k-1} V_{k-1} - \EE(U_{k-1}V_{k-1} \mid \cF_{k-2}) \big]
     + \sum_{k=2}^{\nt} \EE(U_{k-1}V_{k-1} \mid \cF_{k-2}) \\
  &= \sum_{k=2}^{\nt} \big[ U_{k-1} V_{k-1} - \EE(U_{k-1}V_{k-1} \mid \cF_{k-2}) \big]
     - \beta \sum_{k=2}^{\nt} U_{k-2} V_{k-2} \\
  &\quad
     + \OO(n)
     + \text{linear combination of \ $\sum_{k=2}^{\nt} U_{k-2}$ \ and
             \ $\sum_{k=2}^{\nt} V_{k-2}$.}
 \end{align*}
Consequently
 \begin{align*}
  \sum_{k=2}^{\nt} U_{k-1} V_{k-1}
  &= \frac{1}{1+\beta}
     \sum_{k=2}^{\nt} \big[ U_{k-1} V_{k-1} - \EE(U_{k-1}V_{k-1} \mid \cF_{k-2}) \big]
     + \frac{\beta}{1+\beta} U_{\nt-1} V_{\nt-1} \\
  &\quad
     + \OO(n)
     + \text{linear combination of \ $\sum_{k=2}^{\nt} U_{k-2}$ \ and
             \ $\sum_{k=2}^{\nt} V_{k-2}$.}
 \end{align*}
Using \eqref{seged_UV_UNIFORM4} with \ $(\ell, i, j) = (8, 1, 1)$ \ we have
 \begin{align*}
    n^{-5/2}\sup_{t \in [0,T]}\,
           \Biggl\vert \sum_{k=2}^\nt
                    \big[U_{k-1} V_{k-1}
                         - \EE(U_{k-1} V_{k-1} \mid \cF_{k-2}) \big]
           \Biggr\vert
\stoch 0 \qquad \text{as \ $n\to\infty$.}
 \end{align*}
Thus,  in order to show \eqref{Zcond5a}, it suffices to prove
 \begin{gather}
  n^{-5/2} \sum_{k=1}^{\nT} U_k \stoch 0 , \label{5supsumU} \\
  n^{-5/2} \sum_{k=1}^{\nT} |V_k| \stoch 0 , \label{5supsumV} \\
  n^{-5/2} \sup_{t \in [0,T]} | U_\nt V_\nt | \stoch 0  \label{5supUV}
 \end{gather}
 as \ $n \to \infty$.
\ Using \eqref{seged_UV_UNIFORM1} with \ $(\ell, i, j) = (8, 1, 0)$, \ and
 \ $(\ell, i, j) = (8, 0, 1)$, \ we have \eqref{5supsumU} and \eqref{5supsumV}.
Using \eqref{seged_UV_UNIFORM2} with \ $(\ell, i, j) = (8, 1, 1)$ \ we have
 \eqref{5supUV}, thus we conclude \eqref{Zcond5}.

Convergence \eqref{Zcond6} can be handled in the same way as \eqref{Zcond5}.
For completeness we present all of the details.
By Lemma \ref{Moments} and \eqref{XUV}, we have
 \begin{align*}
  \sum_{k=1}^{\nt} U_{k-1}V_{k-1} \EE( M_k^2 \mid \cF_{k-1} )
  &= \sum_{k=1}^{\nt}
      U_{k-1} V_{k-1}
      \big( \alpha \beta (X_{k-1} + X_{k-2}) + \sigma^2 \big) \\
  &= \frac{2 \alpha \beta}{1 + \beta} \sum_{k=1}^\nt U_{k-1}^2 V_{k-1}
     - \frac{\alpha^2 \beta}{1 + \beta} \sum_{k=1}^\nt U_{k-1} V_{k-1}^2
     + \sigma^2 \sum_{k=1}^\nt U_{k-1} V_{k-1} .
 \end{align*}
Using \eqref{seged_UV_UNIFORM1} with \ $(\ell, i, j) = (8, 1, 2)$, \ and
 \ $(\ell, i, j) = (8, 1, 1)$, \ we have
 \begin{align*}
  n^{-7/2} \sum_{k=1}^{\nT} U_{k-1} V_{k-1}^2 \stoch 0 , \qquad
  n^{-7/2} \sum_{k=1}^{\nT} U_{k-1} |V_{k-1}| \stoch 0 \qquad
  \text{as \ $n \to \infty$,}
 \end{align*}
 hence \eqref{Zcond6} will follow from
 \begin{align}\label{Zcond6a}
  n^{-7/2} \sup_{t \in [0,T]} \left| \sum_{k=1}^{\nt} U_{k-1}^2 V_{k-1} \right|
  \stoch 0 \qquad \text{as \ $n \to \infty$.}
 \end{align}
The aim of the following discussion is to decompose
 \ $\sum_{k=1}^{\nt} U_{k-1}^2 V_{k-1}$ \ as a sum of a martingale and some
 negligible terms.
Using the recursions \eqref{rec_V} and \eqref{rec_U}, we obtain
 \begin{align*}
  \EE(U_{k-1}^2V_{k-1} \mid \cF_{k-2})
  &= \EE\big((U_{k-2} + M_{k-1} + \mu)^2
             (- \beta V_{k-2} + M_{k-1} + \mu)
             \mid \cF_{k-2}\big) \\
  &= - \beta U_{k-2}^2 V_{k-2} + \mu U_{k-2}^2 - \beta \mu^2 V_{k-2}
     - 2 \beta \mu U_{k-2} V_{k-2} + 2 \mu^2 U_{k-2} \\
  &\phantom{=\,}
     + (2U_{k-2} - \beta V_{k-2} + 3 \mu) \EE(M_{k-1}^2 \mid \cF_{k-2})
     + \EE(M_{k-1}^3 \mid \cF_{k-2})
     + \mu^3 .
 \end{align*}
Hence, by Lemma \ref{Moments} and \eqref{XUV},
\begin{align*}
 \EE(U_{k-1}^2 V_{k-1} \mid \cF_{k-2})
  &= - \beta U_{k-2}^2 V_{k-2} + \text{constant} \\
  &\phantom{=\,}
     + \text{linear combination of \ $U_{k-2}$, $V_{k-2}$, $U_{k-2}^2$, $V_{k-2}^2$
             \ and \ $U_{k-2} V_{k-2}$.}
\end{align*}
Thus
 \begin{align*}
  &\sum_{k=1}^{\nt} U_{k-1}^2 V_{k-1}
  = \sum_{k=2}^{\nt}
     \big[ U_{k-1}^2 V_{k-1} - \EE(U_{k-1}^2 V_{k-1} \mid \cF_{k-2}) \big]
    + \sum_{k=2}^{\nt} \EE(U_{k-1}^2 V_{k-1} \mid \cF_{k-2}) \\
  &\quad
  = \sum_{k=2}^{\nt}
     \big[ U_{k-1}^2 V_{k-1} - \EE(U_{k-1}^2 V_{k-1} \mid \cF_{k-2}) \big]
    - \beta \sum_{k=2}^{\nt} U_{k-2}^2 V_{k-2} + \OO(n) \\
  &\quad\quad
    + \text{linear combination of \ $\sum_{k=1}^{\nt} U_{k-2}$,
            $\sum_{k=1}^{\nt} V_{k-2}$, $\sum_{k=1}^{\nt} U_{k-2}^2$,
            $\sum_{k=1}^{\nt} V_{k-2}^2$ \ and \ $\sum_{k=1}^{\nt} U_{k-2} V_{k-2}$.}
 \end{align*}
Consequently
 \begin{align*}
  &\sum_{k=1}^{\nt} U_{k-1}^2 V_{k-1}
  = \frac{1}{1+\beta}
    \sum_{k=2}^{\nt}
     \big[ U_{k-1}^2 V_{k-1} - \EE(U_{k-1}^2 V_{k-1} \mid \cF_{k-2}) \big]
    + \frac{\beta}{1+\beta} U_{\nt-1}^2 V_{\nt-1} \\
  &\quad
    + \OO(n)
    + \text{linear combination of \ $\sum_{k=1}^{\nt} U_{k-2}$,
            $\sum_{k=1}^{\nt} V_{k-2}$, $\sum_{k=1}^{\nt} U_{k-2}^2$,
            $\sum_{k=1}^{\nt} V_{k-2}^2$ \ and \ $\sum_{k=1}^{\nt} U_{k-2} V_{k-2}$.}
 \end{align*}
Using \eqref{seged_UV_UNIFORM4} with \ $(\ell, i, j) = (8, 2, 1)$ \ we have
 \begin{align*}
  n^{-7/2} \sup_{t\in[0,T]} \,
           \Biggl| \sum_{k=2}^\nt
                    \big[ U_{k-1}^2 V_{k-1}
                          - \EE(U_{k-1}^2 V_{k-1} \mid \cF_{k-2}) \big]
           \Biggr|
  \stoch 0 \qquad \text{as \ $n \to \infty$.}
 \end{align*}
Thus,  in order to show \eqref{Zcond6a}, it suffices to prove
 \begin{gather}
  n^{-7/2} \sum_{k=1}^{\nT} U_k \stoch 0 , \label{6supsumU} \\
  n^{-7/2} \sum_{k=1}^{\nT} U_k^2 \stoch 0 , \label{6supsumUU}
 \end{gather}
 \begin{gather}
  n^{-7/2} \sum_{k=1}^{\nT} |V_k| \stoch 0 , \label{6supsumV} \\
  n^{-7/2} \sum_{k=1}^{\nT} V_k^2 \stoch 0 , \label{6supsumVV} \\
  n^{-7/2} \sum_{k=1}^{\nT} \vert U_k V_k\vert \stoch 0 , \label{6supsumUV} \\
  n^{-7/2} \sup_{t \in [0,T]} | U_\nt^2 V_\nt | \stoch 0  \label{6supUUV}
 \end{gather}
 as \ $n \to \infty$.
\ Here \eqref{6supsumU}, \eqref{6supsumUU}, \eqref{6supsumV},
 \eqref{6supsumVV} and \eqref{6supsumUV} follow by \eqref{seged_UV_UNIFORM1},
 and \eqref{6supUUV} by \eqref{seged_UV_UNIFORM2}, thus we conclude
 \eqref{Zcond6}.

Finally, we check condition (ii) of Theorem \ref{Conv2DiffThm}, i.e., the
 conditional Lindeberg condition
 \begin{equation}\label{Zcond3_new}
  \sum_{k=1}^{\lfloor nT \rfloor}
   \EE \big( \|\bZ^{(n)}_k\|^2 \bbone_{\{\|\bZ^{(n)}_k\|>\theta\}}
             \bmid \cF_{k-1} \big)
    \stoch 0 \qquad \text{as \ $n\to\infty$ \ for all \ $\theta > 0$ \ and \ $T > 0$.}
 \end{equation}
We have
 \ $\EE \big( \|\bZ^{(n)}_k\|^2 \bbone_{\{\|\bZ^{(n)}_k\|>\theta\}}
              \bmid \cF_{k-1} \big)
    \leq \theta^{-2} \EE \big( \|\bZ^{(n)}_k\|^4 \bmid \cF_{k-1} \big)$
 \ and
 \[
   \|\bZ^{(n)}_k\|^4
   \leq 3\left( n^{-4} M_k^4 + n^{-8} M_k^4 U_{k-1}^4
                + n^{-6} M_k^4 V_{k-1}^4 \right) .
 \]
Hence
 \[
   \sum_{k=1}^{\nT}
    \EE \big( \|\bZ^{(n)}_k\|^2 \bbone_{\{\|\bZ^{(n)}_k\|>\theta\}} \big)
   \to 0
   \qquad \text{as \ $n \to \infty$ \ for all \ $\theta > 0$ \ and \ $T > 0$,}
 \]
 since \ $\EE(M_k^4) = \OO(k^2)$,
 \ $\EE( M_k^4 U_{k-1}^4 ) \leq \sqrt{\EE(M_k^8) \EE(U_{k-1}^8)} = \OO(k^6)$ \ and
 \ $ \EE( M_k^4 V_{k-1}^4 ) \leq \sqrt{\EE(M_k^8) \EE(V_{k-1}^8)} = \OO(k^4)$
 \ by Corollary \ref{EEX_EEU_EEV}.
Here we call the attention that our eight order moment condition
 \ $\EE(\vare_1^8) < \infty$ \ is used for applying Corollary
 \ref{EEX_EEU_EEV}.
This yields \eqref{Zcond3_new}.

\section{Proof of Theorem \ref{10main_Ad}}
\label{section_proof_10main}

We have
 \[
   (\tbA_n,\tbd_n)
   = \left( \sum_{k=1}^n
             \begin{bmatrix}
              n^{-3} X_{k-1}^2 & - n^{-2} X_{k-1} V_{k-1} & n^{-2} X_{k-1} \\
              - n^{-2} X_{k-1} V_{k-1} & n^{-1} V_{k-1}^2 & - n^{-1} V_{k-1} \\
              n^{-2} X_{k-1} & - n^{-1} V_{k-1} & n^{-1}
             \end{bmatrix} , \,
     \sum_{k=1}^n
      \begin{bmatrix}
       n^{-3/2} M_k X_{k-1} \\
       - n^{-1/2} M_k V_{k-1} \\
       n^{-1/2} M_k
      \end{bmatrix} \right) .
 \]
Theorem \ref{10main_Ad} will follow from the following statement (using also
 Slutsky's lemma).

\begin{Thm}\label{10main_conv}
Under the assumptions of Theorem \ref{10main} we have
 \begin{gather*}
  n^{-2} \sum_{k=1}^n X_{k-1} \as \frac{\mu}{2} , \qquad
  n^{-1} \sum_{k=1}^n V_{k-1} \as \mu , \qquad
  n^{-3} \sum_{k=1}^n X_{k-1}^2 \as \frac{\mu^2}{3} , \\
  n^{-2} \sum_{k=1}^n X_{k-1} V_{k-1} \as \frac{\mu^2}{2} , \qquad
  n^{-1} \sum_{k=1}^n V_{k-1}^2 \as  \mu^2 + \sigma^2 ,
 \end{gather*}
 and
 \[
   \sum_{k=1}^n
    \begin{bmatrix}
     n^{-3/2} M_k X_{k-1} \\
     - n^{-1/2} M_k V_{k-1} \\
     n^{-1/2} M_k
   \end{bmatrix}
   \distr
   \cN_3\left( \begin{bmatrix} 0 \\ 0 \\ 0 \end{bmatrix} , \,
             \sigma^2
             \begin{bmatrix}
              \frac{1}{3} \mu^2 &  - \frac{1}{2} \mu^2 &  \frac{1}{2} \mu \\
              - \frac{1}{2} \mu^2 &  \mu^2 + \sigma^2 &  - \mu \\
              \frac{1}{2} \mu &  - \mu & 1
             \end{bmatrix} \right) .
 \]
\end{Thm}

\noindent
\textbf{Proof.} \
In this case equation \eqref{INAR2} has the form \ $X_k = X_{k-1} + \vare_k$,
 \ $k \in \NN$, \ and hence \ $X_k = \vare_1 + \cdots + \vare_k$,
 \ $M_k = X_k - X_{k-1} - \mu = \vare_k - \mu$ \ and
 \ $V_k = X_k - X_{k-1} = \vare_k$, \ $k \in \NN$.

The first statement follows from \eqref{SLLN1} by Toeplitz theorem,
 where we used that 
 \[
   \lim_{n\to\infty}\sum_{k=1}^n\frac{k}{n^2} = \frac{1}{2}.
 \]
Again by \eqref{SLLN1},
 \[
   n^{-1} \sum_{k=1}^n V_k = n^{-1} X_k \as \mu
  \qquad \text{as \ $n \to \infty$.}
 \]
We have already shown the third statement, see \eqref{SLLN2}.
By the strong law of large numbers we have
 \begin{equation}\label{SLLN3}
  n^{-1} \sum_{k=1}^n V_k^2
  = n^{-1} \sum_{k=1}^n \vare_k^2 \as \sigma^2 + \mu^2
  \qquad \text{as \ $n \to \infty$.}
 \end{equation}
Moreover,
 \begin{align*}
  \sum_{k=1}^n X_{k-1} V_{k-1}
  &= \sum_{k=1}^n X_{k-1} \vare_{k-1}
   = \sum_{k=1}^n \vare_{k-1}^2 + \sum_{k=1}^n X_{k-2} \vare_{k-1}
   = \sum_{k=1}^n \vare_{k-1}^2 + \sum_{k=1}^n \vare_{k-1}\sum_{i=1}^{k-2} \vare_i \\
  &= \sum_{k=1}^n \vare_{k-1}^2 + \sum_{1 \leq i < j \leq n-1} \vare_i \vare_j
   = \frac{1}{2}
     \left( \left( \sum_{k=1}^n \vare_{k-1} \right)^2
                    + \sum_{k=1}^n \vare_{k-1}^2 \right)
 \end{align*}
  with \ $\varepsilon_0:=0$, \ and hence by \eqref{SLLN1} and \eqref{SLLN3},
 \begin{equation}\label{SLLN4}
  n^{-2} \sum_{k=1}^n X_{k-1} V_{k-1} \as \frac{\mu^2}{2}.
 \end{equation}
The last statement can be proved by the multidimensional martingale central
 limit theorem (see, e.g., Jacod and Shiryaev
 \cite[Chapter VIII, Theorem 3.33]{JSh}) for the sequence
 \ $(\bY^{(n)}_{k}, \cF_k)_{k \in \NN}$, \ $n \in \NN$, \ of square-integrable
 martingale differences given by
 \[
   \bY^{(n)}_k := \begin{bmatrix}
                  n^{-3/2} M_k X_{k-1} \\
                  - n^{-1/2} M_k V_{k-1} \\
                  n^{-1/2} M_k
                 \end{bmatrix}
               = \begin{bmatrix}
                  n^{-3/2} (\vare_k - \mu) X_{k-1} \\
                  - n^{-1/2} (\vare_k - \mu) \vare_{k-1} \\
                  n^{-1/2} (\vare_k - \mu)
                 \end{bmatrix} , \qquad n, k \in \NN .
 \]
We have
 \[
   \EE(\bY^{(n)}_{k} (\bY^{(n)}_{k})^\top \mid \cF_{k-1})
   = \sigma^2
     \begin{bmatrix}
      n^{-3} X_{k-1}^2 & - n^{-2} X_{k-1} \vare_{k-1} & n^{-2} X_{k-1} \\
      - n^{-2} X_{k-1} \vare_{k-1} & n^{-1} \vare_{k-1}^2 & - n^{-1} \vare_{k-1} \\
      n^{-2} X_{k-1} & - n^{-1} \vare_{k-1} & n^{-1}
     \end{bmatrix} , \qquad n, k \in \NN ,
 \]
 hence by \eqref{SLLN2}, \eqref{SLLN3}, and \eqref{SLLN4}, we have the
 asymptotic covariance matrix
 \[
   \sum_{k=1}^n \EE(\bY^{(n)}_{k} (\bY^{(n)}_{k})^\top \mid \cF_{k-1})
   \as \sigma^2
       \begin{bmatrix}
        \frac{\mu^2}{3} &  - \frac{\mu^2}{2} & \frac{\mu}{2} \\
        - \frac{\mu^2}{2} & \sigma^2 + \mu^2 & - \mu \\
        \frac{\mu}{2} & - \mu & 1
       \end{bmatrix} 
       \qquad \text{as \ $n\to\infty$.}
 \]
The conditional Lindeberg condition
 \[
   \sum_{k=1}^n
    \EE( \|\bY^{(n)}_k\|^2 \bbone_{\{\|\bY^{(n)}_{k}\|>\theta\}} \mid \cF_{k-1} )
   \stoch 0 \qquad \text{as \ $n\to\infty$}
 \]
 is satisfied for all \ $\theta > 0$, \ since using that
 \ $\EE(\vare_1^4) < \infty$,
 \begin{align*}
  &\sum_{k=1}^n
   \EE( \|\bY^{(n)}_{k}\|^2 \bbone_{\{\|\bY^{(n)}_{k}\|>\theta\}} \mid \cF_{k-1} )
   \leq \frac{1}{\theta^2} \sum_{k=1}^n \EE(\|\bY^{(n)}_k\|^4 \mid \cF_{k-1}) \\
  &\leq \frac{3}{\theta^2}
        \sum_{k=1}^n \EE\bigl( n^{-6} (\vare_k - \mu)^4 X_{k-1}^4
                             + n^{-2} (\vare_k - \mu)^4 (\vare_{k-1}^4 + 1)
                             \mid \cF_{k-1} \bigr) \\
  &\leq \frac{3 \EE\big((\vare_1 - \mu)^4\big)}{\theta^2}
        \sum_{k=1}^n \bigl( n^{-6} X_{k-1}^4 + n^{-2} (\vare_{k-1}^4 + 1) \bigr)
   \stoch 0
 \end{align*}
 as \ $n \to \infty$, \ where the last but one step follows by that \ $\vare_k$ \ and \ $\cF_{k-1}$ \ are
 independent, \ $\vare_{k-1}$ \ is measurable with respect to the \ $\sigma$-algebra \ $\cF_{k-1}$ \ 
 (since \ $\vare_{k-1} = X_{k-1} - X_{k-2}$), \ and the last step follows by
 \ $\EE(X_k^4) = \OO(k^4)$ \ (see Corollary \ref{EEX_EEU_EEV}).
\proofend

\section{Proof of Theorem \ref{01main_Ad}}
\label{section_proof_01main}

We have
 \[
   \tbA_n = \sum_{k=1}^n
             \begin{bmatrix}
              n^{-3} X_{k-1}^2 & - n^{-5/2} X_{k-1} V_{k-1} & n^{-2} X_{k-1} \\
              - n^{-5/2} X_{k-1} V_{k-1} & n^{-2} V_{k-1}^2 & - n^{-3/2} V_{k-1} \\
              n^{-2} X_{k-1} & - n^{-3/2} V_{k-1} & n^{-1}
             \end{bmatrix} , \quad
   \tbd_n = \sum_{k=1}^n
             \begin{bmatrix}
              n^{-3/2} M_k X_{k-1} \\
              - n^{-1} M_k V_{k-1} \\
              n^{-1/2} M_k
             \end{bmatrix} .
 \]

\begin{Lem}\label{XX_XV}
Under the assumptions of Theorem \ref{01main}, as $n\to\infty$, we have
 \begin{gather}
  n^{-2} \sum_{k=1}^n X_{k-1} \as \frac{\mu}{4} , \label{X_mu} \\
  n^{-1} \sum_{k=1}^n V_{k-1} \as \frac{\mu}{2} , \label{V_mu} \\
  n^{-3} \sum_{k=1}^n X_{k-1}^2 \as \frac{\mu^2}{12} , \label{XX} \\
  n^{-5/2} \sum_{k=1}^n X_{k-1} V_{k-1} \stoch 0 , \label{main_XV} \\
  n^{-2} \sum_{k=1}^n \bigl[\EE(V_{k-1})\bigr]^2 \to 0 ,
  \label{main_EVEV} \\
  n^{-2} \sum_{k=1}^n (V_{k-1}-\EE(V_{k-1}))\EE(V_{k-1}) \stoch 0 ,
  \label{main_VEV} \\
  n^{-3/2}
  \sum_{k=1}^n M_k \big(X_{k-1}-\EE(X_{k-1})\bigr) \stoch 0 , \label{main_MX} \\
  n^{-1}
  \sum_{k=1}^n M_k \EE(V_{k-1}) \stoch 0 . \label{main_MV}
 \end{gather}
\end{Lem}

\noindent
\textbf{Proof.} \
In this case equation \eqref{INAR2} has the form \ $X_k = X_{k-2} + \vare_k$,
 \ $k \in \NN$, \ and hence \ $X_{2k} = \vare_2 + \vare_4 + \cdots + \vare_{2k}$,
 \ $X_{2k-1} = \vare_1 + \vare_3 + \cdots + \vare_{2k-1}$,
 \ $M_k = X_k - X_{k-2} - \mu = \vare_k - \mu$ \ and
 \ $V_{2k} = X_{2k} - X_{2k-1}
    = (\vare_2 - \vare_1) + \cdots + (\vare_{2k} - \vare_{2k-1})$,
 \ $V_{2k-1} = X_{2k-1} - X_{2k-2}
    = (\vare_1 - \vare_2) + \cdots + (\vare_{2k-3} - \vare_{2k-2}) + \vare_{2k-1}$,
 \ $k \in \NN$.

Convergence \eqref{X_mu} follows from \eqref{SLLN1_01} by Toeplitz theorem.
Again by \eqref{SLLN1_01}, we obtain
 \[
   n^{-1} \sum_{k=1}^n V_{k-1} = n^{-1} X_{n-1}  \as \frac{\mu}{2} ,
 \]
 yielding \eqref{V_mu}.
We have already shown \eqref{XX}, see \eqref{SLLN2_01}.

In order to show \eqref{main_XV}, we use \eqref{XV_VV}.
Clearly we have
 \[
   \EE\left(\sum_{k=1}^n X_k V_k \right)
   = \frac{1}{2} \EE(X_n^2) + \frac{1}{2} \sum_{k=1}^n \EE(V_k^2)
   = \OO(n^2) ,
 \]
 since, by Corollary \ref{EEX_EEU_EEV}, \ $\EE(X_n^2) = \OO(n^2)$,
 \ $n \in \NN$, \ and \ $\EE(V_k^2) = \OO(k)$, \ $k \in \NN$, \ and hence we
 obtain \eqref{main_XV}.
For each \ $k \in \NN$, \ we have \ $\EE(V_{2k}) = 0$ \ and
 \ $\EE(V_{2k-1}) = \mu$, \ hence we conclude \eqref{main_EVEV}, and
 \begin{align*}
  \EE\left( \left| \sum_{k=1}^n (V_{k-1}-\EE(V_{k-1}))\EE(V_{k-1}) \right| \right)
  &\leq \sum_{k=1}^n \mu \EE\bigl( | V_{k-1} - \EE(V_{k-1}) | \bigr) \\
  &\leq \sum_{k=1}^n \mu \sqrt{\EE\bigl( ( V_{k-1} - \EE(V_{k-1}) )^2 \bigr)}
   = \OO(n^{3/2}) ,
 \end{align*}
 since \ $\EE(V_k - \EE(V_k))^2\leq \EE(V_k^2) = \OO(k)$, \ $k \in \NN$
 \ (by Corollary \ref{EEX_EEU_EEV}), which implies \eqref{main_VEV}.

Moreover, using that \ $M_k(X_{k-1} - \EE(X_{k-1}))$, \ $k \in \{1, \ldots, n\}$,
 \ are uncorrelated,
 \begin{align*}
  \EE\left( \left( \sum_{k=1}^n M_k \big(X_{k-1}-\EE(X_{k-1})\bigr) \right)^2
      \right)
  &= \sum_{k=1}^n
      \EE\left( (\vare_k - \mu)^2 \big(X_{k-1}-\EE(X_{k-1})\bigr)^2 \right)\\
  &= \sigma^2 \sum_{k=1}^n \EE\left( \big(X_{k-1}-\EE(X_{k-1})\bigr)^2 \right)
   = \OO(n^2) ,
 \end{align*}
 since \ $\EE(X_{k-1} - \EE(X_{k-1}))^2 \leq \lfloor k/2\rfloor \sigma^2$,
 \ $k \in \NN$ \ (by Corollary \ref{EEX_EEU_EEV}), thus we get \eqref{main_MX}.

Since \ $\EE(V_{2k}) = 0$ \ and \ $\EE(V_{2k-1}) = \mu$, \ $k \in \NN$,
 \ we have
 \begin{align*}
  \EE\left( \left( \sum_{k=1}^n M_k \EE(V_{k-1}) \right)^2 \right)
  = \sum_{k=1}^n
     \bigl[\EE(V_{k-1})\bigr]^2 \EE\bigl((\vare_k - \mu)^2\bigr)
  = \sigma^2 \sum_{k=1}^n \bigl[\EE(V_{k-1})\bigr]^2
  = \OO(n) ,
 \end{align*}
 which implies \eqref{main_MV}.
\proofend

Theorem \ref{01main_Ad} will follow from Lemma \ref{XX_XV} and the following
 statement (using Slutsky's lemma).

\begin{Thm}\label{01main_conv}
Under the assumptions of Theorem \ref{01main} we have
 \[
   \sum_{k=1}^n
    \begin{bmatrix}
     n^{-2} \big(V_{k-1}-\EE(V_{k-1})\bigr)^2 \\
     n^{-3/2} M_k \EE(X_{k-1}) \\
     - n^{-1} M_k \big(V_{k-1} - \EE(V_{k-1})\bigr) \\
     n^{-1/2} M_k
    \end{bmatrix}
   \distr
   \begin{bmatrix}
    \sigma^2 \int_0^1 \cW_t^2 \, \dd t \\
    \frac{1}{2} \mu \sigma \int_0^1 t \, \dd \tcW_t \\
    \sigma^2 \int_0^1 \cW_t \, \dd \cW_t \\
    \sigma \tcW_1
   \end{bmatrix}
   \qquad \text{as \ $n\to\infty$},
 \]
 where \ $(\cW_t)_{t\in\RR_+}$ \ and \ $(\tcW_t)_{t\in\RR_+}$ \ are independent
 standard Wiener processes.
\end{Thm}

\noindent
\textbf{Proof.} \
Consider the sequence
 \[
   \begin{bmatrix} \cS_t^{(n)} \\ \cT_t^{(n)} \end{bmatrix}
   := \begin{bmatrix}
       n^{-1/2} \bigl(X_{2\nt} - \EE(X_{2\nt})\bigr) \\
       n^{-1/2} \bigl(X_{2\nt-1} - \EE(X_{2\nt-1})\bigr)
      \end{bmatrix} , \qquad t \in \RR_+ , \quad n \in \NN ,
 \]
 of stochastic processes.
Then, by the multidimensional martingale central limit theorem,
 \begin{equation}\label{Donsker}
   \begin{bmatrix} \cS^{(n)} \\ \cT^{(n)} \end{bmatrix}
   \distr
   \sigma
   \begin{bmatrix} \cB \\ \tcB \end{bmatrix} \qquad
   \text{as \ $n \to \infty$,}
 \end{equation}
 where \ $(\cB_t)_{t\in\RR_+}$ \ and \ $(\tcB_t)_{t\in\RR_+}$ \ are independent
 standard Wiener processes.
Indeed, with the notation
 \[
   \bY^{(n)}_k := \begin{bmatrix}
                  n^{-1/2}(\vare_{2k} - \mu) \\
                  n^{-1/2}(\vare_{2k-1} - \mu) \\
                 \end{bmatrix} ,
   \qquad n, k \in \NN ,
 \]
 we have that \ $(\bY^{(n)}_k, \cF_{2k})_{k\in\NN}$, \ $n \in \NN$, \ is a
 sequence of square integrable martingale differences such that
 \[
   \sum_{k=1}^\nt \bY^{(n)}_k
   = \begin{bmatrix}
      \cS^{(n)}_t \\
      \cT^{(n)}_t \\
     \end{bmatrix} ,
   \qquad n \in \NN , \quad t \in \RR_+ ,
 \]
 $\EE\big(\bY^{(n)}_k \mid \cF_{2(k-1)}\big) = \bzero \in \RR^2$ \ and
 \[
   \EE \big( \bY^{(n)}_k (\bY^{(n)}_k)^\top \mid \cF_{2(k-1)} \big)
   = \sigma^2
     n^{-1}
     \bI_2 ,
  \qquad n, k \in \NN ,
 \]
 where \ $\bI_2$ \ denotes the $2\times 2$ identity matrix.
Then the asymptotic covariance matrix
 \[
   \sum_{k=1}^\nt
    \EE \big( \bY^{(n)}_k (\bY^{(n)}_k)^\top \mid \cF_{2(k-1)} \big)
   \as \sigma^2 t \bI_2
   \qquad \text{as \ $n \to \infty$ \ for \ $t \in \RR_+$.}
 \]
The conditional Lindeberg condition
 \begin{align}\label{01seged4}
  \sum_{k=1}^\nt
   \EE\big(\|\bY^{(n)}_k\|^2 \bbone_{\{\|\bY^{(n)}_k\|>\theta\}} \mid \cF_{2(k-1)}\big)
  \stoch 0 \qquad \text{as \ $n \to \infty$}
 \end{align}
 is satisfied for all \ $t \in \RR_+$ \ and \ $\theta > 0$.
\ Indeed, we have
 \begin{align*}
  \sum_{k=1}^\nt
  & \EE \big( \|\bY^{(n)}_k\|^2 \bbone_{\{\|\bY^{(n)}_k\|>\theta\}} \big) \\
  &= \frac{1}{n}
     \sum_{k=1}^\nt
      \EE\Big(\left( (\vare_{2k} - \mu)^2 + (\vare_{2k-1} - \mu)^2 \right)
              \bbone_{\{(\vare_{2k}-\mu)^2+(\vare_{2k-1}-\mu)^2>n\theta^2\}}
         \Big) \\
  &= \frac{\nt}{n}
     \EE\Big(\left( (\vare_2 - \mu)^2 + (\vare_1 - \mu)^2 \right)
             \bbone_{\{(\vare_2-\mu)^2+(\vare_1-\mu)^2>n\theta^2\}}\Big)
     \to 0 ,
   \end{align*}
 by dominated convergence theorem.
This yields that the convergence in \eqref{01seged4} holds in fact in
 $L_1$-sense.
Thus we obtain \eqref{Donsker}.
We are going to prove that convergence \eqref{Donsker} implies
 \begin{equation}\label{Ito}
   \sum_{k=1}^n
    \begin{bmatrix}
     n^{-2} \big(V_{k-1} - \EE(V_{k-1})\bigr)^2 \\
     n^{-3/2} M_k \EE(X_{k-1}) \\
     - n^{-1} M_k \big(V_{k-1} - \EE(V_{k-1})\bigr) \\
     n^{-1/2} M_k
    \end{bmatrix}
   \distr
   \begin{bmatrix}
    \frac{1}{2} \sigma^2 \int_0^1 (\cB_t - \tcB_t)^2 \, \dd t \\
    \frac{1}{2^{3/2}} \mu \sigma
    \left( \cB_1 + \tcB_1 - \int_0^1 (\cB_t + \tcB_t) \, \dd t \right) \\
    \frac{1}{4} \sigma^2 \bigl[ (\cB_1 - \tcB_1)^2 - 2 \bigr] \\
    \frac{1}{2^{1/2}} \sigma (\cB_1 + \tcB_1)
   \end{bmatrix}
 \end{equation}
 as $n\to\infty$, \ which yields the statement.
Indeed, \ $\bigl(2^{-1/2}(\cB_t + \tcB_t)\bigr)_{t\in\RR_+}$ \ and
 \ $\bigl(2^{-1/2}(\cB_t - \tcB_t)\bigr)_{t\in\RR_+}$ \ are independent
 standard Wiener processes, and by It\^o's formula,
 \ $\int_0^1 t \, \dd \tcW_t = \tcW_1 - \int_0^1 \tcW_t \, \dd t$ \ and
 \ $\int_0^1 \cW_t \, \dd \cW_t = 2^{-1} (\cW_1^2 - 1)$, \ which yield the
 statement with the choices \ $\tcW_t := 2^{-1/2} (\cB_t + \tcB_t)$,
 \ $t \geq 0$, \ and \ $\cW_t := 2^{-1/2} (\cB_t - \tcB_t)$, \ $t \geq 0$.

Applying Lemmas \ref{Conv2Funct} and \ref{Marci} as in the proof of Proposition
 \ref{ExUn} and using Slutsky's lemma, \eqref{Ito} will follow from
 \begin{gather}
  \frac{1}{n^2} \sum_{k=1}^n \bigl(V_{k-1}-\EE(V_{k-1})\bigr)^2
  - \frac{1}{n}
    \sum_{k=1}^{\nhalf} \Bigl(\cS_{2k/n}^{(\nhalf)} - \cT_{2k/n}^{(\nhalf)}\Bigr)^2
  \stoch 0 , \label{VV_ST} \\
  \frac{1}{n^{3/2}} \sum_{k=1}^n M_k \EE(X_{k-1})
  - \frac{\mu}{2^{3/2}}
    \biggl( \cS_1^{(\nhalf)} + \cT_1^{(\nhalf)}
            - \frac{2}{n}
              \sum_{k=1}^{\nhalf}
               \Bigl(\cS_{2k/n}^{(\nhalf)} + \cT_{2k/n}^{(\nhalf)}\Bigr) \biggr)
  \stoch 0 , \label{MX_ST} \\
  \frac{1}{n} \sum_{k=1}^n M_k \bigl(V_{k-1}-\EE(V_{k-1})\bigr)
      + \frac{1}{4}
    \left[ (\cS_1^{(\nhalf)} - \cT_1^{(\nhalf)})^2 - 2 \sigma^2 \right]
  \stoch 0 , \label{MV_ST} \\
  \frac{1}{n^{1/2}} \sum_{k=1}^n M_k
  - \frac{1}{2^{1/2}} \bigl( \cS_1^{(\nhalf)} + \cT_1^{(\nhalf)} \bigr)
  \stoch 0 . \label{M_ST}
 \end{gather}
Indeed, first considering the subsequence \ $(2n)_{n\in\NN}$, \ let us apply
 Lemmas \ref{Conv2Funct} and \ref{Marci} with the special choices \ $d := 2$,
 \ $p := 2$, \ $q := 2$, \ $h : \RR^2 \to \RR^2$,
 \[
   h(x_1, x_2)
   := \left(x_1 + x_2, \frac{1}{4} (x_1 - x_2)^2 - \frac{\sigma^2}{2}\right) ,
   \qquad (x_1, x_2) \in \RR^2 ,
 \]
 $K : [0,1] \times \RR^4 \to \RR^2$,
 \[
   K(s, x_1, x_2, x_3, x_4)
   := \left( \frac{1}{2} (x_1 - x_2)^2 , x_1 + x_2 \right) ,
   \qquad (s, x_1, x_2, x_3, x_4) \in [0, 1] \times \RR^4 ,
 \]
 and
 \[
   \cU := \sigma \begin{bmatrix}
                  \cB \\
                  \tcB
                 \end{bmatrix},
   \qquad
   \cU^{(n)} := \begin{bmatrix}
                \cS^{(n)} \\
                \cT^{(n)}
               \end{bmatrix} ,
   \qquad n \in \NN .
 \]
Then
 \begin{align*}
  &\|K(s, x_1, x_2, x_3, x_4) - K(t, y_1, y_2, y_3, y_4)\| \\
  &= \left( \frac{1}{4} \big( (x_1 - x_2)^2 - (y_1 - y_2)^2 \big)^2
            + (x_1 - y_1 + x_2 - y_2)^2 \right)^{1/2}\\
  &= \left( \frac{1}{4}(x_1-y_1 + y_2-x_2)^2(x_1-x_2 + y_1 - y_2)^2
            + (x_1 - y_1 + x_2 - y_2)^2 \right)^{1/2} \\
  &\leq 2 \left( \bigl( (x_1 - y_1)^2 + (y_2 - x_2)^2 \bigr)
                 \bigl( (x_1 - x_2)^2 + (y_1 - y_2)^2 \bigr)
                 + \bigl( (x_1 - y_1)^2 + (x_2 - y_2)^2 \bigr) \right)^{1/2} \\
  &\leq 8 R \, \|(x_1, x_2, x_3, x_4) - (y_1, y_2, y_3, y_4)\|
 \end{align*}
 for all \ $s, t \in [0, 1]$ \ and
 \ $(x_1, x_2, x_3, x_4), (y_1, y_2, y_3, y_4) \in \RR^4$ \ with
 \ $\|(x_1, x_2, x_3, x_4)\| \leq R$ \ and \ $\|(y_1, y_2, y_3, y_4)\| \leq R$,
 \ where \ $R > 0$, \ since
 \[
   (x_1 - x_2)^2 + (y_1 - y_2)^2
   \leq 2(x_1^2 + x_2^2) + 2(y_1^2 + y_2^2)
   \leq 8 R^2 .
 \]
Further, using the definitions of \ $\Phi$ \ and \ $\Phi_n$, \ $n \in \NN$,
 \ given in Lemma \ref{Marci}, we have
 \begin{align*}
  & \Phi_n\left(\begin{bmatrix}
                 \cS^{(n)} \\
                 \cT^{(n)} \\
                \end{bmatrix} \right) \\
  &\qquad
    = \left( \cS^{(n)}_1 + \cT^{(n)}_1 ,
             \frac{1}{4} (\cS^{(n)}_1 - \cT^{(n)}_1)^2 - \frac{\sigma^2}{2} ,
             \frac{1}{n}
             \sum_{k=1}^n \frac{1}{2} (\cS^{(n)}_{k/n} - \cT^{(n)}_{k/n})^2 ,
             \frac{1}{n} \sum_{k=1}^n (\cS^{(n)}_{k/n} + \cT^{(n)}_{k/n}) \right)
 \end{align*}
 and
 \begin{align*}
  \Phi\left( \sigma
             \begin{bmatrix}
              \cB \\
              \widetilde\cB
             \end{bmatrix} \right)
  = \left( \sigma (\cB_1 + \tcB_1) ,
           \frac{\sigma^2}{4} (\cB_1 - \tcB_1)^2 - \frac{\sigma^2}{2} ,
           \int_0^1 \frac{\sigma^2}{2} (\cB_u - \tcB_u)^2 \, \dd u ,
           \int_0^1 \sigma (\cB_u + \tcB_u) \, \dd u \right) .
 \end{align*}
Since the process \ $\sigma [\cB_t \; \tcB_t]^\top_{t\in\RR_+}$ \ admits
 continuous paths with probability one, \eqref{Donsker},
 Lemma \ref{Conv2Funct} (with the choice \ $C := \CC(\RR_+, \RR^2)$), \ and
 Lemma \ref{Marci} yield that
 \[
   \Phi_n\left( \begin{bmatrix}
                 \cS^{(n)} \\
                 \cT^{(n)}
                \end{bmatrix} \right)
   \distr
   \Phi\left( \sigma
              \begin{bmatrix}
               \cB \\
               \tcB
              \end{bmatrix} \right)
   \qquad \text{as \ $n \to \infty$.}
 \]
By another easy application of continuous mapping theorem (one can again apply
 Lemmas \ref{Conv2Funct} and \ref{Marci}) we have
 \begin{align*}
  \begin{bmatrix}
   \frac{1}{2n} \sum_{k=1}^n (\cS^{(n)}_{k/n} - \cT^{(n)}_{k/n})^2 \\
   \cS^{(n)}_1 + \cT^{(n)}_1
    - \frac{1}{n} \sum_{k=1}^n (\cS^{(n)}_{k/n} + \cT^{(n)}_{k/n}) \\
   \frac{1}{4} (\cS^{(n)}_1 - \cT^{(n)}_1)^2 - \frac{\sigma^2}{2} \\
   \cS^{(n)}_1 + \cT^{(n)}_1
  \end{bmatrix}
  \distr
  \begin{bmatrix}
   \frac{\sigma^2}{2} \int_0^1 (\cB_u - \tcB_u)^2 \, \dd u \\
   \sigma \left(\cB_1 + \tcB_1 - \int_0^1 (\cB_u + \tcB_u) \, \dd u \right) \\
   \frac{\sigma^2}{4} \left( (\cB_1 - \tcB_1)^2 -2 \right)  \\
   \sigma(\cB_1 + \tcB_1)
  \end{bmatrix}
  \qquad \text{as \ $n \to \infty$.}
\end{align*}
Hence, using \eqref{VV_ST}, \eqref{MX_ST}, \eqref{MV_ST}, and Slutsky's lemma,
 we have \eqref{Ito} for the subsequence \ $(2n)_{n\in\NN}$.
\ To prove \eqref{Ito} for the subsequence \ $(2n - 1)_{n\in\NN}$,
 \ by Slutsky's lemma, it is enough to check that as \ $n\to\infty$,
 \begin{gather}
  \frac{1}{n^2} (V_n - \EE(V_n))^2  \stoch 0, \label{seged_paratlan1} \\
  \frac{1}{n^{3/2}} M_n \EE(X_{n-1}) \stoch 0, \label{seged_paratlan2} \\
  \frac{1}{n} M_n (V_{n-1} - \EE(V_{n-1})) \stoch 0 \label{seged_paratlan3} \\
  \frac{1}{n^{1/2}} M_n  \stoch 0 \label{seged_paratlan4}.
 \end{gather}
By Corollary \ref{EEX_EEU_EEV}, \ $\EE((V_n - \EE(V_n))^2) = \OO(n)$,
 \[
   \EE(|M_n \EE(X_{n-1})|) = \EE(|M_n|) \EE(X_{n-1})
   = \EE(|\vare_1 - \mu|) \EE(X_{n-1}) = \OO(n) ,
 \]
 \begin{align*}
  \EE(M_n^2 (V_{n-1} - \EE(V_{n-1}))^2)
  &= \EE\big((V_{n-1} - \EE(V_{n-1}))^2 \EE(M_n^2 \mid \cF_{n-1}) \big) \\
  &= \sigma^2 \EE((V_{n-1} - \EE(V_{n-1}))^2)
   = \OO(n) ,
 \end{align*}
 and 
 \[
   \EE(M_n^2) = \EE((\vare_n - \mu)^2) = \sigma^2,\qquad n\in\NN,
 \]
 thus we obtain \eqref{seged_paratlan1}, \eqref{seged_paratlan2}, \eqref{seged_paratlan3},
 and \eqref{seged_paratlan4}.

First we will prove \eqref{VV_ST}, \eqref{MX_ST},
 \eqref{MV_ST} and \eqref{M_ST} for the
 subsequence \ $(2n)_{n\in\NN}$ \ and then the subsequence \ $(2n-1)_{n\in\NN}$.
\ In order to prove \eqref{VV_ST} first observe that, for all \ $k \in \NN$,
 \begin{align*}
  &V_{2k} - \EE(V_{2k})
   = (X_{2k} - \EE(X_{2k})) - (X_{2k-1} - \EE(X_{2k-1}))
   = n^{1/2} (\cS_{k/n}^{(n)} - \cT_{k/n}^{(n)}) , \\
  &V_{2k-1} - \EE(V_{2k-1})
   = (\vare_{2k} - \EE(\vare_{2k})) - (V_{2k} - \EE(V_{2k}))
   = (\vare_{2k} - \mu) - n^{1/2} (\cS_{k/n}^{(n)} - \cT_{k/n}^{(n)}) .
 \end{align*}
Then
 \begin{align*}
  &\frac{1}{(2n)^2} \sum_{k=1}^{2n} \big(V_{k-1}-\EE(V_{k-1})\bigr)^2
   = \frac{1}{4n^2} \sum_{k=1}^{n-1} \big(V_{2k}-\EE(V_{2k})\bigr)^2
     + \frac{1}{4n^2} \sum_{k=1}^n \big(V_{2k-1}-\EE(V_{2k-1})\bigr)^2 \\
  &= \frac{1}{4n} \sum_{k=1}^{n-1} (\cS_{k/n}^{(n)} - \cT_{k/n}^{(n)})^2
     + \frac{1}{4n^2} \sum_{k=1}^n
        \big[(\vare_{2k} - \mu)
             - n^{1/2} (\cS_{k/n}^{(n)} - \cT_{k/n}^{(n)}) \bigr]^2 \\
  &= \frac{1}{2n} \sum_{k=1}^n (\cS_{k/n}^{(n)} - \cT_{k/n}^{(n)})^2
     - \frac{1}{4n} (\cS_1^{(n)} - \cT_1^{(n)})^2 \\
  &\quad
     - \frac{1}{2n^{3/2}}
       \sum_{k=1}^n (\vare_{2k} - \mu) (\cS_{k/n}^{(n)} - \cT_{k/n}^{(n)})
     + \frac{1}{4n^2} \sum_{k=1}^n (\vare_{2k} - \mu)^2 \\
  &= \frac{1}{2n} \sum_{k=1}^n (\cS_{k/n}^{(n)} - \cT_{k/n}^{(n)})^2
     - \frac{1}{4n^2} (V_{2n} - \EE(V_{2n}))^2 \\
  &\quad
     - \frac{1}{2n^2}
       \sum_{k=1}^n (\vare_{2k} - \mu) (V_{2k} - \EE(V_{2k}))
     + \frac{1}{4n^2} \sum_{k=1}^n (\vare_{2k} - \mu)^2 .
 \end{align*}
Thus, in order to prove \eqref{VV_ST} for the subsequence \ $(2n)_{n\in \NN}$,
 \ it suffices to prove
 \begin{gather}
  \frac{1}{n^2} (V_{2n} - \EE(V_{2n}))^2 \stoch 0 , \label{VV_ST_ST} \\
  \frac{1}{n^2}
  \sum_{k=1}^n (\vare_{2k} - \mu) (V_{2k} - \EE(V_{2k}))
  \stoch 0 , \label{VV_ST_EST} \\
  \frac{1}{n^2} \sum_{k=1}^n (\vare_{2k} - \mu)^2 \stoch 0  \label{VV_ST_EE}
 \end{gather}
 as $n\to\infty$.
By Corollary \ref{EEX_EEU_EEV}, we have
 \ $\EE\bigl((V_{2n} - \EE(V_{2n}))^2\bigr) = \OO(n)$ \ and
 \ $\EE\bigl((\vare_{2k} - \mu)^2\bigr) = \sigma^2$, \ thus we obtain
 \eqref{VV_ST_ST} and \eqref{VV_ST_EE}.
Further,
 \ $V_{2k} - \EE(V_{2k}) = (\vare_{2k} - \mu) - (V_{2k-1} - \EE(V_{2k-1}))$,
 \ hence \eqref{VV_ST_EST} follows from \eqref{VV_ST_EE} and from
 \begin{align*}
  \EE\left( \left( \sum_{k=1}^n (\vare_{2k} - \mu) (V_{2k-1} - \EE(V_{2k-1}))
            \right)^2
     \right)
  &= \sum_{k=1}^n
      \EE\left( (\vare_{2k} - \mu)^2 (V_{2k-1}-\EE(V_{2k-1}))^2 \right)\\
  &= \sigma^2 \sum_{k=1}^n \EE\left( (V_{2k-1}-\EE(V_{2k-1}))^2 \right)
   = \OO(n^2) ,
 \end{align*}
 and we finish the proof of \eqref{VV_ST} for the subsequence
 \ $(2n)_{n\in \NN}$.

Now we turn to prove \eqref{MX_ST} for the subsequence \ $(2n)_{n\in \NN}$.
\ First observe that
 \begin{align*}
  \sum_{k=1}^{2n} M_k \EE(X_{k-1})
  = \mu \sum_{k=1}^{n-1} (\vare_{2k+1} - \mu) k
    + \mu \sum_{k=1}^n (\vare_{2k} - \mu) k .
 \end{align*}
We have
 \begin{align*}
  \sum_{k=1}^n (\vare_{2k} - \mu) k
  &= \sum_{k=1}^n \sum_{j=1}^k (\vare_{2k} - \mu)
   = \sum_{j=1}^n \sum_{k=j}^n (\vare_{2k} - \mu)
   = \sum_{j=1}^n
      \left( \sum_{k=1}^n (\vare_{2k} - \mu)
             - \sum_{k=1}^{j-1} (\vare_{2k} - \mu) \right) \\
  &= \sum_{j=1}^n
      \bigl[ (X_{2n} - \EE(X_{2n})) - (X_{2j-2} - \EE(X_{2j-2})) \bigr]
   = n^{3/2} \cS_1^{(n)} - n^{1/2} \sum_{j=1}^n \cS_{(j-1)/n}^{(n)} ,
 \end{align*}
 and, in a similar way,
\begin{align*}
 \sum_{k=1}^{n-1} (\vare_{2k+1} - \mu) k
  & = \sum_{k=1}^{n-1} \sum_{j=1}^k (\vare_{2k+1} - \mu)
   = \sum_{j=1}^{n-1} \sum_{k=j}^{n-1} (\vare_{2k+1} - \mu) \\
  & = \sum_{j=1}^{n-1} \left( \sum_{k=1}^{n-1}(\vare_{2k+1} - \mu)
                             - \sum_{k=1}^{j-1}(\vare_{2k+1} - \mu) \right) \\
  & = \sum_{j=1}^{n-1} \Big( X_{2n-1} - \vare_1 - (n-1)\mu
                             - (X_{2j-1} - \vare_1 - (j-1)\mu) \Big) \\
  & = \sum_{j=1}^{n-1} \Big( (X_{2n-1} - \EE(X_{2n-1}) - \vare_1 + \mu)
                           - (X_{2j-1} - \EE(X_{2j-1})  - \vare_1 + \mu) \Big) \\
  & = \sum_{j=1}^{n-1}
       \Big( X_{2n-1} - \EE(X_{2n-1}) - (X_{2j-1} - \EE(X_{2j-1})) \Big) \\
  & = n^{1/2} (n-1) \cT_1^{(n)} - n^{1/2} \sum_{j=1}^{n-1} \cT_{j/n}^{(n)}
    = n^{3/2} \cT_1^{(n)} - n^{1/2} \sum_{j=1}^n \cT_{j/n}^{(n)} .
 \end{align*}
Hence
 \begin{align*}
  \frac{1}{(2n)^{3/2}} \sum_{k=1}^{2n} M_k \EE(X_{k-1})
  &= \frac{\mu}{2^{3/2}}
     \biggl( \cS_1^{(n)} + \cT_1^{(n)}
             - \frac{1}{n}
               \sum_{k=1}^{n} \Bigl(\cS_{k/n}^{(n)} + \cT_{k/n}^{(n)}\Bigr) \biggr)
     + \frac{\mu}{2^{3/2} n} \cS_1^{(n)} .
 \end{align*}
Convergence \eqref{Donsker} implies \ $\cS_1^{(n)} \distr \sigma \cB_1$ \ and
 hence \ $n^{-1} \cS_1^{(n)} \stoch 0$, \ thus we obtain \eqref{MX_ST} for the
 subsequence \ $(2n)_{n\in\NN}$.

Now we turn to prove \eqref{MV_ST} for the subsequence \ $(2n)_{n\in \NN}$.
\ First observe that
 \begin{align*}
  \sum_{k=1}^{2n} M_k \bigl(V_{k-1}-\EE(V_{k-1})\bigr)
  &= \sum_{k=1}^{n-1}
      (\vare_{2k+1} - \mu)
      \bigl[ (X_{2k}-\EE(X_{2k})) - (X_{2k-1}-\EE(X_{2k-1})) \bigr] \\
  &\quad
     + \sum_{k=1}^n
        (\vare_{2k} - \mu)
        \bigl[ (X_{2k-1}-\EE(X_{2k-1})) - (X_{2k-2}-\EE(X_{2k-2})) \bigr] \\
  &= \sum_{k=1}^{n-1} (\vare_{2k+1} - \mu) \sum_{j=1}^k (\vare_{2j}-\mu)
     - \sum_{k=1}^{n-1} (\vare_{2k+1}-\mu) \sum_{j=1}^k (\vare_{2j-1}-\mu)  \\
  &\quad
    + \sum_{k=1}^n (\vare_{2k} - \mu) \sum_{j=1}^k (\vare_{2j-1}-\mu)
    - \sum_{k=1}^n (\vare_{2k} - \mu) \sum_{j=1}^{k-1} (\vare_{2j}-\mu) .
 \end{align*}
Here the sum of the first and third summands is
 \begin{align*}
  &\sum_{k=1}^{n-1} (\vare_{2k+1} - \mu) \sum_{j=1}^k (\vare_{2j} - \mu)
   + \sum_{k=1}^n (\vare_{2k} - \mu) \sum_{j=1}^k (\vare_{2j-1} - \mu) \\
  &\qquad
   = \sum_{k=2}^n \sum_{j=1}^{k-1} (\vare_{2k-1} - \mu) (\vare_{2j} - \mu)
     + \sum_{j=1}^n \sum_{k=1}^j  (\vare_{2j} - \mu)  (\vare_{2k-1} - \mu) \\
  &\qquad
   = \sum_{j=1}^{n} \sum_{k=j+1}^{n} (\vare_{2k-1} - \mu)  (\vare_{2j} - \mu)
     + \sum_{j=1}^n \sum_{k=1}^j (\vare_{2j} - \mu)  (\vare_{2k-1} - \mu) \\
  &\qquad
   = \sum_{j=1}^n \sum_{k=1}^n (\vare_{2k-1} - \mu) (\vare_{2j} - \mu) \\
  &\qquad
   = \sum_{j=1}^n (\vare_{2j} - \mu) \sum_{k=1}^n (\vare_{2k-1} - \mu)
   = n \cS_1^{(n)} \cT_1^{(n)} ,
 \end{align*}
 the second summand is
 \begin{multline*}
  \sum_{k=1}^{n-1} (\vare_{2k+1} - \mu) \sum_{j=1}^k (\vare_{2j-1} - \mu)
  = \sum_{1 \leq j < \ell \leq n} (\vare_{2j-1} - \mu) (\vare_{2\ell-1} - \mu) \\
  = \frac{1}{2}
    \left[ \left( \sum_{k=1}^n (\vare_{2k-1} - \mu) \right)^2
           - \sum_{k=1}^n (\vare_{2k-1} - \mu)^2 \right]
  = \frac{1}{2}
    \left[ n ( \cT_1^{(n)} )^2
           - \sum_{k=1}^n (\vare_{2k-1} - \mu)^2 \right] ,
 \end{multline*}
 and similarly, the forth summand is
 \[
   \sum_{k=1}^n (\vare_{2k} - \mu) \sum_{j=1}^{k-1} (\vare_{2j} - \mu)
   = \frac{1}{2}
     \left[ n ( \cS_1^{(n)} )^2 - \sum_{k=1}^n (\vare_{2k} - \mu)^2 \right] .
 \]
Consequently,
 \[
   \frac{1}{2n} \sum_{k=1}^{2n} M_k \bigl(V_{k-1} - \EE(V_{k-1})\bigr)
   = - \frac{1}{4} \left[ (\cS_1^{(n)} - \cT_1^{(n)})^2 - 2 \sigma^2 \right]
     + \frac{1}{4n} \sum_{k=1}^{2n} (\vare_k - \mu)^2
     - \frac{1}{2} \sigma^2 .
 \]
By the strong law of large numbers
 \ $(2n)^{-1} \sum_{k=1}^{2n} (\vare_k - \mu)^2 \as \sigma^2$ \ as
 \ $n \to \infty$, \ hence we obtain \eqref{MV_ST} for the subsequence
 \ $(2n)_{n\in\NN}$.
\ Note also that the convergence in \eqref{MV_ST} holds almost surely, too.

Now we turn to prove \eqref{M_ST} for the subsequence \ $(2n)_{n\in \NN}$.
\ First observe that
 \begin{align*}
  \sum_{k=1}^{2n} M_k
  = \sum_{k=1}^{2n} (\vare_k - \mu)
  = \sum_{k=1}^n (\vare_{2k} - \mu) + \sum_{k=1}^n (\vare_{2k-1} - \mu) .
 \end{align*}
Hence
 \begin{align*}
  \frac{1}{(2n)^{1/2}} \sum_{k=1}^{2n} M_k
  = \frac{1}{2^{1/2}}\left[ \frac{1}{n^{1/2}}(X_{2k} - \EE(X_{2k})) 
                            + \frac{1}{n^{1/2}}(X_{2k-1} - \EE(X_{2k-1}))  \right]
  = \frac{1}{2^{1/2}} \biggl( \cS_1^{(n)} + \cT_1^{(n)} \biggr) ,
 \end{align*}
 thus we obtain \eqref{M_ST} for the subsequence \ $(2n)_{n\in\NN}$.

Finally, one can show \eqref{VV_ST}, \eqref{MX_ST}, \eqref{MV_ST}, and
 \eqref{M_ST} for the
 subsequence \ $(2n - 1)_{n\in\NN}$ \ in the same way.
\proofend

\section{Estimations of moments}
\label{section_moments}

In the proofs of Theorem \ref{main}, Theorem \ref{10main} and Theorem
 \ref{01main} good bounds for moments of the random variables
 \ $(M_k)_{k\in\ZZ_+}$, \ $(X_k)_{k\in\ZZ_+}$, \ $(U_k)_{k\in\ZZ_+}$ \ and
 \ $(V_k)_{k\in\ZZ_+}$ \ are extensively used.
First note that, for all \ $k \in \NN$, \ $\EE( M_k \mid \cF_{k-1} ) = 0$ \ and
 \ $\EE(M_k) = 0$, \ since \ $M_k = X_k - \EE(X_k \mid \cF_{k-1})$.

\begin{Lem}\label{Moments}
Let \ $(X_k)_{k \geq -1}$ \ be an \INARtwo\ process.
Suppose that \ $X_0 = X_{-1} = 0$ \ and \ $\EE(\vare_1^2) < \infty$.
\ Then, for all \ $k, \ell \in \NN$,
 \begin{gather}
  \EE( M_k M_\ell \mid \cF_{\max\{k,\ell\}-1} )
  = \begin{cases}
     \alpha(1-\alpha) X_{k-1} + \beta(1-\beta) X_{k-2} + \sigma^2
      & \text{if \ $k = \ell$,} \\
     0 & \text{if \ $k \ne \ell$,}
    \end{cases} \label{Mcond} \\
  \EE( M_k M_\ell )
  = \begin{cases}
     \alpha(1-\alpha) \EE(X_{k-1}) + \beta(1-\beta) \EE(X_{k-2})
      + \sigma^2 & \text{if \ $k = \ell$,} \\
     0 & \text{if \ $k \ne \ell$,}
    \end{cases} \label{Cov} \\
  \EE( M_k^3 \mid \cF_{k-1})
  = X_{k-1} \EE\big[(\xi_{1,1} - \EE(\xi_{1,1}))^3\big]
    + X_{k-2} \EE\big[(\eta_{1,1} - \EE(\eta_{1,1}))^3\big]
    + \EE\big[(\vare_1 - \EE(\vare_1))^3\big] , \label{M3cond} \\
  \EE( M_k^3) = \EE\big[(\xi_{1,1} - \EE(\xi_{1,1}))^3\big] \EE(X_{k-1})
                + \EE\big[(\eta_{1,1} - \EE(\eta_{1,1}))^3\big] \EE(X_{k-2})
                + \EE\big[(\vare_1 - \EE(\vare_1))^3\big] . \label{M3}
 \end{gather}
\end{Lem}

\noindent
\textbf{Proof.} \
By \eqref{INAR2} and \eqref{Mk},
 \begin{equation}\label{Mdeco}
  M_k = \sum_{j=1}^{X_{k-1}} \big( \xi_{k,j} - \EE(\xi_{k,j}) \big)
        + \sum_{j=1}^{X_{k-2}} \big( \eta_{k,j} - \EE(\eta_{k,j}) \big)
        + \big( \vare_k - \EE(\vare_k) \big) , \qquad k\in\NN.
 \end{equation}
For all \ $k\in\NN$, \ the random variables
 \ $\big\{\xi_{k,j} - \EE(\xi_{k,j}) , \, \eta_{k,j} - \EE(\eta_{k,j}) , \,
          \vare_k - \EE(\vare_k) : j \in \NN \big\}$
 \ are independent of each other, independent of \ $\cF_{k-1}$, \ and have
 zero mean, thus in case \ $k = \ell$ \ we conclude \eqref{Mcond} and hence
 \eqref{Cov}.
If \ $k < \ell$, \ then
 \ $\EE( M_k M_\ell \mid \cF_{\ell-1} ) = M_k \EE( M_\ell \mid \cF_{\ell-1} )
    = 0$.
\ Thus we obtain \eqref{Mcond} and \eqref{Cov} in case \ $k \ne \ell$.
\ Shedding more light we give more details for deriving \eqref{M3cond} and
 \eqref{M3}.
Namely, using multinomial theorem the above mentioned properties of the random
 variables
 \ $\big\{\xi_{k,j} - \EE(\xi_{k,j}) , \, \eta_{k,j} - \EE(\eta_{k,j}) , \,
          \vare_k - \EE(\vare_k) : j \in \NN \big\}$
 \ yield that
 \begin{align*}
   \EE(M_k^3\mid\cF_{k-1})
   & = \EE\left( \sum_{j=1}^{X_{k-1}} (\xi_{k,j} - \EE(\xi_{k,j}))^3
                + \sum_{j=1}^{X_{k-2}} (\eta_{k,j} - \EE(\eta_{k,j}))^3
                + (\vare_k - \EE(\vare_k))^3 \;\Big\vert \;\cF_{k-1}\right)\\
   & = X_{k-1} \EE[(\xi_{1,1} - \EE(\xi_{1,1}))^3]
       + X_{k-2} \EE[(\eta_{1,1} - \EE(\eta_{1,1}))^3]
       + \EE[(\vare_1 - \EE(\vare_1))^3] .
 \end{align*}
This readily implies \eqref{M3}.
\proofend

\begin{Lem}\label{LEM_moments_seged}
Let \ $(\zeta_k)_{k\in\NN}$ \ be independent and identically distributed random
 variables such that \ $\EE\bigl(|\zeta_1|^\ell\bigr) < \infty$ \ for some
 \ $\ell \in \NN$.
\begin{enumerate}
 \item[\textup{(i)}]
  If \ $\EE(\zeta_1) \ne 0$, \ then there exists a polynomial \ $Q_\ell$ \ of
   degree \ $\ell$ \ such that its leading coefficient is
   \ $\bigl[\EE(\zeta_1)\bigr]^\ell$ \ and
   \[
     \EE\bigl((\zeta_1 +\cdots + \zeta_N)^\ell\bigr) = Q_\ell(N) , \qquad
     N \in \NN .
   \]
 \item[\textup{(ii)}]
  If \ $\EE(\zeta_1) = 0$, \ then there exists a polynomial \ $R_\ell$ \ of
   degree at most \ $\ell/2$ \ such that
   \[
     \EE\bigl((\zeta_1 +\cdots + \zeta_N)^\ell\bigr) = R_\ell(N) , \qquad
     N \in \NN .
   \]
\end{enumerate}
The coefficients of the polynomials in question depend on the moments
 \ $\EE(\zeta_1^j)$, \ $j \in \{1, \ldots, \ell\}$.
\end{Lem}

\noindent
\textbf{Proof.}
(i) By multinomial theorem,
 \begin{align*}
  \EE\bigl((\zeta_1 + \cdots + \zeta_N)^\ell\bigr)
  &= \sum_{\underset{\ell_1, \ldots, \ell_N \in \ZZ_+}{\ell_1 + \cdots + \ell_N = \ell,}}
      \frac{\ell!}{\ell_1!\cdots \ell_N!}
      \EE(\zeta_1^{\ell_1} \cdots \zeta_N^{\ell_N})\\
  &= \sum_{\underset{\ell_1, \ldots, \ell_N \in \ZZ_+}{\ell_1 + \cdots + \ell_N = \ell,}}
      \frac{\ell!}{\ell_1! \cdots \ell_N!}
      \EE(\zeta_1^{\ell_1}) \cdots \EE(\zeta_1^{\ell_N}) \\
  &= \sum_{\underset{k_1, \ldots, k_s \in \ZZ_+, \; 1 \leq s \leq \ell}
                    {k_1 + 2 k_2 + \cdots + s k_s = \ell,}}
      \binom{N}{k_1} \binom{N - k_1}{k_2} \cdots
      \binom{N - k_1 - \cdots - k_{s-1}}{k_s} \\
  &\phantom{=\sum_{\underset{k_1, \ldots, k_s \in \ZZ_+, \; 1 \leq s \leq \ell}
                            {k_1 + 2 k_2 + \cdots + s k_s = \ell,}}\;}
     \times \frac{\ell!}{(2!)^{k_2}(3!)^{k_3} \cdots (s!)^{k_s}}
            \bigl[\EE(\zeta_1)\bigr]^{k_1} \cdots
            \bigl[\EE(\zeta_1^s)\bigr]^{k_s} .
 \end{align*}
Since
 \[
   \binom{N}{k_1}\binom{N - k_1}{k_2} \cdots
   \binom{N - k_1 - \cdots - k_{s-1}}{k_s}
   = \frac{N (N-1) \cdots (N - k_1 - k_2 - \cdots - k_s + 1)}
          {k_1! k_2! \cdots k_s!}
 \]
 is a polynomial of the variable \ $N$ \ having degree
 \ $k_1 + \cdots + k_s \leq \ell$, \ there is a polynomial \ $Q_\ell$ \ of
 degree at most \ $\ell$ \ such that
 \ $\EE\bigl((\zeta_1 + \cdots + \zeta_N)^\ell\bigr) = Q_\ell(N)$,
 \ $N \in \NN$.
\ Note that a term of degree \ $\ell$ \ can occur only in the case
 \ $k_1 + \cdots + k_s = \ell$.
\ Since \ $k_1 + 2k_2 + \cdots + sk_s = \ell$, \ we have \ $s = 1$ \ and
 \ $k_1 = \ell$, \ and the corresponding term of degree \ $\ell$ \ is
 \ $N (N-1) \cdots (N-\ell+1) \bigl[\EE(\zeta_1)\bigr]^\ell$.
\ Hence \ $Q_\ell$ \ is polynomial of degree \ $\ell$ \ having leading
 coefficient \ $\bigl[\EE(\zeta_1)\bigr]^\ell$.

(ii) Using again the multinomial theorem we have
 \begin{align*}
  \EE\bigl((\zeta_1 + \cdots + \zeta_N)^\ell\bigr)
  &= \sum_{\underset{\ell_1, \ldots, \ell_N \in \ZZ_+}{\ell_1 + \cdots + \ell_N = \ell,}}
      \frac{\ell!}{\ell_1!\cdots \ell_N!}
      \EE(\zeta_1^{\ell_1} \cdots \zeta_N^{\ell_N}) \\
  &= \sum_{\underset{\ell_1, \ldots, \ell_N \in \ZZ_+ \setminus \{1\}}%
                    {\ell_1 + \cdots + \ell_N = \ell,}}
      \frac{\ell!}{\ell_1! \cdots \ell_N!}
      \EE(\zeta_1^{\ell_1}) \cdots \EE(\zeta_1^{\ell_N}) \\
  &= \sum_{\underset{k_2, \ldots, k_s \in \ZZ_+, \; 2 \leq s \leq \ell}
                    {2 k_2 + 3 k_3 + \cdots + s k_s = \ell,}}
      \binom{N}{k_2} \binom{N - k_2}{k_3} \cdots
      \binom{N - k_2 - \cdots - k_{s-1}}{k_s} \\
  &\phantom{=\sum_{\underset{k_2, \ldots, k_s \in \ZZ_+, \; 2 \leq s \leq \ell}
                            {2 k_2 + 3 k_3 + \cdots + s k_s = \ell,}}\;}
     \times \frac{\ell!}{(2!)^{k_2}(3!)^{k_3} \cdots (s!)^{k_s}}
            \bigl[\EE(\zeta_1^2)\bigr]^{k_2} \cdots
            \bigl[\EE(\zeta_1^s)\bigr]^{k_s} .
 \end{align*}
Here
 \[
   \binom{N}{k_2}\binom{N - k_2}{k_3} \cdots
   \binom{N - k_2 - \cdots - k_{s-1}}{k_s}
   = \frac{N (N-1) \cdots (N - k_2 - k_3 - \cdots - k_s + 1)}
          {k_2! k_3! \cdots k_s!}
 \]
 is a polynomial of the variable \ $N$ \ having degree \ $k_2 + \cdots + k_s$.
\ Since
 \[
     \ell=2k_2 + 3k_3+ \cdots + sk_s \geq 2(k_2+k_3+\cdots+k_s),
 \]
 we have \ $k_2 + \cdots + k_s \leq \ell/2$ \ yielding part (ii).
Note that if \ $\ell$ \ is even and \ $\EE(\zeta_1^2) \ne 0$, \ then the degree
 of \ $R_\ell$ \ is \ $\ell/2$; \ if \ $\ell$ \ is odd and
 \ $\EE(\zeta_1^2) \ne 0$, \ $\EE(\zeta_1^3) \ne 0$, \ then the degree of
 \ $R_\ell$ \ is also \ $\ell/2$.
\proofend

\begin{Rem}
In what follows using the proof of Lemma \ref{LEM_moments_seged} we give a bit
 more explicit form of the polynomial \ $R_\ell$ \ in part (ii) of Lemma
 \ref{LEM_moments_seged} for the special cases
 \ $\ell \in \{1, 2, 3, 4, 5, 6\}$.
\ If \ $\ell = 1$, \ then \ $\EE(\zeta_1 + \cdots + \zeta_N) = 0 $ \ and
 \ $R_1 : \RR \to \RR$, \ $R_1(x) := 0$, \ $x \in \RR$.

\noindent If \ $\ell = 2$, \ then
 \[
   \EE((\zeta_1 + \cdots + \zeta_N)^2) = N \EE(\zeta_1^2) ,
 \]
 and \ $R_2 : \RR \to \RR$, \ $R_2(x) := \EE(\zeta_1^2) x$, \ $x \in \RR$.

\noindent If \ $\ell = 3$, \ then
 \[
   \EE((\zeta_1 + \cdots + \zeta_N)^3) = N \EE(\zeta_1^3) ,
 \]
 and \ $R_3 : \RR \to \RR$, \ $R_3(x) := \EE(\zeta_1^3) x$, \ $x \in \RR$.

\noindent If \ $\ell = 4$, \ then
 \[
   \EE((\zeta_1 + \cdots + \zeta_N)^4)
   = N \EE(\zeta_1^4) + \binom{N}{2} \frac{4!}{2!2!} (\EE(\zeta_1^2))^2 ,
 \]
 and \ $R_4 : \RR \to \RR$,
 \ $R_4(x) := \EE(\zeta_1^4) x + 3 (\EE(\zeta_1^2))^2 x (x-1)$, \ $x \in \RR$.

\noindent If \ $\ell = 5$, \ then
 \[
   \EE((\zeta_1 + \cdots + \zeta_N)^5)
   = N \EE(\zeta_1^5)
     + 2 \binom{N}{2} \frac{5!}{2!3!} \EE(\zeta_1^3) \EE(\zeta_1^2) ,
 \]
 and \ $R_5 : \RR \to \RR$,
 \ $R_5(x) := \EE(\zeta_1^5) x + 10 \EE(\zeta_1^3) \EE(\zeta_1^2) x (x-1)$,
 \ $x \in \RR$.

\noindent If \ $\ell = 6$, \ then
 \[
   \EE((\zeta_1 + \cdots + \zeta_N)^6)
   = N \EE(\zeta_1^6)
     + 2 \binom{N}{2} \frac{6!}{2!4!} \EE(\zeta_1^4) \EE(\zeta_1^2)
     + \binom{N}{2} \frac{6!}{3!3!} (\EE(\zeta_1^3))^2
     + \binom{N}{3} \frac{6!}{2!2!2!} (\EE(\zeta_1^2))^3 ,
 \]
 and \ $R_6 : \RR \to \RR$,
 \[
   R_6(x) := \EE(\zeta_1^6) x + 15 \EE(\zeta_1^4) \EE(\zeta_1^2) x (x-1)
             + 10 (\EE(\zeta_1^3))^2 x (x-1)
             + 15 (\EE(\zeta_1^2))^3 x (x-1) (x-2), \quad x \in \RR .
 \]
\proofend
\end{Rem}

\begin{Lem}\label{LEM_Putzer}
If \ $\alpha + \beta = 1$, \ then the matrix \ $A$ \ defined in \eqref{bA} has
 eigenvalues \ $1$ \ and \ $\alpha - 1 = -\beta$, \ and the powers of \ $A$
 \ take the following form
 \begin{align*}
  A^k = \frac{1}{1+\beta}
        \begin{bmatrix}
         1 & \beta \\
         1 & \beta \\
        \end{bmatrix}
        + \frac{(-\beta)^k}{1+\beta}
          \begin{bmatrix}
           \beta & -\beta \\
           -1 & 1 \\
          \end{bmatrix}
        = \bu \tbu^\top + (-\beta)^k \bv \tbv^\top ,
          \qquad k \in \ZZ_+ ,
 \end{align*}
 with
 \[
   \bu := \frac{1}{1 + \beta}
          \begin{bmatrix} 1 \\ 1 \end{bmatrix}, \quad
   \tbu := \begin{bmatrix} 1 \\ \beta \end{bmatrix}, \quad
   \bv  := \frac{1}{1 + \beta}
           \begin{bmatrix} \beta \\ -1 \end{bmatrix}, \quad
   \tbv := \begin{bmatrix} 1 \\ -1 \end{bmatrix} .
 \]
\end{Lem}

\noindent\textbf{Proof.}
The formula for the powers of \ $A$ \ follows by the so-called Putzer's
 spectral formula, see, e.g., Putzer \cite{Put}.
\proofend

\begin{Rem}\label{REM_Putzer}
Using Lemma \ref{LEM_Putzer} we obtain the decomposition
 \begin{equation}\label{Xdeco}
  \begin{bmatrix} X_k \\ X_{k-1} \end{bmatrix}
  = U_k \bu + V_k \bv
  = \frac{1}{1 + \beta}
    \begin{bmatrix} U_k + \beta V_k \\ U_k - V_k \end{bmatrix}
  = \frac{1}{1+\beta}
    \begin{bmatrix} 1 & \beta \\ 1 & -1 \end{bmatrix}
    \begin{bmatrix} U_k \\ V_k \end{bmatrix} , \qquad k \in \NN ,
 \end{equation}
 with
 \begin{equation}%\label{U_V}
  U_k = X_k + \beta X_{k-1} , \qquad
  V_k = X_k - X_{k-1} ,
  \qquad k \in \NN .
 \end{equation}
Note that \eqref{Xdeco} is valid for \ $k = 0$ \ with the convention
 \ $U_0 := 0$ \ and \ $V_0 := 0$.
\ The decomposition \eqref{Xdeco} can be considered as a motivation for the
 definition of \ $U_k$ \ and \ $V_k$, \ $k \in \NN$, \ given in Sections
 \ref{section_estimators} and \ref{section_proof_main}.
\proofend
\end{Rem}

\begin{Lem}\label{LEM_moments_X}
Let \ $(X_k)_{k\geq-1}$ \ be an \INARtwo\ process with autoregressive parameters
 \ $(\alpha, \beta) \in [0, 1]^2$ \ such that \ $\alpha + \beta = 1$
 \ (hence it is unstable).
Suppose that \ $X_0 = X_{-1} = 0$ \ and \ $\EE(\vare_1^\ell) < \infty$ \ with
 some \ $\ell \in \NN$.
\ Then there exists a constant \ $c_\ell$ \ such that
 \ $\EE(X_n^{\ell_1}X_{n-1}^{\ell_2}) \leq c_\ell n^\ell$, \ $n \in \NN$, \ for all
 \ $\ell_1, \ell_2 \in \ZZ_+$ \ with \ $\ell_1 + \ell_2 \leq \ell$.
\end{Lem}

\noindent
\textbf{First proof.}
Observe that the statement is equivalent with the following: for each
 polynomial $P$ of two variables having degree at most $\ell$, there
 exists a constant $c_P$ such that
 $\EE\bigl(|P(X_n, X_{n-1})|\bigr) \leq c_P n^\ell$, \ $n \in \NN$.

First let us suppose that \ $(\alpha, \beta) \in (0, 1)^2$.
\ If \ $\ell = 1$, \ i.e., \ $(\ell_1, \ell_2) = (1, 0)$ \ or
 \ $(\ell_1, \ell_2) = (0, 1)$, \ then to conclude the statement we show that
 \begin{align}\label{seged_LEM_moments_X}
  \EE(X_n) = \frac{\mu}{1+\beta} n
             + \frac{\mu \beta}{(1+\beta)^2} (1 - (-\beta)^n) ,
  \qquad n \in \NN .
 \end{align}
Since \ $\EE(X_n \mid \cF_{n-1}) = \alpha X_{n-1} + \beta X_{n-2} + \mu$,
 \ $n \in \NN$, \ we have
 \ $\EE(X_n) = \alpha \EE(X_{n-1}) + \beta \EE(X_{n-2}) + \mu$,
 \ $n \in \NN$, \ yielding that
 \begin{align*}
  \begin{bmatrix}
    \EE(X_n) \\
    \EE(X_{n-1}) \\
  \end{bmatrix}
   = \begin{bmatrix}
       \alpha & \beta \\
       1 & 0 \\
     \end{bmatrix}
     \begin{bmatrix}
       \EE(X_{n-1}) \\
       \EE(X_{n-2}) \\
     \end{bmatrix}
     +
     \begin{bmatrix}
       \mu \\
       0 \\
     \end{bmatrix}
   = A\begin{bmatrix}
       \EE(X_{n-1}) \\
       \EE(X_{n-2}) \\
      \end{bmatrix}
      +
      \begin{bmatrix}
        \mu \\
        0 \\
      \end{bmatrix}, \qquad n \in \NN .
 \end{align*}
By Lemma \ref{LEM_Putzer},  we get
 \begin{align*}
  \begin{bmatrix}
    \EE(X_n) \\
    \EE(X_{n-1}) \\
  \end{bmatrix}
   = \sum_{j=1}^n A^{n-j}
       \begin{bmatrix}
        \mu \\
        0 \\
       \end{bmatrix}
  = \left(\frac{n}{1+\beta}
          \begin{bmatrix}
            1 & \beta \\
            1 & \beta \\
          \end{bmatrix}
          + \frac{1-(-\beta)^n}{(1+\beta)^2}
            \begin{bmatrix}
            \beta & -\beta \\
            -1 & 1 \\
          \end{bmatrix}
    \right)
    \begin{bmatrix}
        \mu \\
          0 \\
        \end{bmatrix},\quad n\in\NN,
 \end{align*}
 which yields \eqref{seged_LEM_moments_X}.

\noindent
Let us suppose now that the statement holds for \ $1, \ldots, \ell-1$.
\ By multinomial theorem,
 \begin{align}\label{Xmoment_seged_uj}
  X_n^k = \sum_{\underset{k_1, k_2,k_3 \in \ZZ_+}
                         {k_1 + k_2 + k_3 = k,}}
           \frac{k!}{k_1!k_2!k_3!}
           \left( \sum_{j=1}^{X_{n-1}} \xi_{n,j}\right)^{k_1}
           \left( \sum_{j=1}^{X_{n-2}} \eta_{n,j}\right)^{k_2}
           \vare_n^{k_3} , \qquad k \in \NN .
 \end{align}
Since for all \ $n \in \NN$ \ the random variables
 \ $\{\xi_{n,j}, \eta_{n,j}, \vare_n : j \in \NN\}$ \ are independent of each
 other and of the $\sigma$-algebra \ $\cF_{n-1}$, \ we have for all
 \ $\ell_1, \ell_2 \in \ZZ_+$ \ with \ $\ell_1 + \ell_2 = \ell$
 \begin{align*}
  &\EE(X_n^{\ell_1} X_{n-1}^{\ell_2} \mid \cF_{n-1}) \\
  &\quad
   = X_{n-1}^{\ell_2} \sum_{\underset{k_1,k_2,k_3 \in \ZZ_+}
                                  {k_1 + k_2 + k_3 = \ell_1,}}
     \frac{\ell_1!}{k_1!k_2!k_3!}
     \EE\left(\left(\sum_{j=1}^M \xi_{n,j}\right)^{k_1}\right)\Bigg|_{M=X_{n-1}}
     \EE\left(\left(\sum_{j=1}^N \eta_{n,j}\right)^{k_2}\right)\Bigg|_{N=X_{n-2}}
     \EE(\vare_1^{k_3}) .
 \end{align*}
Using part (i) of Lemma \ref{LEM_moments_seged} and separating the terms
 having degree \ $\ell$ \ and less than \ $\ell$, \ we have
 \begin{align*}
  \EE(X_n^{\ell_1}X_{n-1}^{\ell_2}\mid \cF_{n-1})
  = \sum_{\underset{k_1,k_2\in \ZZ_+}
                   {k_1 + k_2 = \ell_1,}}
     \frac{\ell_1!}{k_1! k_2!}
     \alpha^{k_1} X_{n-1}^{\ell_2+k_1} \beta^{k_2} X_{n-2}^{k_2}
     + Q_{\ell_1,\ell_2}(X_{n-1}, X_{n-2}) ,
 \end{align*}
 where \ $Q_{\ell_1,\ell_2}$ \ is a polynomial of two variables having degree at
 most $\ell - 1$.
\ Hence
 \[
   \EE(X_n^{\ell_1} X_{n-1}^{\ell_2})
   = \sum_{\underset{k_1,k_2\in \ZZ_+}
                    {k_1 + k_2 = \ell_1,}}
      \frac{\ell_1!}{k_1!k_2!}
      \alpha^{k_1} \beta^{k_2} \EE\big(X_{n-1}^{\ell_2+k_1} X_{n-2}^{k_2}\big)
     + \EE\big(Q_{\ell_1,\ell_2}(X_{n-1}, X_{n-2})\big) .
 \]
By the induction hypothesis (used for polynomials, see the beginning of the
 proof), there exists a constant \ $c_{Q_{\ell_1, \ell_2}}$ \ such that
 \ $\EE\big(|Q_{\ell_1,\ell_2}(X_n,X_{n-1})|\big) \leq c_{Q_{\ell_1,\ell_2}} n^{\ell-1}$,
 \ $n \in \NN$.
\ In fact, we have
 \begin{align}\label{univ_const}
  \EE\big(|Q_{\ell_1,\ell_2}(X_n, X_{n-1})|\big) \leq \tc_\ell n^{\ell-1}
 \end{align}
 for \ $n \in \NN$ \ and \ $\ell_1, \ell_2 \in \ZZ_+$ \ with
 \ $\ell_1 + \ell_2 = \ell$, \ where
 \ $\tc_\ell := \max_{0 \leq i \leq \ell} c_{Q_{i,\ell-i}}$.
\ Consequently, we have
 \[
   \EE(X_n^{\ell_1} X_{n-1}^{\ell_2})
   \leq \sum_{\underset{k_1,k_2\in \ZZ_+}
                       {k_1 + k_2 = \ell_1,}}
         \frac{\ell_1!}{k_1! k_2!}
         \alpha^{k_1} \beta^{k_2} \EE\big(X_{n-1}^{\ell_2+k_1} X_{n-2}^{k_2}\big)
         + \tc_\ell (n-1)^{\ell-1} .
 \]
Similarly, for all \ $k_1, k_2 \in \ZZ_+$ \ with \ $k_1 + k_2 = \ell_1$, \ we
 have
 \begin{align*}
  \EE\big(X_{n-1}^{\ell_2+k_1} X_{n-2}^{k_2}\big)
  = \sum_{\underset{j_1,j_2\in \ZZ_+}
                   {j_1 + j_2 = \ell_2+k_1,}}
     \frac{(\ell_2+k_1)!}{j_1! j_2!}
     \alpha^{j_1} \beta^{j_2} \EE\big(X_{n-2}^{k_2+j_1} X_{n-3}^{j_2}\big)
    + \EE\big(Q_{\ell_2+k_1,k_2}(X_{n-2}, X_{n-3})\big) .
 \end{align*}
Hence we have
 \begin{align*}
  \EE(X_n^{\ell_1}X_{n-1}^{\ell_2})
  &= \sum_{\underset{k_1, k_2 \in \ZZ_+}
                    {k_1 + k_2 = \ell_1,}}
      \frac{\ell_1!}{k_1! k_2!}  \alpha^{k_1} \beta^{k_2}
      \sum_{\underset{j_1, j_2 \in \ZZ_+}
                     {j_1 + j_2 = \ell_2 + k_1,}}
       \frac{(\ell_2+k_1)!}{j_1!j_2!} \alpha^{j_1} \beta^{j_2}
       \EE\big(X_{n-2}^{k_2+j_1} X_{n-3}^{j_2}\big)\\
  &\quad
     + \sum_{\underset{k_1, k_2 \in \ZZ_+}
                      {k_1 + k_2 = \ell_1,}}
        \frac{\ell_1!}{k_1!k_2!}
        \alpha^{k_1} \beta^{k_2} \EE\big(Q_{\ell_2+k_1,k_2}(X_{n-2},X_{n-3})\big)
       + \EE\big(Q_{\ell_1, \ell_2}(X_{n-1},X_{n-2})\big) .
 \end{align*}
Applying \eqref{univ_const} and
 \begin{align*}
  \sum_{\underset{k_1,k_2\in \ZZ_+}
                 {k_1 + k_2 = \ell_1,}}
   \frac{\ell_1!}{k_1!k_2!}  \alpha^{k_1} \beta^{k_2}
  = (\alpha + \beta)^{\ell_1} = 1 ,
 \end{align*}
 we conclude
 \begin{align*}
  \EE(X_n^{\ell_1}X_{n-1}^{\ell_2})
  &\leq \sum_{\underset{k_1, k_2 \in \ZZ_+}
                       {k_1 + k_2 = \ell_1,}}
         \frac{\ell_1!}{k_1!k_2!}  \alpha^{k_1} \beta^{k_2}
         \sum_{\underset{j_1, j_2 \in \ZZ_+}
                        {j_1 + j_2 = \ell_2+k_1,}}
          \frac{(\ell_2+k_1)!}{j_1!j_2!} \alpha^{j_1} \beta^{j_2}
          \EE\big(X_{n-2}^{k_2+j_1} X_{n-3}^{j_2}\big) \\
  &\quad + \tc_\ell (n-2)^{\ell-1} + \tc_\ell (n-1)^{\ell-1} .
 \end{align*}
Using that \ $\EE(X_1^r X_0^q) = 0$, \ $r, q \in \ZZ_+$ \ (since \ $X_0 = 0$),
 after \ $n-1$ \ steps, one can derive
 \begin{align*}
  \EE(X_n^{\ell_1} X_{n-1}^{\ell_2})
  \leq \tc_\ell \sum_{i=1}^{n-1} i^{\ell-1}
  \leq \tc_\ell n \cdot n^{\ell-1} = \OO(n^\ell) , \qquad n \in \NN ,
 \end{align*}
 that is, \ $\EE(P(X_n, X_{n-1})) \leq \tc_\ell n^\ell$ \ for all monomials
 \ $P(x, y) := x^{\ell_1} y^{\ell_2}$, \ $x, y \in \RR$, \ with
 \ $\ell_1 + \ell_2 = \ell$, \ $\ell_1, \ell_2 \in \ZZ_+$.
\ If \ $P$ \ has the form
 \[
   P(x, y) := \sum_{i=0}^\ell p_i x^i y^{\ell-i} + Q(x, y), \qquad x, y \in \RR ,
 \]
 where \ $p_i \in \RR$, \ $i \in \{0, \ldots, \ell\}$, \ and \ $Q$ \ is a
 polynomial of two variables having degree at most \ $\ell - 1$, \ then for
 all \ $n \in \NN$,
 \begin{align*}
  \EE(|P(X_n, X_{n-1})|)
  \leq \sum_{i=0}^\ell |p_i| \EE(X_n^i X_{n-1}^{\ell-i}) + \EE(Q(X_n, X_{n-1}))
  \leq \left( \sum_{i=0}^\ell |p_i| c_\ell \right)n^\ell + c_Q n^{\ell-1}
  \leq c_P n^\ell ,
 \end{align*}
 where \ $c_P := c_Q + c_\ell \sum_{i=0}^\ell |p_i|$, \ as desired.

Next let us suppose that \ $(\alpha, \beta) = (1, 0)$.
\ Then \ $X_n = X_{n-1} + \vare_n$, \ $n \in \NN$, \ which implies that
 \ $X_n = \sum_{i=1}^n \vare_i$, \ $n \in \NN$.
\ By part (i) of Lemma \ref{LEM_moments_seged},
 \begin{align}\label{seged3_LEM_moments_X}
  \EE(X_n^\ell) = Q_\ell(n), \qquad n \in \NN ,
 \end{align}
 where \ $Q_\ell$ \ is a polynomial of degree \ $\ell$.
\ If \ $\ell_1, \ell_2 \in \ZZ_+$ \ with $\ell_1 + \ell_2 \leq \ell$, \ then
 using independence of \ $X_{n-1}$ \ and \ $\vare_n$ \ we have
 \begin{align*}
  \EE(X_n^{\ell_1}X_{n-1}^{\ell_2})
  &= \EE((X_{n-1} + \vare_n)^{\ell_1} X_{n-1}^{\ell_2})
   = \EE\left( \sum_{j=0}^{\ell_1}
                \binom{\ell_1}{j} X_{n-1}^j \vare_n^{\ell_1-j}
                X_{n-1}^{\ell_2} \right) \\
  &= \sum_{j=0}^{\ell_1}
      \binom{\ell_1}{j} \EE(X_{n-1}^{j+\ell_2}) \EE(\vare_n^{\ell_1-j}) , \qquad
  n \in \NN .
 \end{align*}
Using \eqref{seged3_LEM_moments_X},
 \[
   \EE(X_n^{\ell_1} X_{n-1}^{\ell_2})
   = \sum_{j=0}^{\ell_1} \binom{\ell_1}{j} Q_{j+\ell_2}(n-1) \EE(\vare_1^{\ell_1-j})
   = \OO(n^\ell) , \qquad n \in \NN ,
 \]
 since for each \ $j \in \{0, \ldots, \ell_1\}$, \ the polynomial
 \ $Q_{j+\ell_2}$ \ is of degree \ $j + \ell_2 \leq \ell$, \ which yields the
 statement in case \ $(\alpha, \beta) = (1, 0)$.

Finally, let us suppose that \ $(\alpha, \beta) = (0, 1)$.
\ Then \ $X_n = X_{n-2} + \vare_n$, \ $n \in \NN$, \ which implies that
 \[
    X_{2n} = \sum_{i=1}^n \vare_{2i}, \qquad \qquad
    X_{2n-1} = \sum_{i=1}^n \vare_{2i-1}, \qquad n \in \NN .
 \]
By part (i) of Lemma \ref{LEM_moments_seged}, we have
 \begin{align*}%\label{seged3_LEM_moments_X}
  \EE(X_{2n}^\ell) = Q_\ell(n), \qquad n \in \NN ,
  \qquad \qquad
  \EE(X_{2n-1}^\ell) = Q_\ell(n), \qquad n \in \NN ,
 \end{align*}
 where \ $Q_\ell$ \ is a polynomial of degree \ $\ell$.
\ Using the independence of \ $X_{2n}$ \ and \ $X_{2n-1}$, \ for
 \ $\ell_1 + \ell_2 \leq \ell$, \ $\ell_1, \ell_2 \in \ZZ_+$, \ we have
 \[
   \EE(X_{2n}^{\ell_1} X_{2n-1}^{\ell_2})
   = \EE(X_{2n}^{\ell_1}) \EE(X_{2n-1}^{\ell_2})
   = Q_{\ell_1}(n) Q_{\ell_2}(n)
   = \OO(n^\ell) , \qquad n \in \NN ,
 \]
 as desired.
Similarly,
 \[
   \EE(X_{2n-1}^{\ell_1} X_{2n-2}^{\ell_2})
   = \EE(X_{2n-1}^{\ell_1}) \EE(X_{2n-2}^{\ell_2})
   = Q_{\ell_1}(n) Q_{\ell_2}(n-1) = \OO(n^\ell) , \qquad n \in \NN .
 \]
 Hence we have the assertion.

\medskip

\noindent
\textbf{Second proof.}
It is enough to prove that there exists some \ $c_\ell \in \RR_+$ \ such that
 \ $\EE(X_n^{\ell_1} X_{n-1}^{\ell_2}) \leq c_\ell n^\ell$ \ for all \ $n \in \NN$
 \ and \ $\ell_1, \ell_2 \in \ZZ_+$ \ with \ $\ell_1 + \ell_2 = \ell$.
\ Let us introduce the notation
 \[
   \bX_n^{(k)}
   := \begin{bmatrix}
       X_n^k & X_n^{k-1}X_{n-1} & X_n^{k-2} X_{n-1}^2 & \cdots & X_n X_{n-1}^{k-1}
        & X_{n-1}^k
      \end{bmatrix}^\top \in \RR_+^{k+1} ,
   \qquad n, k \in \NN .
 \]
First we check that
 \begin{align}\label{Xmoment_recursion}
  \EE(\bX_n^{(k)} \mid \cF_{n-1})
  = A_k \bX_{n-1}^{(k)} + \sum_{j=1}^{k-1} B_{k,j} \bX_{n-1}^{(j)} + \bmu_k , \qquad
  n \in \NN , \quad k \in \{1, \ldots, \ell\} ,
 \end{align}
 where
 \begin{align*}
  A_k
  := \begin{bmatrix}
      \alpha^k & \binom{k}{1} \alpha^{k-1}\beta & \cdots
       & \binom{k}{k-1}\alpha\beta^{k-1} & \beta^k \\
      \alpha^{k-1} & \binom{k-1}{1}\alpha^{k-2}\beta  & \cdots & \beta^{k-1}
       & 0 \\
      \vdots & \vdots & \ddots & \vdots & \vdots \\
      \alpha & \beta & \cdots & 0 & 0 \\
      1 & 0 & \cdots & 0 & 0
     \end{bmatrix}
     \in \RR_+^{(k+1)\times(k+1)}
 \end{align*}
 and \ $B_{k,j} \in \RR_+^{(k+1)\times(j+1)}$ \ are appropriate matrices of which
 the entries are non-negative and depend only on \ $\alpha$ \ and the moments
 of \ $\vare_1$ \ of order less than or equal to \ $(k-j)$ \ and
 \[
   \bmu_k := \begin{bmatrix}
              \EE(\varepsilon_1^k) & 0 & \cdots & 0 \\
             \end{bmatrix}^\top \in \RR_+^{k+1} .
 \]
For a better understanding, first we give a proof for
 \eqref{Xmoment_recursion} in the case of \ $k = 1$ \ and \ $k = 2$.
\ If \ $k = 1$, \ then
 \begin{align*}
  \EE(\bX_n^{(1)} \mid \cF_{n-1})
  = \begin{bmatrix}
     \EE(X_n \mid \cF_{n-1}) \\
     \EE(X_{n-1} \mid \cF_{n-1} )
    \end{bmatrix}
  = \begin{bmatrix}
     \alpha & \beta \\
     1 & 0
    \end{bmatrix}
    \begin{bmatrix}
     X_{n-1} \\
     X_{n-2}
    \end{bmatrix}
    + \begin{bmatrix}
       \mu \\
       0
      \end{bmatrix}
  = A_1 \bX_{n-1}^{(1)} + \bmu_1,
  \qquad n \in \NN .
 \end{align*}
If \ $k = 2$, \ then, by \eqref{Xmoment_seged_uj}, we have
 \begin{align*}
  \EE(X_n^2\mid\cF_{n-1})
  & = \EE\left( \sum_{\underset{k_1, k_2,k_3 \in \ZZ_+}
                     {k_1 + k_2 + k_3 = 2,}}
                 \frac{2!}{k_1!k_2!k_3!}
                 \left( \sum_{j=1}^{X_{n-1}} \xi_{n,j} \right)^{k_1}
                 \left( \sum_{j=1}^{X_{n-2}} \eta_{n,j} \right)^{k_2}
                 \vare_n^{k_3}
                \; \Bigg| \; \cF_{n-1} \right) \\
  & = \alpha X_{n-1} + \alpha^2 (X_{n-1}^2 - X_{n-1}) + \beta X_{n-2}
      + \beta^2 (X_{n-2}^2 - X_{n-2})
      + 2 \alpha \beta X_{n-1} X_{n-2} \\
  &\phantom{=\;}
      + 2 \alpha X_{n-1} \EE(\vare_1)
      + 2\beta X_{n-2} \EE(\varepsilon_1)
      + \EE(\vare_1^2) , \qquad n \in \NN ,
 \end{align*}
 and hence, using also that
 \[
   \EE(X_n X_{n-1} \mid \cF_{n-1})
   = X_{n-1}\EE(X_n \mid \cF_{n-1})
   = \alpha X_{n-1}^2 + \beta X_{n-1} X_{n-2}+ X_{n-1} \EE(\vare_1) , \qquad
   n \in \NN ,
 \]
 we have
 \begin{align*}
  \EE(\bX_n^{(2)} \mid \cF_{n-1})
  = \EE\left( \begin{bmatrix}
               X_n^2 \\
               X_nX_{n-1} \\
               X_{n-1}^2
              \end{bmatrix}
              \; \Bigg| \; \cF_{n-1} \right)
  = A_2 \begin{bmatrix}
         X_{n-1}^2 \\
         X_{n-1}X_{n-2} \\
         X_{n-2}^2
        \end{bmatrix}
    + B_{2,1} \begin{bmatrix}
              X_{n-1} \\
              X_{n-2}
             \end{bmatrix}
    + \bmu_2,
  \qquad n \in \NN ,
 \end{align*}
 where
 \begin{align*}
  A_2 = \begin{bmatrix}
         \alpha^2 & 2 \alpha \beta & \beta^2 \\
         \alpha &  \beta & 0 \\
         1 & 0 & 0
        \end{bmatrix}
  \qquad \text{and} \qquad
  B_{2,1}
  = \begin{bmatrix}
     \alpha \beta + 2\alpha \EE(\vare_1)
      & \alpha \beta + 2 \beta \EE(\vare_1) \\
     \EE(\vare_1) & 0 \\
     0 & 0
    \end{bmatrix}
 \end{align*}
 as desired.
In the general case using part (i) of Lemma \ref{LEM_moments_seged} one can
 prove \eqref{Xmoment_recursion}
 (giving also explicit forms for the matrices \ $B_{k,j}$).

Taking expectation of \eqref{Xmoment_recursion}, we have
 \begin{align}\label{Xmoment_seged_uj2}
  \EE(\bX_n^{(k)})
  = A_k \EE(\bX_{n-1}^{(k)}) + \sum_{j=1}^{k-1} B_{k,j} \EE(\bX_{n-1}^{(j)})
    + \bmu_k , \qquad n \in \NN , \quad k \in \{1, \ldots, \ell\} .
 \end{align}
For a $d$-dimensional vector \ ${\mathbf v} = (v_i)_{i=1}^d \in \RR^d$ \ and a
 $d \times d$ matrix \ $M = (m_{i,j})_{i,j=1}^d \in \RR^{d\times d}$, \ let us
 introduce the notations
 \begin{align*}
  \|{\mathbf v}\|_\infty := \max_{1\leq i\leq d} |v_i|
  \qquad \text{and} \qquad
  \|M\|_\infty  := \max_{1\leq i\leq d} \sum_{j=1}^d |m_{i,j}| .
 \end{align*}
By the binomial theorem one can easily have \ $\|A_k\|_\infty = 1$,
 \ $k \in \{1, \ldots, \ell\}$.
\ We prove the statement using a double induction with respect to
 \ $k \in \{1, \ldots, \ell\}$ \ and \ $n \in \NN$.
\ First we show that the statement holds for \ $k = 1$ \ using induction with
 respect to \ $n$.
\ Namely, we show that
 \[
   \|\EE(\bX_n^{(1)})\|_\infty  \leq c_1 n , \qquad n \in \NN ,
 \]
 where \ $c_1 := \|\bmu_1\|_\infty$.
\ If \ $n = 1$, \ then
 \begin{align*}
  \EE(\bX_1^{(1)}) = \begin{bmatrix}
                     \EE(X_1) \\
                     \EE(X_0)
                    \end{bmatrix}
                  = \begin{bmatrix}
                     \EE(\vare_1) \\
                     0
                    \end{bmatrix}
                  = \bmu_1 ,
 \end{align*}
 which implies that \ $\|\EE(\bX_1^{(1)})\|_\infty = c_1$.
\ Let us suppose now that \ $\|\EE(\bX_m^{(1)})\|_\infty  \leq c_1 m$ \ holds for
 \ $m \in \{1, \ldots, n-1\}$ \ with \ $n \geq 2$.
\ Then, \eqref{Xmoment_seged_uj2},
 \begin{align*}
  \|\EE(\bX_n^{(1)})\|_\infty
  &= \| A \EE(\bX_{n-1}^{(1)}) + \bmu_1 \|_\infty
     \leq \| A \EE(\bX_{n-1}^{(1)}) \|_\infty + \|\bmu_1\|_\infty \\
  &\leq \|A\|_\infty \|\EE(\bX_{n-1}^{(1)})\|_\infty + \|\bmu_1\|_\infty
   \leq c_1 (n-1) + c_1 = c_1 n ,
 \end{align*}
 as desired.

Let us suppose now that the statement holds for \ $j = 1, \ldots, \ell-1$,
 \ i.e.,
 \[
   \|\EE(\bX_n^{(j)})\|_\infty \leq c_j n^j , \qquad n \in \NN , \quad
   j \in \{1, \ldots, \ell-1\} .
 \]
Next, using induction with respect to \ $n \in \NN$ \ we prove that
 \[
   \|\EE(\bX_n^{(\ell)})\|_\infty \leq c_\ell n^\ell , \qquad n \in \NN ,
 \]
 where
 \[
   c_\ell := \sum_{j=1}^{\ell-1} c_j \|B_{\ell,j}\|_\infty + \|\bmu_\ell\|_\infty .
 \]
If \ $n = 1$, \ then, using that \ $X_0 = 0$ \ and \ $X_1 = \vare_1$, \ we have
 \[
   \EE(\bX_1^{(\ell)})
   = \begin{bmatrix}
      \EE(\vare_1^\ell) & 0 & \cdots & 0
     \end{bmatrix}^\top
   = \bmu_\ell ,
 \]
 which yields that
 \ $ \|\EE(\bX_1^{(\ell)})\|_\infty =  \|\bmu_\ell\|_\infty \leq c_\ell$.
\ Let us suppose now that
 \[
   \|\EE(\bX_m^{(\ell)})\|_\infty \leq c_\ell m^\ell , \qquad
   m \in \{1, \ldots, n-1\} ,
 \]
 where \ $n \geq 2$.
\ Then, by \eqref{Xmoment_seged_uj2},
 \begin{align*}
  \|\EE(\bX_n^{(\ell)})\|_\infty
  &\leq \|A_\ell\|_\infty \|\EE(\bX_{n-1}^{(\ell)})\|_\infty
        + \sum_{j=1}^{\ell-1} \|B_{\ell,j}\|_\infty \|\EE(\bX_{n-1}^{(j)})\|_\infty
        + \|\bmu_\ell\|_\infty \\
  &\leq c_\ell (n-1)^\ell + \sum_{j=1}^{\ell-1} \|B_{\ell,j}\|_\infty c_j (n-1)^j
        + \|\bmu_\ell\|_\infty \\
  &\leq c_\ell (n-1)^\ell
        + \left( \sum_{j=1}^{\ell-1} c_j \|B_{\ell,j}\|_\infty
                 + \|\bmu_\ell\|_\infty \right)
          (n-1)^{\ell-1} \\
  &= c_\ell (n-1)^{\ell-1} (n-1+1) \\
  &\leq c_\ell n^\ell ,
 \end{align*}
 as desired.
\proofend

\begin{Cor}\label{EEX_EEU_EEV}
Let \ $(X_k)_{k \geq -1}$ \ be an \INARtwo\ process with autoregressive
 parameters \ $(\alpha, \beta) \in [0,1]^2$ \ such that \ $\alpha + \beta = 1$
 \ (hence it is unstable).
Suppose that \ $X_0 = X_{-1} = 0$ \ and \ $\EE(\vare_1^\ell) < \infty$ \ with
 some \ $\ell \in \NN$.
\ Then
 \begin{gather*}
   \EE(X_k^i) = \OO(k^i),\qquad
   \EE(M_k^i) = \OO(k^{\lfloor i/2 \rfloor}), \qquad
   \EE(U^i_k ) = \OO(k^i), \qquad
   \EE(V^{2j}_k ) = \OO(k^j),\qquad k\in\NN,
 \end{gather*}
 for \ $i,j\in\ZZ_+$ \ with \ $i\leq \ell$ \ and \ $2j\leq \ell$.
\end{Cor}

\noindent
\textbf{Proof.}
The estimate \ $\EE(X_k^i) = \OO(k^i)$ \ readily follows by Lemma
 \ref{LEM_moments_X}.
Next we turn to prove \ $\EE(M_k^i) = \OO(k^{\lfloor i/2 \rfloor})$.
\ Using \eqref{Mdeco} and that the random variables
 \ $\{\xi_{n,j}, \eta_{n,j}, \vare_n : j \in \NN\}$ \ are independent of each
 other and of the $\sigma$-algebra \ $\cF_{n-1}$, \ we have for all
 \ $n \in \NN$,
 \begin{align*}
  \EE(M_n^i\mid \cF_{n-1})
  &= \sum_{\underset{i_1, i_2, i_3 \in \ZZ_+}
                    {i_1 + i_2 + i_3 = i,}}
      \frac{i!}{i_1!i_2!i_3!}
      \EE\left(\left( \sum_{j=1}^M (\xi_{n,j}-\EE(\xi_{n,j})) \right)^{i_1}
               \right)\Bigg|_{M=X_{n-1}} \\
  &\phantom{= \sum_{\underset{i_1, i_2, i_3 \in \ZZ_+}
                             {i_1 + i_2 + i_3 = i,}}}
      \times \EE\left(\left( \sum_{j=1}^N (\eta_{n,j}-\EE(\eta_{n,j}))
                      \right)^{i_2}\right)\Bigg|_{N=X_{n-2}}
             \EE\left((\vare_n - \EE(\vare_n))^{i_3}\right) .
 \end{align*}
By part (ii) of Lemma \ref{LEM_moments_seged}, there exist polynomials
 \ $Q_{i_1}$, \ $i_1 \in \NN$, \ of degree at most \ $i_1/2$, \ and
 \ $\tQ_{i_2}$, \ $i_2 \in \NN$, \ of degree at most \ $i_2/2$ \ such
 that
 \begin{align*}
  \EE(M_n^i \mid \cF_{n-1})
  = \sum_{\underset{i_1, i_2, i_3 \in \ZZ_+}
                   {i_1 + i_2 + i_3 = i,}}
     \frac{i!}{i_1!i_2!i_3!}
     Q_{i_1}(X_{n-1}) \widetilde Q_{i_2}(X_{n-2})
     \EE\left((\vare_1 - \EE(\vare_1))^{i_3}\right) .
 \end{align*}
Hence
 \begin{align*}
  \EE(M_n^i)
  = \sum_{\underset{i_1, i_2, i_3 \in \ZZ_+}
                   {i_1 + i_2 + i_3 = i,}}
     \frac{i!}{i_1!i_2!i_3!}
     \EE\big(Q_{i_1}(X_{n-1}) \widetilde Q_{i_2}(X_{n-2})\big)
     \EE\left((\vare_1 - \EE(\vare_1))^{i_3}\right) , \qquad n \in \NN .
 \end{align*}
Clearly,
 \ $Q_{i_1}(X_{k-1}) \widetilde Q_{i_2}(X_{k-2}) = Q^*_{i_1+i_2}(X_{k-1},X_{k-2})$,
 \ where \ $Q^*_{i_1+i_2}$ \ is a polynomial of two variables having degree
 at most \ $(i_1 + i_2)/2 \leq i/2$, \ and hence, at most
 \ $\lfloor i/2 \rfloor$. 
\ By Lemma \ref{LEM_moments_X}, there exists a constant \ $c_{Q^*_{i_1+i_2}}$
 \ such that $\EE\big(|Q^*_{i_1+i_2}(X_{k-1},X_{k-2})\big|)
  \leq c_{Q^*_{i_1+i_2}} (k-1)^{\lfloor i/2 \rfloor}$.
\ Hence
 \begin{align*}
  |\EE(M_k^i)|
  \leq (k-1)^{\lfloor i/2 \rfloor}
       \sum_{\underset{i_1, i_2,i_3 \in \ZZ_+}
            {i_1 + i_2 + i_3 = i,}}
        \frac{i!}{i_1!i_2!i_3!}
        c_{Q^*_{i_1+i_2}}
        \left\vert \EE\left((\vare_1 - \EE(\vare_1))^{i_3}\right)\right\vert
 \end{align*}
 for all \ $k \in \NN$, \ as desired.

Next we turn to prove \ $\EE(U^i_k) = \OO(k^i)$, \ $i, k \in \NN$ \ with
 \ $i \leq \ell$.
\ First note that, by power mean inequality, for all \ $i \in \NN$,
 \begin{align*}
  \frac{a+b}{2} \leq \left(\frac{a^i+b^i}{2}\right)^{\frac{1}{i}} , \qquad
  a, b \geq 0 ,
 \end{align*}
 yielding that \ $(a + b)^i \leq 2^{i-1} (a^i + b^i)$, \ $a, b \geq 0$.
\ Hence, by Lemma \ref{LEM_moments_X},
 \begin{align*}
  \EE(U_k^i) = \EE((X_k + \beta X_{k-1})^i)
             \leq 2^{i-1} (\EE(X_k^i) + \beta^i \EE(X_{k-1}^i))
             \leq 2^{i-1} (P_i(k) + \beta^i P_i(k-1)) ,
 \end{align*}
 where \ $P_i$ \ is a polynomial of degree at most \ $i$, \ which yields that
 \ $\EE(U^i_k) = \OO(k^i)$.

Finally, for \ $2 j \leq \ell$, \ $j \in \ZZ_+$, \ we prove
 \ $\EE(V^{2j}_k) = \OO(k^{j})$, \ $k \in \NN$, \ using induction in \ $k$.
\ By the recursion \ $V_k = - \beta V_{k-1} + M_k + \mu$, \ $k \in \NN$,
 \ we have \ $\EE(V_k) = - \beta \EE(V_{k-1}) + \mu$, \ $k \in \NN$,
 \ with initial value \ $\EE(V_0) = 0$, \ hence
 \[
   \EE(V_k) = \mu \sum_{i=0}^{k-1} (-\beta)^i , \qquad k \in \NN ,
 \]
 which yields that \ $|\EE(V_k)| = \OO(1)$.
\ Indeed, for all \ $k \in \NN$,
 \[
   \left|\sum_{i=0}^{k-1} (-\beta)^i \right|
   \leq \begin{cases}
         \frac{1}{1-\beta} & \text{if \ $0 \leq \beta < 1$,} \\
         1 & \text{if \ $\beta=1$,}
        \end{cases}
 \]
 where the inequality for the case \ $\beta = 1$ \ follows by that the
 sequence of partial sums in question is nothing else but the alternating one
 \ $1, 0, 1, 0, 1, 0, \ldots$.
\ Let us introduce the notation \ $\tV_k := V_k - \EE(V_k)$, \ $k \in \NN$.
\ Since, by the triangular inequality for the $L_{2j}$-norm,
 \[
   \left( \EE(V_k^{2j}) \right)^{\frac{1}{2j}}
   \leq \left( \EE(\tV_k^{2j}) \right)^{\frac{1}{2j}} + \EE(|V_k|) ,
 \]
 and \ $|\EE(V_k)| = \OO(1)$, \ for proving \ $\EE(V^{2j}_k) = \OO(k^{j})$,
 \ $k \in \NN$, \ it is enough to show that \ $\EE(\tV^{2j}_k) = \OO(k^{j})$,
 \ $k \in \NN$.
\ Using again the recursion \ $V_k = - \beta V_{k-1} + M_k + \mu$,
 \ $k \in \NN$, \ we get \ $\tV_k = - \beta \tV_{k-1} + M_k$, \ $k \in \NN$.
\ Hence
 \begin{align*}
  \left( \EE(\tV_k^{2j}) \right)^{\frac{1}{2j}}
  \leq \beta \left( \EE(\tV_{k-1}^{2j}) \right)^{\frac{1}{2j}}
       + \left( \EE(M_k^{2j}) \right)^{\frac{1}{2j}}
  = \left(\OO((k-1)^j)\right)^{\frac{1}{2j}} + \left(\OO(k^j)\right)^{\frac{1}{2j}}
  = \OO(k^{1/2}) ,
 \end{align*}
 where the first inequality follows by the triangular inequality for the
 $L_{2j}$-norm, and the second one by the induction hypothesis and that
 \ $\EE(M_k^{2j}) = \OO(k^j)$.
\ Hence \ $\EE(\tV^{2j}_k) = \OO(k^{j})$, \ $k \in \NN$, \ as desired.
\proofend

\begin{Cor}\label{LEM_UV_UNIFORM}
Let \ $(X_k)_{k\geq-1}$ \ be an \INARtwo\ process with autoregressive parameters
 \ $(\alpha, \beta) \in [0,1]^2$ \ such that \ $\alpha + \beta = 1$
 \ (hence it is unstable).
Suppose that \ $X_0 = X_{-1} = 0$ \ and \ $\EE(\vare_1^\ell) < \infty$ \ with
 some \ $\ell \in \NN$.
\ Then
 \begin{itemize}
  \item[\textup{(i)}]
   for all \ $i, j \in \ZZ_+$ \ with \ $\max\{i, j\} \leq \ell/2$, \ and for
    all \ $\kappa > i + \frac{j}{2} + 1$, \ we have
    \begin{align}\label{seged_UV_UNIFORM1}
     n^{-\kappa} \sum_{k=1}^n |U_k^i V_k^j| \stoch 0 \qquad
     \text{as \ $n \to \infty$,}
    \end{align}
  \item[\textup{(ii)}]
   for all \ $i, j \in \ZZ_+$ \ with \ $\max\{i, j\} \leq \ell$, \ for all
    \ $T > 0$, \ and for all \ $\kappa > i + \frac{j}{2} + \frac{i+j}{\ell}$,
    \ we have
    \begin{align}\label{seged_UV_UNIFORM2}
     n^{-\kappa} \sup_{t\in[0,T]} |U_\nt^i V_\nt^j| \stoch 0 \qquad
     \text{as \ $n \to \infty$,}
    \end{align}
  \item[\textup{(iii)}]
   for all \ $i, j \in \ZZ_+$ \ with \ $\max\{i, j\} \leq \ell/4$, \ for all
    \ $T > 0$, \ and \ for all \ $\kappa > i + \frac{j}{2} + \frac{1}{2}$,
    \ we have
    \begin{align}\label{seged_UV_UNIFORM4}
     n^{-\kappa} \sup_{t\in[0,T]}
     \left| \sum_{k=1}^\nt [U_k^i V_k^j - \EE(U_k^i V_k^j \mid \cF_{k-1})] \right|
     \stoch 0
     \qquad \text{as \ $n \to \infty$.}
    \end{align}
 \end{itemize}
\end{Cor}

\noindent
\textbf{Proof.}
By Cauchy-Schwartz's inequality and Lemma \ref{EEX_EEU_EEV}, we have
 \begin{align*}
  \EE\left(\sum_{k=1}^n |U_k^i V_k^j|\right)
  \leq \sum_{k=1}^n \sqrt{\EE(U_k^{2i}) \EE(V_k^{2j})}
  = \sum_{k=1}^n \sqrt{\OO(k^{2i}) \OO(k^{j})}
  = \sum_{k=1}^n \OO(k^{i+j/2})
  = \OO(n^{1+i+j/2}) .
 \end{align*}
Using Slutsky's lemma this implies \eqref{seged_UV_UNIFORM1}.

Now we turn to prove \eqref{seged_UV_UNIFORM2}.
First note that
 \begin{align}\label{seged_UV_UNIFORM3}
  \sup_{t\in[0,T]} |U_\nt^i V_\nt^j|
  \leq \sup_{t\in[0,T]} |U_\nt^i| \sup_{t\in[0,T]} |V_\nt^j| ,
 \end{align}
 and for all \ $\vare > 0$ \ and \ $\delta > 0$, \ we have, by Markov's
 inequality,
 \begin{align*}
  \PP\left( n^{-\vare} \sup_{t\in[0,T]} |U_\nt^i| > \delta \right)
  &= \PP\left( n^{-\ell\vare/i} \sup_{t\in[0,T]} |U_\nt^\ell| > \delta^{\ell/i} \right)
   \leq \sum_{k=1}^\nT \PP(U_k^\ell > \delta^{\ell/i} n^{\ell\vare/i}) \\
  &\leq \sum_{k=1}^\nT \frac{\EE(U_k^\ell)}{\delta^{\ell/i} n^{\ell\vare/i}}
   =\sum_{k=1}^\nT \frac{\OO(k^\ell)}{\delta^{\ell/i} n^{\ell\vare/i}}
   = \OO(n^{\ell+1-\ell\vare/i}),
   \qquad i \in \{1, 2, \ldots, \ell\} ,
 \end{align*}
 and
 \begin{align*}
  &\PP\left( n^{-\vare} \sup_{t\in[0,T]} |V_\nt^j| >\delta \right)
   = \PP\left( n^{-\ell\vare/j} \sup_{t\in[0,T]} |V_\nt^\ell| > \delta^{\ell/j} \right)
   \leq \sum_{k=1}^\nT \PP(|V_k^\ell| > \delta^{\ell/j} n^{\ell\vare/j}) \\
  &\qquad
   \leq \sum_{k=1}^\nT \frac{\EE(|V_k^\ell|)}{\delta^{\ell/j} n^{\ell\vare/j}}
   \leq \sum_{k=1}^\nT \frac{\sqrt{\EE(V_k^{2\ell})}}{\delta^{\ell/j} n^{\ell\vare/j}}
   = \sum_{k=1}^\nT \frac{\OO(k^{\ell/2})}{\delta^{\ell/j} n^{\ell\vare/j}}
   = \OO(n^{\ell/2+1-\ell\vare/j}) ,
   \quad j \in \{1, 2, \ldots, \ell\} .
 \end{align*}
Hence, if \ $\ell + 1 - \ell\vare/i < 0$, \ i.e.,
 \ $\vare > \frac{\ell+1}{\ell}i$, \ then
 \[
   n^{-\vare} \sup_{t\in[0,T]} |U_\nt^i| \stoch 0
   \qquad \text{as \ $n \to \infty$,}
 \]
 and if \ $\ell/2 + 1 - \ell\vare/j < 0$, \ i.e.,
 \ $\vare > \frac{\ell/2+1}{\ell}j$, \ then
 \[
   n^{-\vare} \sup_{t\in[0,T]} |V_\nt^j| \stoch 0
   \qquad \text{as \ $n \to \infty$.}
 \]
By \eqref{seged_UV_UNIFORM3}, we get \eqref{seged_UV_UNIFORM2}.

Finally, we show \eqref{seged_UV_UNIFORM4}.
Applying Doob's maximal inequality (see, e.g., Revuz and Yor
 \cite[Chapter II, Theorem 1.7]{RevYor}) for the martingale
 \[
   \sum_{k=1}^n \big[U_k^i V_k^j - \EE(U_k^i V_k^j \mid \cF_{k-1}) \big] ,
   \qquad n \in \NN ,
  \]
 (with the filtration \ $(\cF_k)_{k\in\NN}$) \ and then \eqref{cond_var}, we
 obtain
 \begin{align*}
  \EE\left(\sup_{t\in[0,T]}
            \Biggl( \sum_{k=1}^\nt
                     \big[U_k^i V_k^j
                          - \EE(U_k^i V_k^j \mid \cF_{k-1}) \big]
            \Biggr)^2 \right)
  &\leq 4 \EE\left( \Biggl( \sum_{k=1}^\nT
                            \big[U_k^i V_k^j
                            - \EE(U_k^i V_k^j \mid \cF_{k-1}) \big]
                   \Biggr)^2 \right)\\
  &\leq 4 \sum_{k=1}^\nT \EE(U_k^{2i} V_k^{2j})
   = \sum_{k=1}^\nT \OO(k^{2i+j}) = \OO(n^{2i+j+1}) ,
 \end{align*}
 since
 \ $\EE(U_k^{2i} V_k^{2j}) \leq \sqrt{\EE(U_k^{4i})\EE(V_k^{4j})} = \OO(k^{2i+j})$
 \ by Corollary \ref{EEX_EEU_EEV}.
\proofend

\begin{Rem}
We note that in the special case \ $(\ell, i, j) = (2, 1, 0)$, \ we also get
 \begin{align}\label{supU_stronger}
  n^{-\kappa} \sup_{t\in[0,T]} U_{\nt} \stoch 0 \qquad
  \text{as \ $n \to \infty$ \ for \ $\kappa > 1$.}
 \end{align}
Indeed, by \eqref{rec_U}, we have
 \begin{equation}\label{U}
  U_n = \sum_{k=1}^n (M_k + \mu) , \qquad n \in \NN ,
 \end{equation}
 and hence convergence \eqref{supU_stronger} will follow from
 \begin{align}\label{supsumM}
  n^{-\kappa} \sup_{t \in [0,T]} \left| \sum_{k=1}^{\nt} M_k \right| \stoch 0 \qquad
  \text{as \ $n \to \infty$ \ for all \ $\kappa > 1$.}
 \end{align}
Doob's maximal inequality
 (see, e.g., Revuz and Yor \cite[Chapter II, Theorem 1.7]{RevYor})
 for the martingale \ $\sum_{i=1}^k M_i$, \ $k \in \NN$,
 \ (with the filtration \ $(\cF_k)_{k\in\NN}$) \ gives
 \begin{align*}
  \EE\left(\sup_{t\in[0,T]} \Biggl( \sum_{k=1}^\nt M_k \Biggr)^2 \right)
  \leq 4 \EE\left( \Biggl( \sum_{k=1}^\nT M_k \Biggr)^2 \right)
  = 4 \sum_{k=1}^\nT \EE(M_k^2)
  = \OO(n^2) ,
 \end{align*}
 since  \ $\EE(M_k^2) = \OO(k)$ \ by Corollary \ref{EEX_EEU_EEV}.
This implies \eqref{supsumM}, hence \eqref{supU_stronger}.

However, it turns out that we do not need this stronger statement.
\proofend
\end{Rem}

%\eject

\appendix

%\appendixpage

\noindent{\bf\Large Appendicies}

\section{Classification of \INARtwo\ processes}
\label{app_A}

An \INARtwo\ process is called \emph{positively regular} if there is a
 positive integer \ $k$ \ such that the entries of \ $A^k$ \ are positive
 (see Kesten and Stigum \cite{KesSti1}).
If \ $\alpha > 0$ \ and \ $\beta > 0$ \ then the \INARtwo\ process is
 positively regular, since
 \[
   A = \begin{bmatrix} \alpha & \beta \\ 1 & 0 \end{bmatrix} , \qquad
   A^2 = \begin{bmatrix}
          \alpha^2 + \beta & \alpha \beta \\
          \alpha & \beta
         \end{bmatrix} .
 \]
If \ $\alpha = 0$, \ then
 \[
   A^{2k+1} = \beta^k A
           = \beta^k \begin{bmatrix} 0 & \beta \\ 1 & 0 \end{bmatrix} , \qquad
   A^{2k} = \beta^k \begin{bmatrix} 1 & 0 \\ 0 & 1 \end{bmatrix} , \qquad
   k \in \ZZ_+ ,
 \]
 hence the process is not positively regular.
If \ $\beta = 0$, \ then
 \[
   A^k = \alpha^{k-1} A
       = \alpha^{k-1} \begin{bmatrix} \alpha & 0 \\ 1 & 0 \end{bmatrix} ,
   \qquad k \in \NN ,
 \]
 hence the process is not positively regular.
Consequently, an \INARtwo\ process is positively regular if and only if
 \ $\alpha > 0$ \ and \ $\beta > 0$.

An \INARtwo\ process is called \emph{decomposable} if the matrix \ $A$ \ is
 decomposable (see Kesten and Stigum \cite{KesSti3}).
Note that an \INARtwo\ process is decomposable if and only if the matrix \ $A$
 \ is reducible (see Horn and Johnson \cite[Definition 6.2.21]{HJ}), that is,
 there exists a permutation matrix \ $P \in \RR^{2\times2}$ \ such that
 \[
   P^\top A P = \begin{bmatrix} b & c \\ 0 & d \end{bmatrix},
 \]
 where \ $b, c, d \in \RR$.
\ Since
 \[
   \begin{bmatrix}
     0 & 1 \\
     1 & 0 \\
   \end{bmatrix}
   \begin{bmatrix}
     \alpha & \beta \\
     1 & 0 \\
   \end{bmatrix}
   \begin{bmatrix}
     0 & 1 \\
     1 & 0 \\
   \end{bmatrix}
   = \begin{bmatrix}
     0 & 1 \\
     \beta & \alpha \\
   \end{bmatrix},
 \]
 we get an \INARtwo\ process is decomposable if and only if \ $\beta = 0$.
\ Moreover, an \INARtwo\ process is indecomposable but not positively regular
 if and only if \ $\alpha = 0$ \ and \ $\beta > 0$.

Note that an \INARtwo\ process is positively regular if and only if the matrix
 \ $A$ \ is primitive
 (see Horn and Johnson \cite[Definition 8.5.0 and Theorem 8.5.2]{HJ}), so this
 case can also be called \emph{primitive}
 (see Barczy et al.\ \cite[Definition 2.4]{BarIspPap0}).
Further we remark that the not positively regular case is also called
 \emph{non-primitive}.

\section{A version of the continuous mapping theorem}
\label{app_B}

A function \ $f : \RR_+ \to \RR^d$ \ is called \emph{c\`adl\`ag} if it is right
 continuous with left limits.
\ Let \ $\DD(\RR_+, \RR^d)$ \ and \ $\CC(\RR_+, \RR^d)$ \ denote the space of
 all $\RR^d$-valued c\`adl\`ag and continuous functions on \ $\RR_+$,
 \ respectively.
Let \ $\cB(\DD(\RR_+, \RR^d))$ \ denote the Borel $\sigma$-algebra on
 \ $\DD(\RR_+, \RR^d)$ \ for the metric defined in Jacod and Shiryaev
 \cite[Chapter VI, (1.26)]{JSh} (with this metric \ $\DD(\RR_+, \RR^d)$ \ is a
 complete and separable metric space and the topology induced by this metric
 is the so-called Skorokhod topology).
For $\RR^d$-valued stochastic processes \ $(\bcY_t)_{t \in \RR_+}$ \ and
 \ $(\bcY^{(n)}_t)_{t \in \RR_+}$, \ $n \in \NN$, \ with c\`adl\`ag paths we write
 \ $\bcY^{(n)} \distr \bcY$ \ if the distribution of \ $\bcY^{(n)}$ \ on the
 space \ $(\DD(\RR_+, \RR), \cB(\DD(\RR_+, \RR^d)))$ \ converges weakly to the
 distribution of \ $\bcY$ \ on the space
 \ $(\DD(\RR_+, \RR), \cB(\DD(\RR_+, \RR^d)))$ \ as \ $n \to \infty$.
\ Concerning the notation \ $\distr$ \ we note that if \ $\xi$ \ and
 \ $\xi_n$, \ $n \in \NN$, \ are random elements with values in a metric space
 \ $(E, d)$, \ then we also denote by \ $\xi_n \distr \xi$ \ the weak
 convergence of the distributions of \ $\xi_n$ \ on the space \ $(E, \cB(E))$
 \ towards the distribution of \ $\xi$ \ on the space \ $(E, \cB(E))$ \ as
 \ $n \to \infty$, \ where \ $\cB(E)$ \ denotes the Borel $\sigma$-algebra on
 \ $E$ \ induced by the given metric \ $d$.

The following version of continuous mapping theorem can be found for example
 in Kallenberg \cite[Theorem 3.27]{K}.

\begin{Lem}\label{Lem_Kallenberg}
Let \ $(S, d_S)$ \ and \ $(T, d_T)$ \ be metric spaces and
 \ $(\xi_n)_{n \in \NN}$, \ $\xi$ \ be random elements with values in \ $S$
 \ such that \ $\xi_n \distr \xi$ \ as \ $n \to \infty$.
\ Let \ $f : S \to T$ \ and \ $f_n : S \to T$, \ $n \in \NN$, \ be measurable
 mappings and \ $C \in \cB(S)$ \ such that \ $\PP(\xi \in C) = 1$ \ and
 \ $\lim_{n \to \infty} d_T(f_n(s_n), f(s)) = 0$ \ if
 \ $\lim_{n \to \infty} d_S(s_n,s) = 0$ \ and \ $s \in C$.
\ Then \ $f_n(\xi_n) \distr f(\xi)$ \ as \ $n \to \infty$.
\end{Lem}

For the case \ $S := \DD(\RR_+, \RR^d)$ \ and \ $T := \RR^q$
 \ ($T := \DD(\RR_+,\RR^q)$), \ where \ $d, q \in \NN$, \ we formulate a
 consequence of Lemma \ref{Lem_Kallenberg}.

For functions \ $f$ \ and \ $f_n$, \ $n \in \NN$, \ in \ $\DD(\RR_+, \RR^d)$,
 \ we write \ $f_n \lu f$ \ if \ $(f_n)_{n\in\NN}$ \ converges to \ $f$
 \ locally uniformly, i.e., if \ $\sup_{t\in[0,T]} \|f_n(t) - f(t)\| \to 0$ \ as
 \ $n \to \infty$ \ for all \ $T > 0$.
\ For measurable mappings
 \ $\Phi : \DD(\RR_+, \RR^d) \to \RR^q \ (\DD(\RR_+, \RR^q))$ \ and
 \ $\Phi_n : \DD(\RR_+, \RR^d) \to \RR^q \ (\DD(\RR_+,\RR^q))$, \ $n \in \NN$,
 \ we will denote by \ $C_{\Phi,(\Phi_n)_{n\in\NN}}$ \ the set of all functions
 \ $f \in \CC(\RR_+, \RR^d)$ \ such that
 \ $\Phi_n(f_n) \to \Phi(f) \ (\lu \Phi(f))$ \ whenever
 \ $f_n \lu f$ \ with \ $f_n \in \DD(\RR_+, \RR^d)$, \ $n \in \NN$.

We will use the following version of the continuous mapping theorem several
 times, see, e.g., Barczy et al. \cite[Lemma 4.2]{BarIspPap1} and Isp\'any and
 Pap \cite[Lemma 3.1]{IspPap}.

\begin{Lem}\label{Conv2Funct}
Let \ $d, q \in \NN$, \ and \ $(\bcU_t)_{t\in\RR_+}$ \ and
 \ $(\bcU^{(n)}_t)_{t\in\RR_+}$, \ $n \in \NN$, \ be $\RR^d$-valued stochastic
 processes with c\`adl\`ag paths such that \ $\bcU^{(n)} \distr \bcU$.
\ Let \ $\Phi : \DD(\RR_+, \RR^d) \to \RR^q \ (\DD(\RR_+, \RR^q))$ \ and
 \ $\Phi_n : \DD(\RR_+, \RR^d) \to \RR^q \ (\DD(\RR_+,\RR^q))$, \ $n \in \NN$,
 \ be measurable mappings such that there exists
 \ $C \subset C_{\Phi,(\Phi_n)_{n\in\NN}}$ \ with \ $C \in \cB(\DD(\RR_+, \RR^d))$
 \ and \ $\PP(\bcU \in C) = 1$.
\ Then \ $\Phi_n(\bcU^{(n)}) \distr \Phi(\bcU)$.
\end{Lem}

In order to apply Lemma \ref{Conv2Funct}, we will use the following statement
 several times.

\begin{Lem}\label{Marci}
Let \ $d, p, q \in \NN$, \ $h : \RR^d \to \RR^q$ \ be a continuous function and
 \ $K : [0,1] \times \RR^{2d} \to \RR^p$ \ be a function such that for all
 \ $R > 0$ \ there exists \ $C_R > 0$ such that
 \begin{equation}\label{Lipschitz}
  \| K(s, x) - K(t, y) \| \leq C_R \left( | t - s | + \| x - y \| \right)
 \end{equation}
 for all \ $s, t \in [0, 1]$ \ and \ $x, y \in \RR^{2d}$ \ with
 \ $\| x \| \leq R$ \ and \ $\| y \| \leq R$.
\ Moreover, let us define the mappings
 \ $\Phi, \Phi_n : \DD(\RR_+, \RR^d) \to \RR^{q+p}$, \ $n \in \NN$, \ by
 \begin{align*}
  \Phi_n(f)
  &:= \left( h(f(1)),
             \frac{1}{n}
             \sum_{k=1}^n
              K\left( \frac{k}{n}, f\left( \frac{k}{n} \right),
                      f\left( \frac{k-1}{n} \right) \right) \right) , \\
  \Phi(f)
  & := \left(  h(f(1)), \int_0^1 K( u, f(u), f(u) ) \, \dd u \right)
 \end{align*}
 for all \ $f \in \DD(\RR_+, \RR^d)$.
\ Then the mappings \ $\Phi$ \ and \ $\Phi_n$, \ $n \in \NN$, \ are
 measurable, and
 \ $C_{\Phi,(\Phi_n)_{n \in \NN}} = \CC(\RR_+, \RR^d) \in \cB(\DD(\RR_+, \RR^d))$.
\end{Lem}

\noindent{\bf Proof.}
For an arbitrary Borel set \ $B \in \cB(\RR^{q + p})$ \ we have
 \[
   \Phi_n^{-1}(B) = \pi_{0,\frac{1}{n},\frac{2}{n},\dots,1}^{-1}(\tK_n^{-1}(B)) ,
   \qquad n \in \NN ,
 \]
 where for all \ $n \in \NN$ \ the mapping
 \ $\tK_n : (\RR^{d})^{n+1} \to \RR^{q + p}$ \ is defined by
 \[
   \tK_n(x_0,x_1, \dots, x_n)
   := \left( h(x_n),
             \frac{1}{n}
             \sum_{k=1}^n K\left( \frac{k}{n}, x_k, x_{k-1} \right) \right) ,
   \qquad x_0, x_1, \dots, x_n \in \RR^d ,
 \]
 and the natural projections
\ $\pi_{t_0,t_1,t_2,\dots,t_n} : \DD(\RR_+, \RR^d) \to (\RR^d)^{n+1}$,
 \ $t_0, t_1, t_2, \dots, t_n \in \RR_+$, \ are given by
 \ $\pi_{t_0,t_1,t_2,\dots,t_n}(f) := (f(t_0), f(t_1), f(t_2), \dots, f(t_n))$,
 \ $f \in \DD(\RR_+, \RR^d)$, \ $t_0, t_1, t_2, \dots, t_n \in \RR_+$.
\ Since \ $h$ \ and \ $K$ \ are continuous, \ $\tK_n$ \ is also continuous,
 and hence \ $\tK_n^{-1}(B) \in \cB((\RR^{d})^{n+1})$.
\ It is known that \ $\pi_{t_0,t_1,t_2,\dots,t_n}$,
 \ $t_0, t_1, t_2, \dots, t_n \in \RR_+$, \ are measurable mappings
 (see, e.g., Billingsley \cite[Theorem 16.6 (ii)]{Bil} or Ethier and Kurtz
 \cite[Proposition 3.7.1]{EthKur}),
 and hence \ $\Phi_n = \tK_n \circ \pi_{0,\frac{1}{n},\frac{2}{n},\dots,1}$ \ is also
 measurable.

Next we show the measurability of \ $\Phi$.
\ Since the natural projection
 \ $\DD(\RR_+, \RR^d) \ni f\mapsto f(1) = \pi_1(f)$ \ is measurable, \ $h$ \ is
 continuous, it is enough to show that the mapping
 \[
   \DD(\RR_+, \RR^d) \ni f
   \mapsto \tPhi(f) := \int_0^1 K( t, f(t) , f(t) ) \, \dd t
 \]
 is measurable.
Namely, we show that \ $\tPhi$ \ is continuous.
We have to check that \ $\tPhi(f_n) \to \tPhi(f)$ \ in \ $\RR^{p}$ \ as
 \ $n \to \infty$ \ whenever \ $f_n \to f$ \ in \ $\DD(\RR_+, \RR^d)$ \ as
 \ $n \to \infty$, \ where \ $f, f_n \in D(\RR_+,\RR^d)$, \ $n \in \NN$.
\ Due to Ethier and Kurtz \cite[Proposition 3.5.3]{EthKur}, for all \ $T > 0$
 \ there exists a sequence \ $\lambda_n : \RR_+ \to \RR_+$, \ $n \in \NN$,
 \ of strictly increasing continuous functions with \ $\lambda_n(0) = 0$ \ and
 \ $\lim_{t \to \infty} \lambda_n(t) = \infty$ \ such that
 \begin{align}\label{01seged15}
  \lim_{n \to \infty} \sup_{t\in[0, T]} | \lambda_n(t) - t | = 0 , \qquad
  \lim_{n \to \infty} \sup_{t\in[0, T]} \| f_n(t) - f(\lambda_n(t)) \| = 0 .
 \end{align}
We check that \ $\lim_{n\to\infty} f_n(t) = f(t)$ \ whenever \ $t \in \RR_+$ \ is
 a continuity point of \ $f$.
\ This readily follows by
 \[
   \|f_n(t)-f(t)\|
   \leq \|f_n(t)-f(\lambda_n(t))\| + \|f(\lambda_n(t))-f(t)\| ,
   \qquad n \in \NN , \quad t \in \RR_+ .
 \]
Using that \ $f$ \ has at most countably many discontinuities
 (see, e.g., Jacod and Shiryaev \cite[page 326]{JSh}), we have
 \ $\lim_{n\to\infty} f_n(t) = f(t)$ \ for all \ $t \in \RR_+$ \ except a
 countable set having Lebesgue measure zero.
In what follows we check that
 \[
   \sup_{n \in \NN} \sup_{t \in [0, 1]} \| K( t, f_n(t), f_n(t)) \| < \infty .
 \]
Since \ $K$ \ is continuous and hence it is bounded on a compact set, it is
 enough to verify that
 \[
  \sup_{n\in\NN} \sup_{t\in[0,1]} \| f_n(t) \| < \infty .
 \]
This follows by Jacod and Shiryaev \cite[Chapter VI, Lemma 1.14 (b)]{JSh},
 since \ $f_n \to f$ \ in \ $D(\RR_+, \RR^d)$ \ yields that
 \ $\{f_n : n \in \NN\}$ \ is a relatively compact set
 (with respect to the Skorokhod topology).
Then Lebesgue dominated convergence theorem yields the continuity of
 \ $\tPhi$.

In order to show \ $C_{\Phi,(\Phi_n)_{n\in\NN}} = \CC(\RR_+, \RR^d)$ \ we have to
 check that \ $\Phi_n(f_n) \to \Phi(f)$ \ whenever \ $f_n \lu f$ \ with
 \ $f \in \CC(\RR_+, \RR^d)$ \ and \ $f_n \in \DD(\RR_+, \RR^d)$,
 \ $n \in \NN$.
\ We have
 \begin{align*}
  \| \Phi_n(f_n) - \Phi(f) \|
  &\leq \| h(f_n(1)) - h(f(1)) \| \\
  &\quad
        + \frac{1}{n}
          \sum_{k = 1}^n
           \left\| K \left( \frac{k}{n}, f_n\left( \frac{k}{n} \right),
                            f_n\left( \frac{k-1}{n} \right) \right)
                   - K \left( \frac{k}{n}, f\left( \frac{k}{n} \right),
                              f\left( \frac{k-1}{n} \right) \right)
           \right\| \\
  &\quad + \sum_{k = 1}^n
            \int_{(k-1)/n}^{k/n}
             \left\| K \left( \frac{k}{n}, f\left( \frac{k}{n} \right) ,
                              f \left( \frac{k-1}{n} \right) \right)
                     - K( t, f(t), f(t) ) \right\| \, \dd t \\
  &=: \| h(f_n(1)) - h(f(1)) \| + A_n^{(1)} + A_n^{(2)} .
 \end{align*}
Since \ $f_n \lu f$ \ implies that \ $f_n(1) \to f(1)$ \ as \ $n \to \infty$,
 \ using the continuity of \ $h$, \ we get
 \[
   \| h(f_n(1)) - h(f(1)) \| \to 0
   \qquad \text{as \ $n \to \infty$.}
 \]
Let us also observe that
 \[
   \sup_{n\in\NN} \sup_{t\in[0,1]} \| f_n(t) \|
   \leq \sup_{n\in\NN} \sup_{t\in[0,1]}
         \| f_n(t) - f(t) \| + \sup_{t\in[0,1]} \| f(t) \|
   =: c < \infty ,
 \]
 hence
 \begin{align*}
  \left\| \left( f_n\left( \frac{k}{n} \right),
                 f_n\left( \frac{k-1}{n} \right) \right) \right\|
  \leq \sqrt{2} c, \qquad n \in \NN , \quad k \in \{1, \ldots, n\} ,
 \end{align*}
 and then, by \eqref{Lipschitz},
 \[
   A_n^{(1)} \leq \sqrt{2}C_{\sqrt{2}c} \sup_{t\in[0,1]} \| f_n(t) - f(t) \| \to 0
 \]
 as \ $n \to \infty$.
\ Moreover,
 \begin{align*}
  A_n^{(2)}
  &\leq C_{\sqrt{2}c}
        \sum_{k = 1}^n
         \int_{(k-1)/n}^{k/n}
          \left( \left| \frac{k}{n} - t \right|
                 + \left\| \left( f\left( \frac{k}{n} \right),
                                  f\left( \frac{k-1}{n} \right) \right)
                           -  ( f(t), f(t) ) \right\|
          \right) \dd t \\
   &\leq \sqrt{2} C_{\sqrt{2}c} ( n^{-1} + \omega_1(f, n^{-1}) ) ,
 \end{align*}
 where
 \[
   \omega_1(f, \vare)
   := \sup_{t, \, s \in [0,1], \, |t-s|<\vare}
       \|\ f(t) - f(s) \| , \qquad \vare > 0 ,
 \]
 denotes the modulus of continuity of \ $f$ \ on \ $[0, 1]$.
\ Since \ $f$ \ is continuous, \ $\omega_1(f, n^{-1}) \to 0$ \ as
 \ $n \to \infty$
 \ (see, e.g., Jacod and Shiryaev \cite[Chapter VI, 1.6]{JSh}),
 and we obtain \ $A_n^{(2)} \to 0$ \ as \ $n \to \infty$.
\ Then \ $C_{\Phi,(\Phi_n)_{n\in\NN}} = \CC(\RR_+, \RR^d)$.

Finally, \ $\CC(\RR_+, \RR^d) \in \cB(\DD(\RR_+, \RR^d))$ \ holds since
 \ $\DD(\RR_+, \RR^d) \setminus \CC(\RR_+, \RR^d)$ \ is open.
Indeed, if \ $f \in \DD(\RR_+, \RR^d) \setminus \CC(\RR_+, \RR^d)$ \ then there
 exists \ $t \in \RR_+$ \ such that
 \ $\vare := \| f(t) - \lim_{s \uparrow t} f(s) \| > 0$, \ and then the open ball
 in \ $\DD(\RR_+, \RR^d)$ \ with centre \ $f$ \ and radius \ $\vare / 2$ \ does
 not contain any continuous function.
We note that for \ $\CC(\RR_+, \RR^d) \in \cB(\DD(\RR_+, \RR^d))$ \ one can
 also simply refer to Ethier and Kurtz \cite[Problem 3.11.25]{EthKur}.
\proofend

\section{Convergence of random step processes}
\label{section_conv_step_processes}

We recall a result about convergence of random step processes towards a
 diffusion process, see Isp\'any and Pap \cite{IspPap}.
This result is used for the proof of convergence \eqref{conv_Z}.

\begin{Thm}\label{Conv2DiffThm}
Let \ $\bgamma : \RR_+ \times \RR^d \to \RR^{d \times r}$ \ be a continuous
 function.
Assume that uniqueness in the sense of probability law holds for the SDE
 \begin{equation}\label{SDE}
  \dd \, \bcU_t
  = \gamma (t, \bcU_t) \, \dd \bcW_t ,
  \qquad t \in \RR_+,
 \end{equation}
 with initial value \ $\bcU_0 = \bu_0$ \ for all \ $\bu_0 \in \RR^d$, \ where
 \ $(\bcW_t)_{t \in \RR_+}$ \ is an $r$-dimensional standard Wiener process.
Let \ $(\bcU_t)_{t \in \RR_+}$ \ be a solution of \eqref{SDE} with initial value
 \ $\bcU_0 = \bzero \in \RR^d$.

For each \ $n \in \NN$, \ let \ $(\bU^{(n)}_k)_{k\in\NN}$ \ be a sequence of
 $d$-dimensional martingale differences with respect to a filtration
 \ $(\cF^{(n)}_k)_{k\in\ZZ_+}$, i.e., \ $\EE(\bU^{(n)}_k \mid \cF^{(n)}_{k-1}) = \bzero$,
 \ $k, n \in \NN$.
\ Let
 \[
   \bcU^{(n)}_t := \sum_{k=1}^{\nt} \bU^{(n)}_k ,
   \qquad t \in \RR_+, \quad n \in \NN .
 \]
Suppose that \ $\EE \big( \|\bU^{(n)}_k\|^2 \big) < \infty$ \ for all
 \ $k, n \in \NN$.
\ Suppose that for each \ $T > 0$,
 \begin{enumerate}
  \item [\textup{(i)}]
        $\sup\limits_{t\in[0,T]}
         \left\| \sum\limits_{k=1}^{\nt}
                  \EE\Bigl(\bU^{(n)}_k (\bU^{(n)}_k)^\top \mid \cF^{(n)}_{k-1}\Bigr)
                 - \int_0^t
                    \bgamma(s,\bcU^{(n)}_s) \bgamma(s,\bcU^{(n)}_s)^\top
                    \dd s \right\|
         \stoch 0$,\\
  \item [\textup{(ii)}]
        $\sum\limits_{k=1}^{\lfloor nT \rfloor}
          \EE \big( \|\bU^{(n)}_k\|^2 \bone_{\{\|\bU^{(n)}_k\| > \theta\}}
                    \bmid \cF^{(n)}_{k-1} \big)
         \stoch 0$
        \ for all \ $\theta>0$,
 \end{enumerate}
 where \ $\stoch$ \ denotes convergence in probability.
Then \ $\bcU^{(n)} \distr \bcU$ \ as \ $n \to \infty$.
\end{Thm}

Note that in (i) of Theorem \ref{Conv2DiffThm}, \ $\|\cdot\|$ \ denotes
 a matrix norm, while in (ii) it denotes a vector norm.

\end{document}